\documentclass[11pt]{article}
\usepackage{fullpage}

\date{}
\usepackage{tikz}
\usetikzlibrary{calc}
\usetikzlibrary{matrix}
\usepackage{amsmath,amssymb,amsfonts,amsthm,subcaption}
\usepackage[pdftex]{hyperref} 
\usepackage{caption}
\usepackage{graphicx}
\allowdisplaybreaks
\usepackage{float}
\usepackage{titlesec}
\usepackage{rotating}

\titleformat{name=\section,numberless}
  {\normalfont\Large\bfseries}
  {}
  {0pt}
  {}

 \def \M{\mathbb{M}}
 \def \J{\mathbb{J}}

\newcommand{\bey}{\begin{eqnarray}}
\newcommand{\eey}{\end{eqnarray}}

\newcommand{\beq}{\begin{equation}}
\newcommand{\eeq}{\end{equation}}
\theoremstyle{plain}

\theoremstyle{definition}

\theoremstyle{remark}
\newtheorem{exam}{\hspace{1mm}Example}[section]

\title{Selection of the Regularization Parameter in the Ambrosio-Tortorelli Approximation of the Mumford-Shah
Functional for Image Segmentation}
\author{Yufei Yu%
\thanks{Department of Mathematics, the University of Kansas, Lawrence, Kansas 66045, U.S.A.
({\it y920y782@ku.edu})},
\and Weizhang Huang%
\thanks{Department of Mathematics, the University of Kansas, Lawrence, Kansas 66045, U.S.A.
({\it whuang@ku.edu})}
}
\begin{document}
\maketitle

\begin{abstract}
The Ambrosio-Tortorelli functional is a phase-field approximation of the Mumford-Shah functional that has
been widely used for image segmentation. The approximation has the advantages of being easy to implement,
maintaining the segmentation ability, and $\Gamma$-converging to the Mumford-Shah functional.
However, it has been observed in actual computation that the segmentation
ability of the Ambrosio-Tortorelli functional varies significantly with different values of the parameter
and it even fails to $\Gamma$-converge to the original functional for some cases.
In this paper we present an asymptotic analysis on the gradient flow equation
of the Ambrosio-Tortorelli functional and show that
the functional can have different segmentation behavior for small but finite values of the regularization parameter
and eventually loses its segmentation ability as the parameter goes to zero when the input image is treated as a continuous function.
This is consistent with the existing observation as well as the numerical examples presented in this work.
A selection strategy for the regularization parameter and a scaling procedure
for the solution are devised based on the analysis. Numerical results show that they lead to good segmentation of
the Ambrosio-Tortorelli functional for real images.
 \end{abstract}



{\bf Abbreviated title:} Selection of regularization parameter in Ambrosio-Tortorelli functional

\section{Introduction}

Segmentation for a given image is a process to find the edges of objects and partition the image into separate parts
that are relatively smooth. This can be achieved by minimizing some objective functionals 
and multiple theories have been developed. One of the most commonly used functionals, proposed by Mumford
and Shah \cite{MS89}, takes the form
\begin{equation}
E[u, \Gamma] = \frac{\alpha}{2} \int_{\Omega \setminus \Gamma} | \nabla u |^{2}dx
+ \beta H^{1}(\Gamma) + \frac{\gamma}{2} \int_{\Omega}(u-g)^{2}dx,
\label{MandS}
\end{equation}
where $\Omega$ is a rectangular domain, $\alpha$, $\beta$, and $\gamma$ are positive parameters,
$g$ is the grey level of the input image,
$u$ is the target image, $\Gamma$ denotes the edges of the objects in the image, and $H^{1}(\Gamma)$
is the one-dimensional Hausdorff measure. Upon minimization,
$u$ is close to $g$, $\nabla u$ is small on $\Omega \setminus \Gamma$, and
$\Gamma$ is as short as possible. An optimal image is thus close to the original one and almost piecewise
constant. Moreover, the terms in (\ref{MandS}) represent different and often conflicting objectives,
making its minimization and thus image segmentation an interesting but challenging topic to study.

To avoid mathematical difficulties caused by the $H^1(\Gamma)$ term, De Giorgi et al. \cite{DGCL89}
propose an alternative functional as
\begin{equation}
F[u] =  \frac{\alpha}{2} \int_{\Omega} | \nabla u|^{2}dx + \beta H^{1}(S_{u}) + \frac{\gamma}{2} \int_{\Omega} |u -g|^{2} dx,
\label{DandL}
\end{equation}
where $S_{u}$ is the jump set of $u$. They show that (\ref{DandL}) has minimizers in $SBV(\Omega)$
(the space of special functions of bounded variation)
and is equivalent to (\ref{MandS}) in the sense that
if $u \in {\Omega}$ is a minimizer of (\ref{DandL}), then $(u, S_{u})$ is a minimizer of (\ref{MandS}).

Although it is a perfectly fine functional to study in mathematics, (\ref{DandL}) is not easy to implement
in actual computation due to the fact that the jump set of the unknown function and its
Hausdorff measure are extremely difficult, if not impossible, to compute.
To avoid this difficulty,  Ambrosio and Tortorelli \cite{AT92} propose a regularized version as
\begin{equation}
AT_{\epsilon}[u,\phi] =  \frac{\alpha}{2} \int_{\Omega}(\phi^{2} + k_{\epsilon}) |\nabla u| ^{2} dx
+ \beta \int_{\Omega}\left (\epsilon |\nabla \phi|^{2} + \frac{1}{4 \epsilon} (1 - \phi)^{2}\right )dx
+ \frac{\gamma}{2} \int_{\Omega}|u - g|^{2} dx,
\label{AandT}
\end{equation}
where $\epsilon > 0$ is the regularization parameter,
$k_{\epsilon} = o(\epsilon)$ is a parameter used to prevent the functional from becoming degenerate,
and $\phi$ is a new unknown variable which ideally is an approximation of the characteristic function
for the jump set of $u$, i.e., 
\begin{equation}
\phi(x)  \approx \chi_u (x) \equiv \begin{cases}
0, \quad  & \text{ if }  x \in S_u\\
1, \quad  &\text{ if } x \notin S_u .
\end{cases}
\label{phidefine}
\end{equation}
They show that $AT_{\epsilon}$ has minimizers $u \in SBV(\Omega)$ and $\phi \in L^2(\Omega)$ and
$\Gamma$-converges to $F(u)$. $\Gamma$-convergence, first introduced by
de Giorgi and Franzoni \cite{Giorgi1975},
is a concept that guarantees the minimizer of $AT_{\epsilon}$ converges to that of $F$ as $\epsilon \to 0$.

The first finite element approximation for the functional $AT_{\epsilon}$ is given by
Bellettini~and~Coscia~\cite{Bellettini1994}. They seek linear finite element approximations
$u_{h}$ and $\phi_{h}$ to minimize 
\begin{align}
AT_{\epsilon,h}[u_{h},\phi_{h}]  & =   \frac{\alpha}{2} \int_{\Omega}(\phi_{h}^{2} + k_{\epsilon})|\nabla u_{h}|^{2} dx
+ \beta \int_{\Omega} \left (\epsilon|\nabla \phi_{h}|^2 + \frac{1}{4 \epsilon} \pi_{h}((1-\phi_{h})^{2})\right )dx
\notag \\
& + \frac{\gamma}{2} \int_{\Omega} \pi_{h}((u_{h}-g_{\epsilon})^2)dx,
\label{GB}
\end{align}
where $\pi_{h}$ is the linear Lagrange interpolation operator and $g_{\epsilon}$ is a smooth function which converges
to $g$ in the $L^{2}$ norm as $\epsilon \to 0$. They show that $AT_{\epsilon,h}$ $\Gamma$-converges to $F(u)$
when the maximum element diameter is chosen as $h = o(\epsilon)$.
It should be pointed out that Feng and Prohl \cite{Feng2004} have established the existence and
uniqueness of the solution to an initial-boundary value problem (IBVP) of the gradient flow equation
of (\ref{AandT}) and proven that a finite element approximation of the IBVP converges to the continuous
solution as the mesh is refined. 

It is noted that the Ambrosio-Tortorelli functional (\ref{AandT}) is actually a phase-field approximation of
the Mumford-Shah functional (\ref{MandS}). Phase-field modeling has been used widely
in science and engineering to handle sharp interfaces, boundaries, and cracks in numerical simulation
of problems such as dendritic crystal growth \cite{Kobayashi-1993,Wheeler-1993},
multiple-fluid hydrodynamics \cite{LiuShen2003,Shen2014,Shen2015c,LiuShen2006},
and brittle fracture \cite{Bourdin-2000,Francfort-1998,Miehe-2010}.
It employs a phase-field variable $\phi$, which depends on a regularization parameter $\epsilon$
describing the actual width of the smeared interfaces, to indicate the location of the interfaces.
Phase-field modeling has the advantage of being able to handle complex interfaces without
relying on their explicit description.
Mathematically, phase-field models such as (\ref{AandT}) have been studied extensively (e.g., see \cite{AT92})
for $\Gamma$-convergence.
However, few studies have been published for the role of the regularization parameter
in actual simulation. It is a common practice that a specific value of $\epsilon$ is used
without discussion or explanation in phase-field modeling. Even worse, it has been observed
\cite{May2015,Pham2011,Vignollet-2014}
that a phase-field model for brittle fracture simulation does not $\Gamma$-converge as $\epsilon \to 0$
and $\epsilon$ can be interpreted as a material parameter since its choice influences the critical stress.
More recently, the choice of $\epsilon$ has been studied in \cite{Nguyen} based on physical arguments
and with experimental validation.

The objective of this paper is to study the effects and selection of the regularization parameter in
the Ambrosio-Tortorelli functional (\ref{AandT}) which
is a special example of phase-field modeling in image segmentation. We consider the gradient flow equation
of $AT_{\epsilon}$ subject to a homogeneous Neumann boundary condition and carry out
an asymptotic analysis for the solution of the corresponding IBVP as $\epsilon \to 0$. We show that,
when $g$ is continuous, the functional can have different segmentation behavior for small but finite $\epsilon$
and eventually loses its segmentation ability for infinitesimal $\epsilon$. This is consistent with the existing
observation in phase-field modeling and with the numerical examples to be presented. The analysis is also
used to devise a selection strategy for $\epsilon$ and a scaling for $u$ and $g$. Numerical results
with real images confirm that the strategies can lead to good segmentation of $AT_{\epsilon}$
in the sense that $\phi$ is close to the characteristic function of $g$ (cf. (\ref{phidefine})).

An outline of the paper is as follows. The asymptotic analysis is given in Section~\ref{SEC:analysis},
followed by the description of a moving mesh finite element method in Section~\ref{SEC:fem}. Illustrative numerical examples
are given in Section~\ref{SEC:numerics}. A selection strategy for $\epsilon$ and a scaling procedure
for $u$ and $g$ as well as examples with several real input images are presented in Section~\ref{SEC:select}. Finally,
Section~\ref{SEC:conclusion} contains conclusions.

\section{Behavior of the minimizer of $AT_{\epsilon}$ as $\epsilon \to 0$ for continuous $g$}
\label{SEC:analysis}

We first explain why we consider $g$ as a continuous function.
In image segmentation,
the function $g$ represents an image and is given the grey-level values at the pixels.
Generally speaking, the values of $g$ at points other than the pixels are needed in finite
element computation. These values are computed commonly
through (linear) interpolation of the values at the pixels.
This means that $g$ is treated as a continuous function in finite element computation
and such a treatment is independent of the regularization parameter in the phase-field modeling.
Thus we consider $g$ as a continuous function and study
the behavior of the minimizer of $AT_{\epsilon}$ as $\epsilon \to 0$ in this section.

To this end, we consider the gradient flow equation of functional (\ref{AandT}), 
\begin{equation}
\begin{cases}
u_{t}= \alpha \nabla \cdot ((k_{\epsilon} + \phi^{2})\nabla u) - \gamma (u - g),
\\
\phi_{t} = 2\beta \epsilon \Delta \phi - \alpha | \nabla u|^{2}\phi + \frac{\beta}{2 \epsilon}(1 - \phi),
\end{cases}
\quad  x \in \Omega, \; t > 0
 \label{grad-flow-1}
\end{equation}
subject to the homogeneous Neumann boundary condition
\begin{equation}
\frac{\partial u}{\partial n} = \frac{\partial \phi}{\partial n} = 0 \quad \text{for} \quad x \in \partial \Omega
\label{bc1}
\end{equation}
and the initial condition
\begin{equation}
u(x,0) = u^0(x), \quad \phi(x,0) = \phi^0(x), \quad x \in \Omega.
\label{ic}
\end{equation}
This IBVP has been studied and used to find the minimizer of (\ref{AandT})
(as a steady-state solution) by a number of researchers. Noticeably, Feng and Prohl \cite{Feng2004} have
established the existence and uniqueness of the solution of the IBVP and proven that a finite element approximation
converges to the continuous solution as the mesh is refined. 

By assumption, $g \in C^0(\Omega)$. Then we can expect that the solution $u$ and $\phi$ of the IBVP
is smooth. To see the behavior of the solution as $\epsilon \to 0$, we consider the asymptotic expansion of
$u$ and $\phi$ as 
\begin{align}
& u = u^{(0)} + \epsilon u^{(1)} + \epsilon^{2} u^{(2)} + \cdot \cdot \cdot,
\label{perturb1}
\\
& \phi  = \phi^{(0)} + \epsilon \phi^{(1)} + \epsilon^{2} \phi^{(2)} + \cdot \cdot \cdot ,
\label{perturb2}
\end{align}
where $u^{(0)},\; u^{(1)},\; ...$ and $\phi^{(0)},\; \phi^{(1)},\; ...$ are functions independent of $\epsilon$.
Inserting these into (\ref{grad-flow-1}), we get 
\begin{align}
u^{(0)}_t + \epsilon u^{(1)}_t + O(\epsilon^2) =  & \;
\alpha \nabla \cdot \left [ \left (o(\epsilon) + ( \phi^{(0)} + \epsilon \phi^{(1)} + o(\epsilon))^2\right )
\nabla (u^{(0)} + \epsilon u^{(1)} + O(\epsilon^2)) \right ]
\notag \\
& \qquad - \gamma (u^{(0)} + \epsilon u^{(1)} + O(\epsilon^2) - g),
\label{u-1} \\
\phi^{(0)}_t + \epsilon \phi^{(1)}_t + O(\epsilon^2) =  &\;
2\beta \epsilon(\Delta \phi^{(0)} + \epsilon \Delta \phi^{(1)} + O(\epsilon^2))
\notag \\
& - \alpha \left | \nabla u^{(0)} + \epsilon \nabla u^{(1)} + O(\epsilon^2)\right |^2 (\phi^{(0)}
+ \epsilon \phi^{(1)} + O(\epsilon^{2}))
\notag \\
&  + \frac{\beta}{2 \epsilon} (1 - \phi^{(0)} - \epsilon \phi^{(1)} - O(\epsilon^2)),
\label{phi-1}
\end{align}
where we have used $k_\epsilon = o(\epsilon)$.
Collecting the $O(1)$ terms in (\ref{u-1}), we have
\begin{equation}
u_{t}^{(0)}= \alpha\Delta u^{(0)} - \gamma(u^{(0)} - g), \quad \text{ in } \Omega  .
\label{u-3}
\end{equation}
Similarly, collecting the $O(1/\epsilon)$ terms and $O(1)$ terms in (\ref{phi-1}) we get
\begin{align*}
\frac{\beta}{2}(1 - \phi^{(0)}) = 0,
\qquad \phi^{(0)}_t =  -\alpha |\nabla u^{(0)}|^{2} \phi^{(0)} - \frac{\beta}{2} \phi^{(1)} .
\end{align*}
From these we obtain
\begin{align}
\phi = 1 - \epsilon \frac{2 \alpha}{\beta}|\nabla u^{(0)}|^{2} + O(\epsilon^2). 
\label{phi-0}
\end{align}
Like $u$,  $u^{(0)}$ also satisfies a homogeneous Neumann boundary condition.
Since $g \in C^0(\Omega)$, it can be shown (e.g., see \cite{Eva98}) that $\nabla u^{(0)}$ is continuous and bounded.
Combining this with (\ref{phi-0}) we conclude that $\phi \to 1$ as $\epsilon \to 0$. 
Since the boundaries between different objects in $u$ are indicated by $\phi = 0$,
this implies that $u$ is a single object and there is no segmentation as $\epsilon \to 0$
when $g$ is continuous. Moreover, $u$ and thus $u^{(0)}$ are kept close to $g$
and we can expect $\nabla u^{(0)}$ to be large in the places where $\nabla g$ is large.
From (\ref{phi-0}) we can see that, for small but not infinitesimal $\epsilon$,
$\phi$ can become zero at places where $\nabla g$ is large. In this case, the functional will have
good segmentation (cf. the numerical examples in Section~\ref{SEC:numerics}).

The above analysis shows that, when $g$ is continuous, the choice of the regularization parameter in (\ref{AandT})
can be crucial for image segmentation:  {\em different values of $\epsilon$ can lead to
very different segmentation behavior of the functional and its segmentation ability will
disappear as $\epsilon \to 0$}.

It should be emphasized that the above observation is not in contradiction with the theoretical analysis made
in \cite{AT92} for the $\Gamma$-convergence and segmentation ability of the functional (\ref{AandT}).
In \cite{AT92}, these properties are analyzed for $u \in SBV(\Omega)$, implicitly implying that
$u$ is discontinuous in general. The above analysis has been made under the assumption
that $g$ and thus $u$ are continuous although they may have large gradient from place to place.

It is instructive to see some transient behavior of the solution to the gradient flow equation.
To simplify, we drop the diffusion term in the second equation in (\ref{grad-flow-1}) and get
\begin{equation}
 \phi_t = -\alpha |\nabla u|^2 \phi + \frac{\beta}{2 \epsilon}(1 - \phi).
 \label{simplephi}
 \end{equation}
It has been proven in \cite{Feng2004} that the solution of (\ref{grad-flow-1}) satisfies $ 0 \le \phi \le 1$.
From this we see that the first term on the right-hand side of (\ref{simplephi}) is nonpositive,
which makes $\phi$ decrease, and the second term is nonnegative, making $\phi$ increase.
These two terms compete and reach an equilibrium state.
Moreover, if $\phi = 1$, we have $\phi_t = - \alpha | \nabla u|^2 \leq 0$,
meaning that as long as $|\nabla u| \neq 0$, the first term decreases $\phi$ until $\phi_t = 0$ is reached.
Similarly, if $\phi = 0$,  we have $\phi_t = \frac{\beta}{2 \epsilon} > 0$, which means $\phi$ increases
until the system reaches its equilibrium. The equilibrium value of $\phi$ can be obtained by setting
the right-hand side of (\ref{simplephi}) to be zero, i.e., 
\begin{equation}
\phi \approx \frac{\beta}{\beta + 2 \epsilon \alpha |\nabla u|^2} .
\label{phi-2}
\end{equation}
Thus, the equilibrium value of $\phi$ is around 1 for smooth regions where $\nabla u$ is small
and around 0 on edges where $\nabla u$ is large.

\section{The adaptive moving mesh finite element method}
\label{SEC:fem}

In this section we describe an adaptive moving mesh finite element method for solving
the gradient flow equation (\ref{grad-flow-1}).
Recall that a crucial requirement for the Ambrosio-Tortorelli approximation (\ref{AandT}) of the Mumford-Shah
functional is that $\epsilon$ must be small. Since the width of object edges is in the same order of $\epsilon$,
the size of the mesh elements around the edges should be in the same order of $\epsilon$ or smaller
for any finite element approximation to be meaningful.
On the other hand, the mesh elements do not have to be that small within each object where $u$ and $\phi$
are smooth. Thus, mesh adaptation is necessary for the efficiency of the finite element computation.
We use here the MMPDE moving mesh method \cite{HRR94a,HR11} that has been specially designed for
time dependent problems.

It should be pointed out that a number of other moving mesh methods have been developed in the past
and there is a vast literature in the area. The interested reader is referred to the books or review articles \cite{Bai94a,Baines-2011,BHR09,HR11,Tan05} and references therein.
Moreover, moving mesh methods have been successfully applied to phase-field models,
e.g., see \cite{DZ08,MR02,SY09,WLT08,YFLS06,YCD08,ZHLZ2017}.

It is remarked that the spatial domain $\Omega$ is typically a rectangular domain
for image segmentation. However, finite element computation is not subject to this restriction.
Moreover, we will consider examples in both one and two dimensions for illustrative purpose in the next section.
For these reasons, we consider $\Omega$ as a general polygonal domain in $d$-dimensions ($d = 1$ and $2$).

\subsection{Finite element discretization}

We now consider the integration of (\ref{grad-flow-1}) up to a finite time $t = T$.
Denote the time instants by
\[
0 = t_0 < t_1 < ... < t_{n_f} = T .
\]
For the moment, we assume that a simplicial mesh for $\Omega$ is given at these time instants, i.e.,
$\mathcal{T}_h^n$, $n = 0, ..., n_f$, which are considered as the deformation from each other
and have the same number of the elements ($N$) and the vertices ($N_v$) and the same connectivity.
Such a mesh is generated using the MMPDE moving mesh strategy to be described in
Section~\ref{SEC:MMPDE}.

For the finite element discretization of (\ref{grad-flow-1}),
the mesh is considered to change linearly between $t_n$ and $t_{n+1}$, i.e., 
\[
x_j(t) = \frac{t - t_n}{t_{n+1} - t_n} x_j^{n+1}+ \frac{t_{n+1} - t}{t_{n+1} - t_n}x_j^n,\quad j = 1,...,N_v,
\quad t \in (t_n, t_{n+1})
\]
where $x_j(t)$, $x_j^n$, and $x_j^{n+1}$ ($j = 1,...,N_v$)
denote the coordinates of the vertices of $\mathcal{T}_h(t)$, $\mathcal{T}_h^n$, and $\mathcal{T}_h^{n+1}$, respectively.
Denote the linear basis function associated with the $j$-th vertex by $\psi_{j}(\cdot,t)$ and let
$V_h(t) = \text{span}\{\psi_1,...,\psi_{N_v}\}$. Then, the weak formulation for the linear finite element
approximation for (\ref{grad-flow-1}) is to find $u_h(\cdot,t)$, $\phi_h(\cdot,t) \in V^h(t)$,
$ 0 < t \le T$ such that
\begin{equation} 
\begin{cases}
\int_{\Omega} \frac{\partial u_h} {\partial t} v dx = -\alpha \int_{\Omega} (k_{\epsilon} + \phi_h^2) \nabla u_h \cdot \nabla v dx
- \gamma \int_{\Omega}(u_h - g)v dx = 0, &\forall v \in V^h(t)\\
\int_{\Omega}\frac{\partial \phi_h}{\partial t} v dx = - 2 \beta \epsilon \int_{\Omega} \nabla \phi_h \cdot \nabla v dx
- \alpha \int_{\Omega} | \nabla u_h |^2 \phi_h v dx + \frac{\beta}{2 \epsilon} \int_{\Omega}( 1 - \phi_{h})v dx,
&\forall v \in V^h(t).
\end{cases}
\label{numericalpde}
\end{equation}
This is almost the same as that for the finite element approximation on a fixed mesh. The main difference lies in
time differentiation. To see this, expressing $u_h$ into
\begin{align}
u_h(x, t) = \sum_{i = 1}^{N_v} u_i(t) \psi_i(x, t)
\label{udis}
\end{align}
and differentiating it with respect to time, we get
\[
\frac{\partial u_h(x, t)}{\partial t} dx = \sum_{ i = 1}^{Nv} \frac{du_i}{dt}\psi_i(x, t)
+ \sum_{i = 1}^{Nv}u_i(t) \frac{\partial \psi_i(x,t)}{\partial t}.
\]
It is known (e.g., see \cite{HR11}) that
\[
\frac{\partial \psi_i}{\partial t} = -\nabla \psi_i \cdot \dot{X}, \quad\text{a.e. in } \Omega
\]
where
\[
\dot{X} = \sum_{i = 1}^{N_v} \dot{x}_i \psi_i(x,t)
\]
and $\dot{x}_i$'s denote the nodal mesh velocities. Combining the above results, we obtain
\[
\frac{\partial u_h}{\partial t} = \sum_{i = 1}^{N_v} \frac{du_i}{dt} \psi_i - \nabla u_h \cdot \dot{X}.
\]
Similarly,
\[
\phi_h(x, t) = \sum_{i = 1}^{Nv} \phi_i(t) \psi_i(x, t) ,\quad
\frac{\partial \phi_h}{\partial t} = \sum_{i = 1}^{N_v} \frac{d \phi_i}{dt} \psi_i - \nabla \phi_h \cdot \dot{X}.
\]
From these we can see that mesh movement introduces an extra convection term.
Inserting these into (\ref{numericalpde}) and taking $v = \psi_j$ successively,  we can rewrite
(\ref{numericalpde}) into an ODE system in the form
\begin{equation}
\begin{cases}
M(X) \dot{U} = F(\dot{X}, X, \Phi, U, X), \\
M(X) \dot{\Phi} = G(\dot{X}, X, \Phi, U, X),
\end{cases}
\label{simode}
\end{equation}
where $M(X)$ is the mass matrix. This system for $U$ and $\Phi$ is integrated from
$t_n$ to $t_{n+1}$ using the fifth-order Radau IIA method (e.g., see Hairer and Wanner \cite{HW96}),
with a variable time step being selected based on a two-step error estimator \cite{Montijano2004}.

\subsection{The MMPDE moving mesh strategy}
\label{SEC:MMPDE}

We now describe the generation of $\mathcal{T}_h^{n+1}$ using the MMPDE moving mesh strategy \cite{HR11}.
For this purpose, we denote the physical mesh by $\mathcal{T}_h = \{ x_1, ..., x_{N_v} \}$,
the reference computational mesh by $\hat{ \mathcal{T}}_{c,h} = \{ \hat{\xi}_1,..., \hat{\xi}_{N_v} \}$
(which is chosen as the very initial physical mesh in our computation), and the computational mesh
$\mathcal{T}_{c,h} = \{ \xi_1, ..., \xi_{N_v}\}$. We assume that all of these meshes have the same
number of elements and vertices and the same connectivity. Then, for any element $K \in \mathcal{T}_h$
there exists a corresponding element $K_c \in \mathcal{T}_{c,h}$. We denote the affine mapping between
$K_c$ and $K$ by $F_K$ and its Jacobian matrix by $F_K'$.

A main idea of the MMPDE moving mesh strategy is to view any adaptive mesh as a uniform one in
the metric specified by a certain tensor.  A metric tensor (denoted by $\M$)
is a symmetric and uniformly positive definite matrix-valued
function defined on $\Omega$. In our computation, we choose $\M$ to be a piecewise constant function
depending on $u_h$ as
\begin{equation}
\M_K = \det(|H_K|)^{-\frac{1}{d+4}}|H_K|, \quad \forall K \in \mathcal{T}_h
\label{M-1}
\end{equation}
where $H_K$ is a recovered Hessian of $u_h$ on element $K$, $|H_K|
= Q \text{diag}(|\lambda_1|,...,|\lambda_d|)Q^T$, assuming that the eigen-decomposition of $H_K$ is
$Q \text{diag}(\lambda_1,...,\lambda_d)Q^T$,
and $\det(|H_K|)$ is the determinant of $|H_K|$.
The recovered Hessian in $K$ is obtained by twice differentiating a local quadratic polynomial
fitting in the least-squares sense to the nodal values of $u_h$ at the neighboring vertices of the element.
The form of (\ref{M-1}) is known \cite{Hua05b} optimal with respect
to the $L^2$ norm of linear interpolation error. With this choice of $\M$, we hope that the mesh elements
are concentrated in the regions of object edges where the curvature of $u$ is large.

The mesh $\mathcal{T}_h$ being uniform in metric $\M$ will mean that the volume of $K$ in $\M$
is proportional to the volume of $K_c$ with the same proportional constant for all $K \in \mathcal{T}_h$
and $K$ measured in $\M$ is similar to $K_c$.
These requirements can be expressed mathematically as the equidistribution and alignment conditions
(e.g., see \cite{HR11}),
\begin{align}
& |K|\sqrt{\det(\M_K)}  = \frac{\sigma_h |K_c|}{|\Omega_c|}, \qquad \forall K \in \mathcal{T}_h
\label{equi}
\\
& \frac{1}{d}\text{tr}\left ( (F_K')^{-1}\M_K^{-1}(F_K')^{-T}\right ) =
\det \left ( (F_K')^{-1}\M_K^{-1}(F_K')^{-T}\right )^{\frac{1}{d}},
\qquad \forall K \in \mathcal{T}_h
\label{align}
\end{align}
where $|K|$ and $|K_c|$ denote the volume of $K$ and $K_c$, respectively,
$d$ is the dimension of $\Omega$, $\text{tr}(\cdot)$ denotes the trace of a matrix, and
\[
|\Omega_c| = \sum_{K_c \in \mathcal{T}_{c,h}} |K_c| ,\qquad
\sigma_h = \sum \limits_{K\in \mathcal{T}_h} |K| \sqrt{\det(\M_K)} .
\]
An energy functional associated with these conditions has been proposed in \cite{Hua01b} as
\begin{align}
I_h(\mathcal{T}_h, \mathcal{T}_{c,h})  = & \; \theta \sum_{K \in \mathcal{T}_h} |K| \sqrt{\det(\M_K)}
\left (\mathbf{\text{tr}}( (F_k')^{-1}\M_K^{-1}(F_K')^{-T}) \right )^{\frac{dp}{2}}
\notag \\
+ & \; (1 - 2\theta)d^{\frac{d p}{2}}   \sum_{K \in \mathcal{T}_h} |K| \sqrt{\det(\M_K)} 
\left (\frac{|K_c|}{|K|\sqrt{\det(\M_K)}}\right )^p  ,
\label{ih}
\end{align}
where $\theta \in (0,0.5]$ and $p>1$ are two dimensionless parameters. In our computation, we take
$\theta = 1/3$ and $p = 3/2$ which are known experimentally to work well for most problems.

Notice that $I_h$ is a function of $\mathcal{T}_h$ and $\mathcal{T}_{c,h}$.
We can take $\mathcal{T}_{c,h}$ as the reference computational mesh
$\hat{\mathcal{T}}_{c,h}$ and minimize $I_h$ with respect to $\mathcal{T}_h$. With the MMPDE strategy,
the minimization is carried out by integrating a modified gradient system of $I_h$,
\begin{equation}
\frac{\partial x_i}{\partial t} = - \frac{P_i}{\tau} \left (\frac{\partial I_h}{\partial x_i}\right )^T, \quad i = 1,..., N_v,
\quad t \in (t_n, t_{n+1}]
\label{xmethod}
\end{equation}
where $\frac{\partial I_h}{\partial x_i}$ is a row vector, $P_i = \det(\M(x_i))^{\frac{p-1}{2}}$ is a positive function
chosen to make (\ref{xmethod}) invariant under the scaling transformation of $\M$,
and $\tau > 0$ is a positive parameter used to adjust the time scale of mesh movement. 
Starting from $\mathcal{T}_h^n$, we can integrate (\ref{xmethod}) (with proper modifications
for the boundary vertices to allow them to slide on the boundary) from $t_n$ to $t_{n+1}$ to obtain $\mathcal{T}_h^{n+1}$.
Special attention may be needed for the computation of the metric tensor that is typically available only at
$\mathcal{T}_h^n$ (the mesh at $t=t_n$).
During the integration of (\ref{xmethod}), the location of the physical vertices changes,
and the values of $\M$ at these vertices should be updated via
interpolation of its values on the vertices of $\mathcal{T}_h^n$.
It is also worth mentioning that the mesh governed by (\ref{xmethod}) is known  \cite{HK2015}
to stay nonsingular if it is nonsingular initially.

To avoid the need of constantly updating the metric tensor $\M$ during the integration of the mesh equation,
we now consider an indirect approach of minimizing $I_h$. In this approach, we choose
$\mathcal{T}_h = \mathcal{T}_h^n$ and minimize $I_h$ with respect to $\mathcal{T}_{c,h}$.
Then the MMPDE for the computational vertices reads as
\begin{equation}
\frac{\partial \xi_i}{\partial t} = - \frac{P_i}{\tau}\left (\frac{\partial I_h}{\partial \xi_i}\right )^T,
\quad i = 1, ..., N_v,  \quad t \in (t_n, t_{n+1}] .
\label{ximethod}
\end{equation}
Starting from $\hat{\mathcal{T}}_{c,h}$, this equation can be integrated from $t_n$ to $t_{n+1}$ to obtain
a new computational mesh $\mathcal{T}_{c,h}^{n+1}$. In our computation, we use
Matlab\textsuperscript \textregistered\, function {\em ode15s},
a Numerical Differentiation Formula based integrator, for this purpose.
Note that $\mathcal{T}_h^n$ and $\M = \M^n$
are fixed during the integration and $\mathcal{T}_h^n$ and $\mathcal{T}_{c,h}^{n+1}$ form a correspondence.
Denote the correspondence by $\Psi_h$, i.e., $\mathcal{T}_h^n = \Psi_h(\mathcal{T}_{c,h}^{n+1})$.
Then, the new physical mesh is defined as 
\begin{equation}
\mathcal{T}_h^{n+1} = \Psi_h( \hat{\mathcal{T}}_c),
\end{equation}
which can be readily computed using linear interpolation.

A benefit of the above $\xi$-formulation is that the derivative $\partial I_h/\partial \xi_i$ in (\ref{ximethod})
can be found analytically using the notion of scalar-by-matrix differentiation \cite{HK2014} and has a
relatively simple, compact matrix form. Using this, we can rewrite (\ref{ximethod}) into
\begin{equation}
\frac{\partial \xi_i}{\partial t} = \frac{P_i}{\tau} \sum_{K \in \omega_i} |K| v_{i_K}^K, \quad i = 1,..., N_v
\label{ximethod-2}
\end{equation}
where $\omega_i$ is the patch of the elements containing $x_i$ as a vertex, the index $i_K$ denotes
the local index of $x_i$ in $K$, and $v_{i_K}^K$ is the local velocity contributed by the element $K$
to the partial derivative $\frac{\partial I_h}{\partial x_i}$. The local velocities on element $K$ are given by
\begin{equation}
\begin{bmatrix}
(v_1^K)^T \\
\vdots \\
(v_d^K)^T
\end{bmatrix} = -E_K^{-1} \frac{\partial G}{\partial \det(\J)}
- \frac{\partial G}{\partial \det(\J)} \frac{\det(\hat{E}_K)}{\det(E_K)} {\hat{E}_K}^{-1},
\quad v_0^{K} = - \sum_{ i = 1}^d v_d^K,
\label{velocity}
\end{equation}
where $E_K = [x_1^K-x_0^K, ..., x_d^K-x_0^K]$ and $\hat{E}_K = [\xi_1^K-\xi_0^K, ..., \xi_d^K-\xi_0^K]$
are the edge matrices of $K$ and $K_c$, respectively, $\J = (F_K)^{-1} = \hat{E}_K {E_K}^{-1}$,
$G = G(\J, \det(\J), \M_K)$ is a function associated with the meshing energy functional, and
$\partial G/\partial \J$ and $\partial G/\partial \det(\J)$ are the partial derivatives of $G$
with respect to the first and second arguments, respectively.
For the meshing energy functional (\ref{ih}), we have
\begin{align*}
& G = \theta \sqrt{\det(\M_K)} (\text{tr}(\J\M_K^{-1}\J^T))^{\frac{dp}{2}}
+ (1 - 2\theta) d^{\frac{dp}{2}} \sqrt{\det(\M_K)}\left (\frac{\det(\J)}{\sqrt{\det(\M_K)}}\right )^p,
\\
& \frac{\partial G}{\partial \J} = dp \theta \sqrt{\det(\M_K)}(\text{tr}(\J\M_K^{-1}\J^T))^{\frac{dp}{2} - 1}\M_K^{-1}\J^T, \\
& \frac{\partial G}{\partial \det(\J)} = p(1 - 2\theta)d^{\frac{dp}{2}} \det(\M_K)^{\frac{1 - p}{2}} \det(\J)^{p-1}.
\end{align*}

\section{Numerical results: behavior of $(u, \phi)$ as $\epsilon \to 0$}
\label{SEC:numerics}

In this section we present numerical results obtained with the moving mesh finite element method
described in the previous section to illustrate the
analysis in Section~\ref{SEC:analysis}. We choose two analytical functions for $g$, with one each
in one dimension and two dimensions, to simulate the grey-level values of images. In particular,
the sharp jumps in $g$ model the object edges in the image. 

\begin{figure}[htb]
\centering
\begin{subfigure}{0.32\textwidth}
\centering
\includegraphics[scale = 0.24]{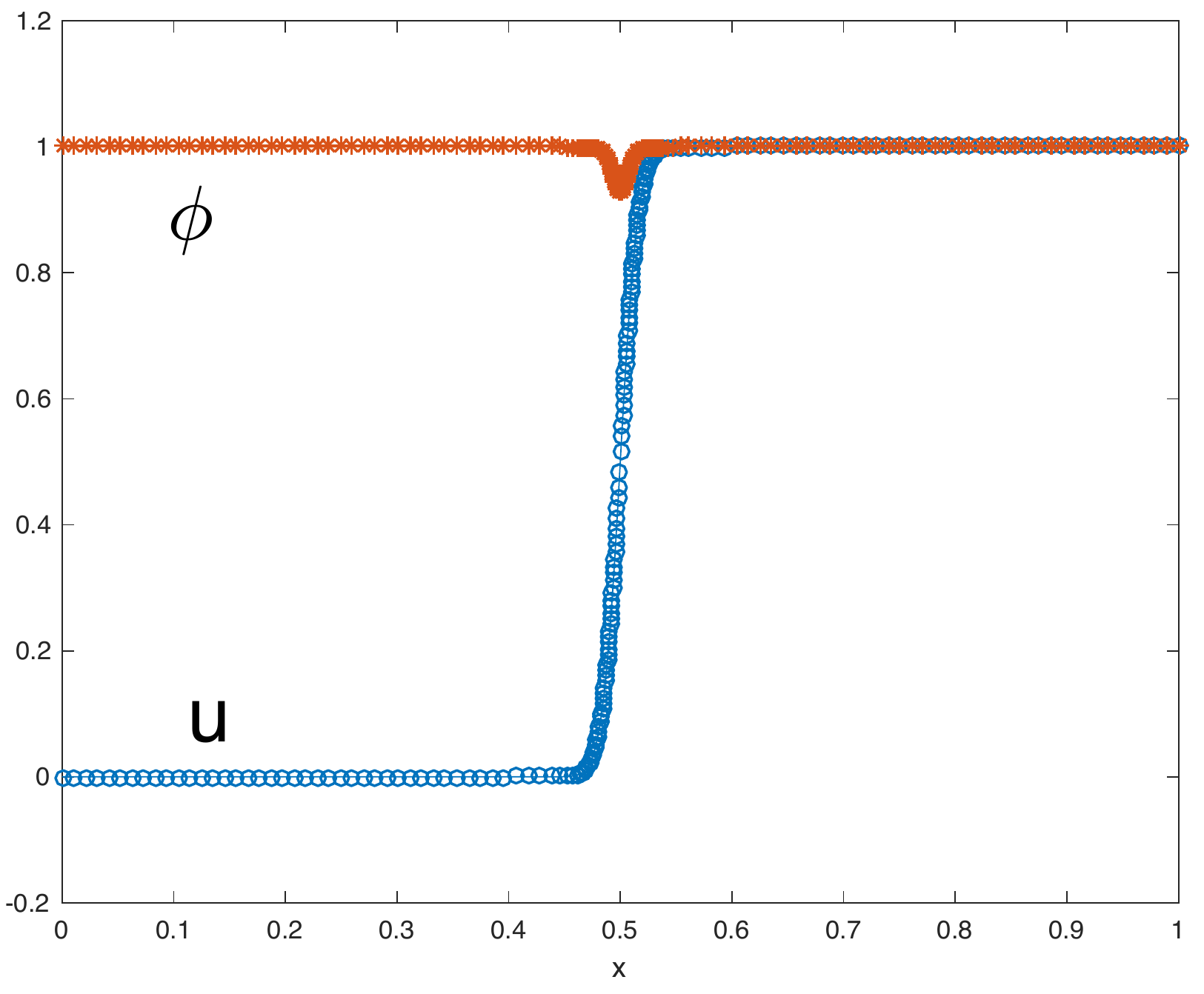}
\caption{ $t = 0.005$, $\epsilon = 0.1$}
\end{subfigure}
\begin{subfigure}{0.32\textwidth}
\centering
\includegraphics[scale = 0.24]{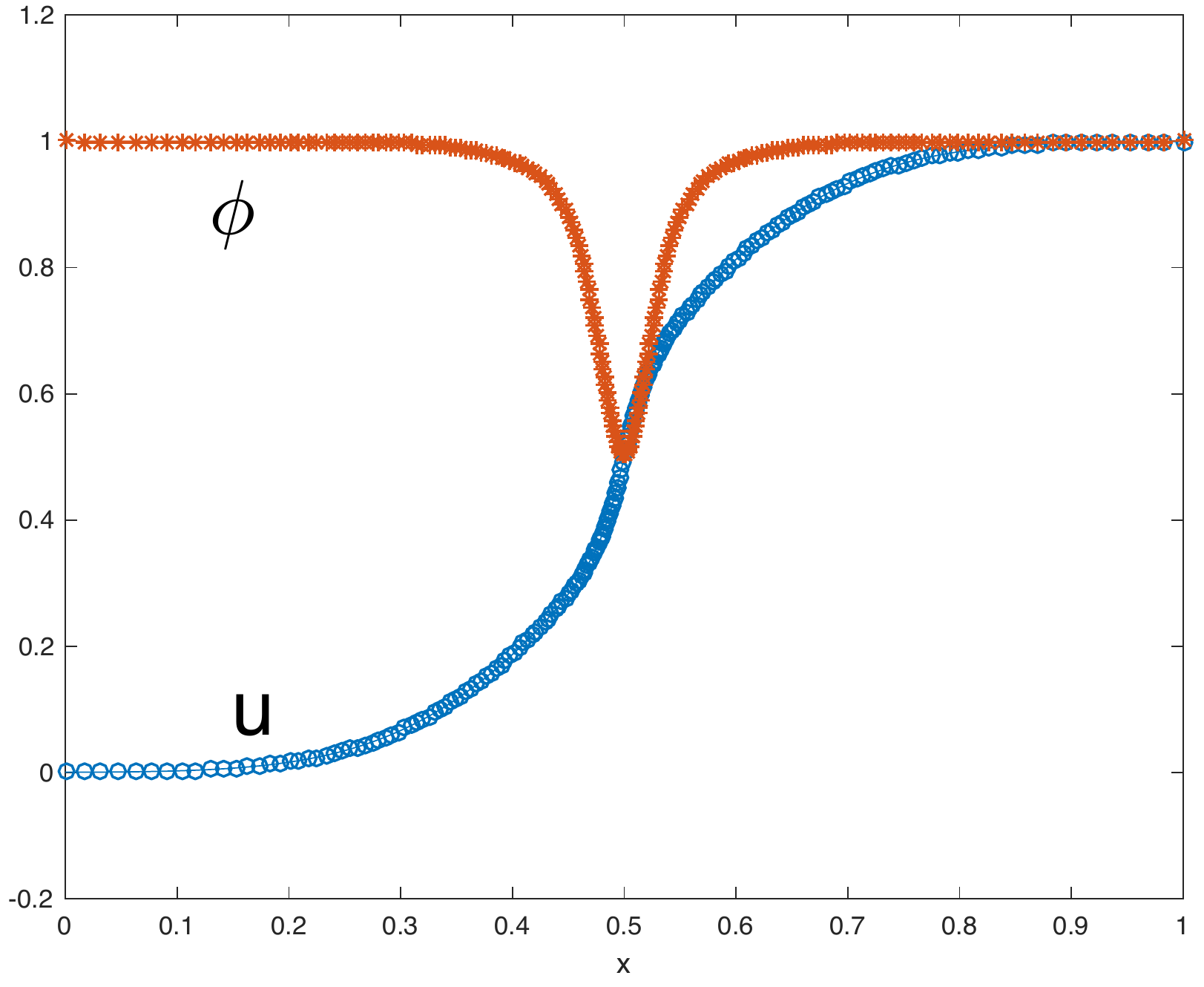}
\caption{ $t = 1$, $\epsilon = 0.1$}
\end{subfigure}
\begin{subfigure}{0.32\textwidth}
\centering
\includegraphics[scale = 0.24]{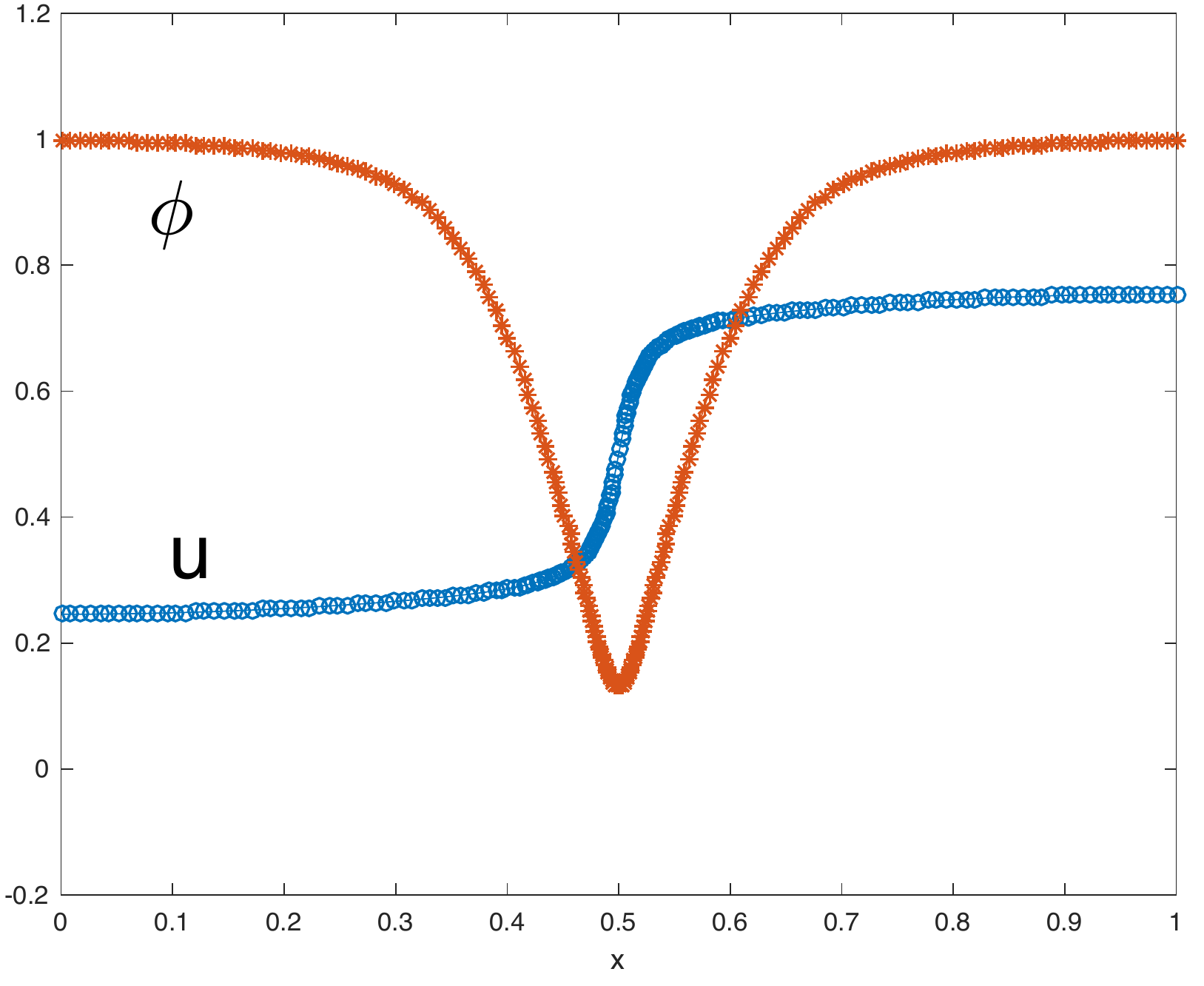}
\caption{ $t = 20$, $\epsilon = 0.1$}
\end{subfigure}
\begin{subfigure}{0.32\textwidth}
\centering
\includegraphics[scale = 0.24]{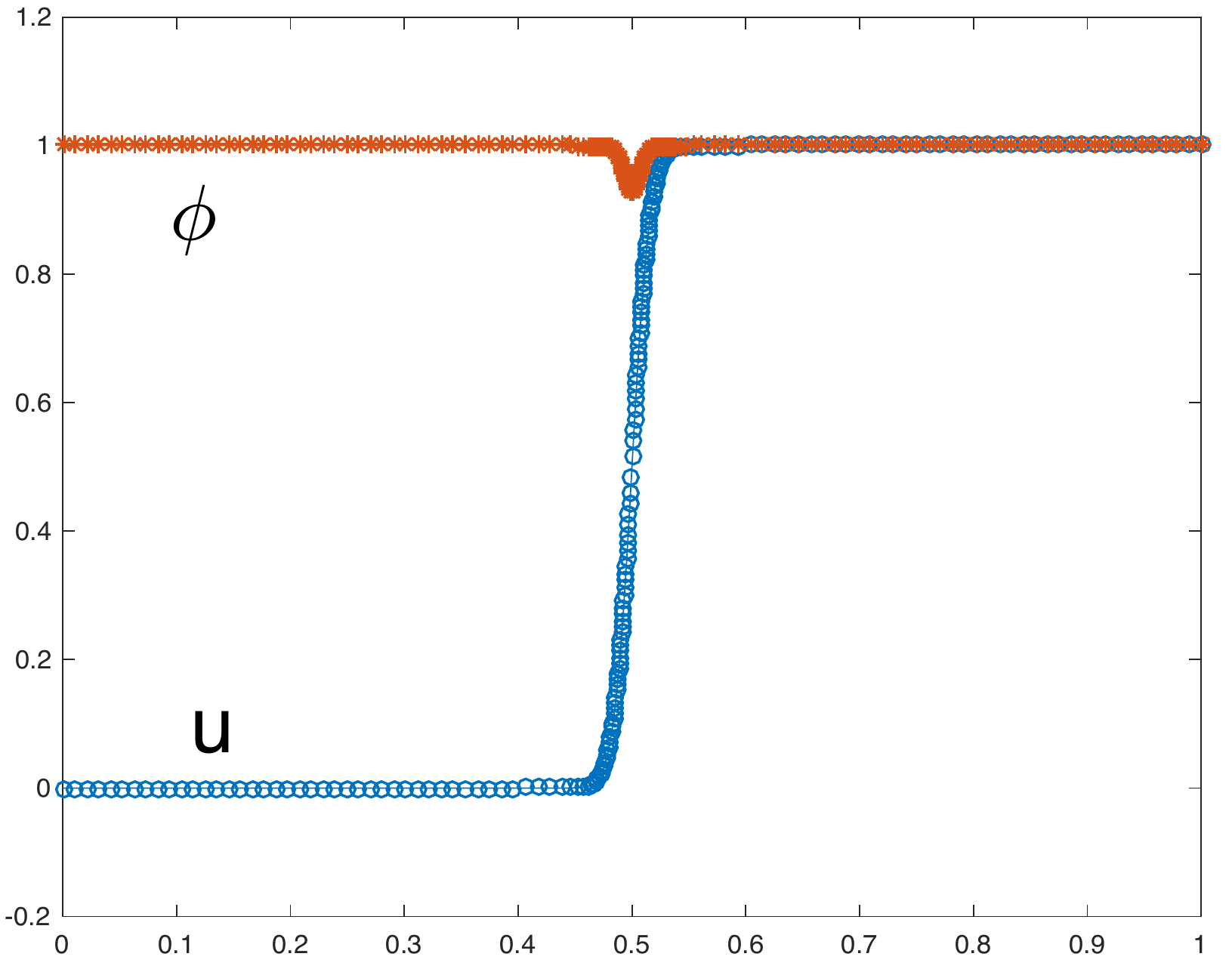}
\caption{ $t = 0.005$, $\epsilon = 0.01$}
\end{subfigure}
\begin{subfigure}{0.32\textwidth}
\centering
\includegraphics[scale = 0.24]{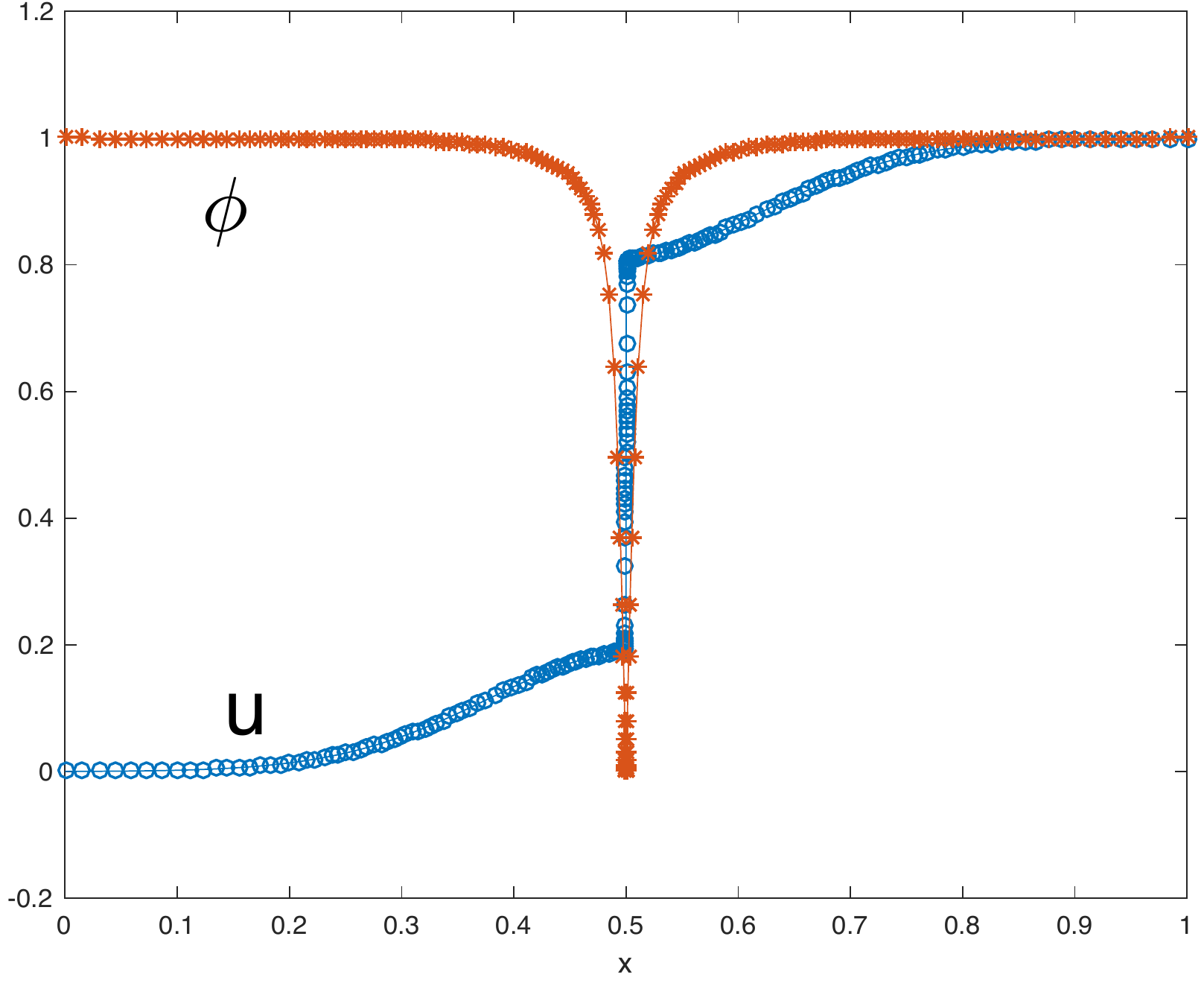}
\caption{ $t = 1$, $\epsilon = 0.01$}
\end{subfigure}
 \begin{subfigure}{0.32\textwidth}
 \centering
\includegraphics[scale = 0.24]{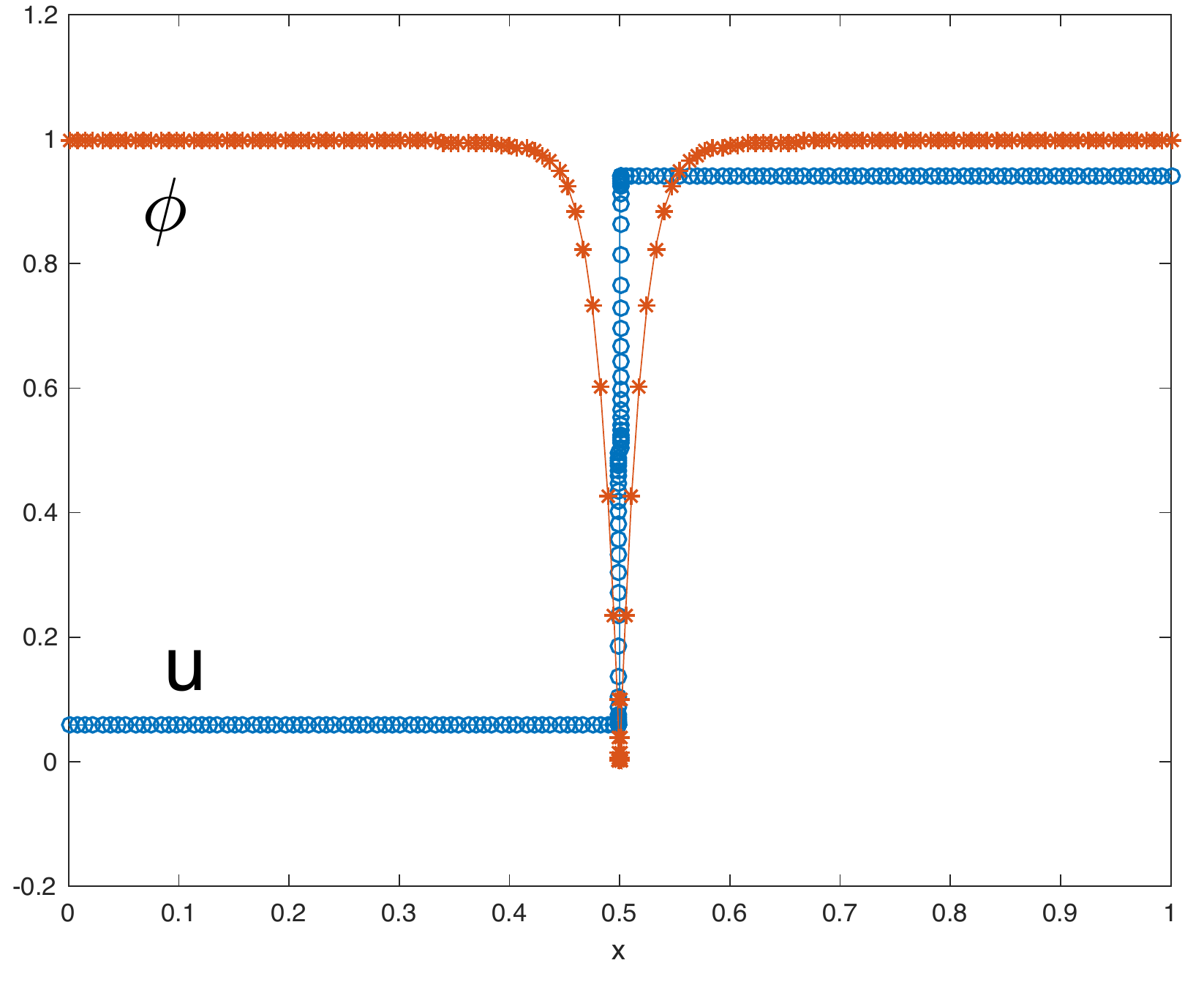}
\caption{ $t = 20$, $\epsilon = 0.01$}
\end{subfigure}
\begin{subfigure}{0.32\textwidth}
\centering
\includegraphics[scale = 0.24]{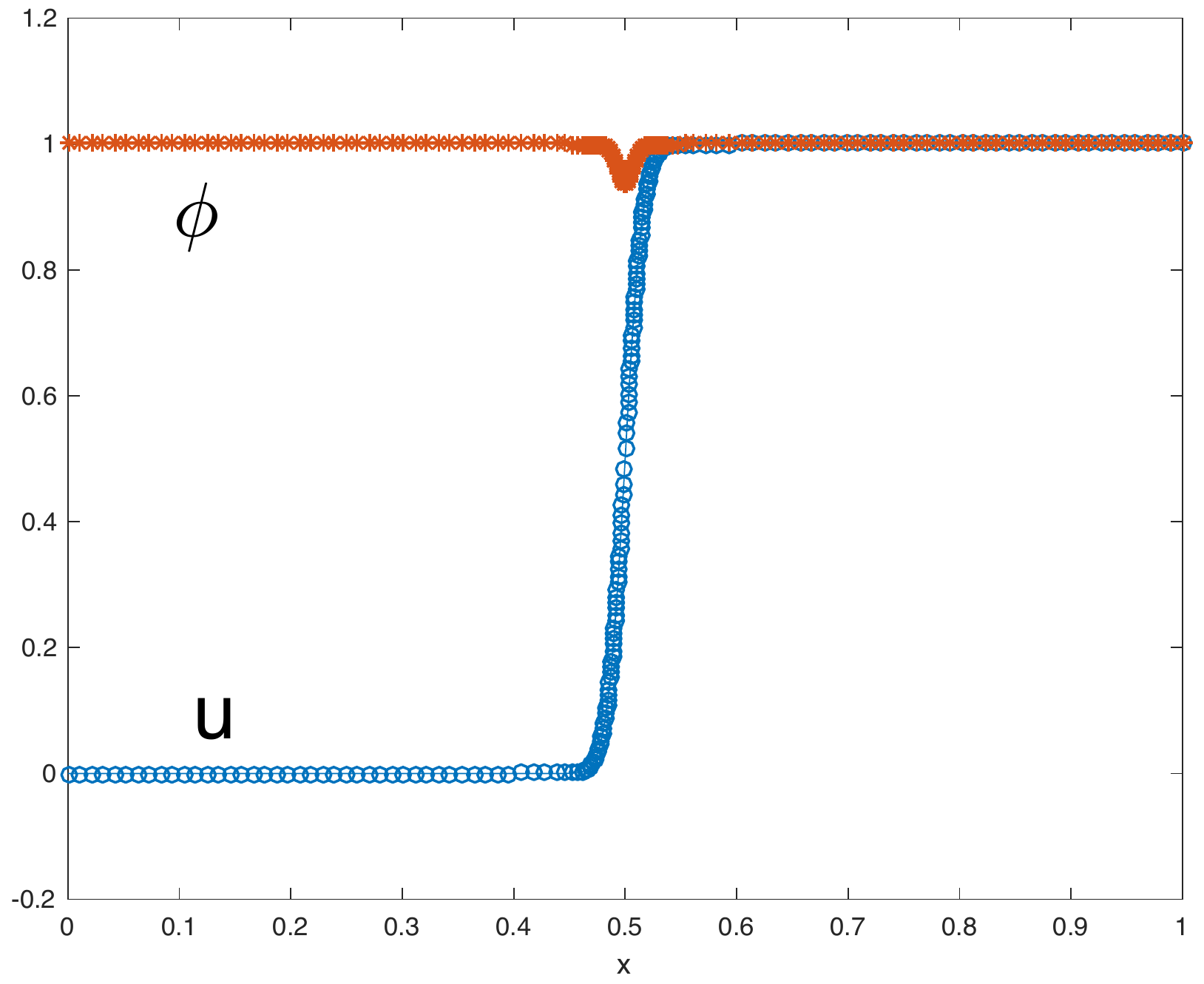}
\caption{ $t = 0.005$, $\epsilon = 10^{-5}$}
\end{subfigure}
\begin{subfigure}{0.32\textwidth}
\centering
\includegraphics[scale = 0.24]{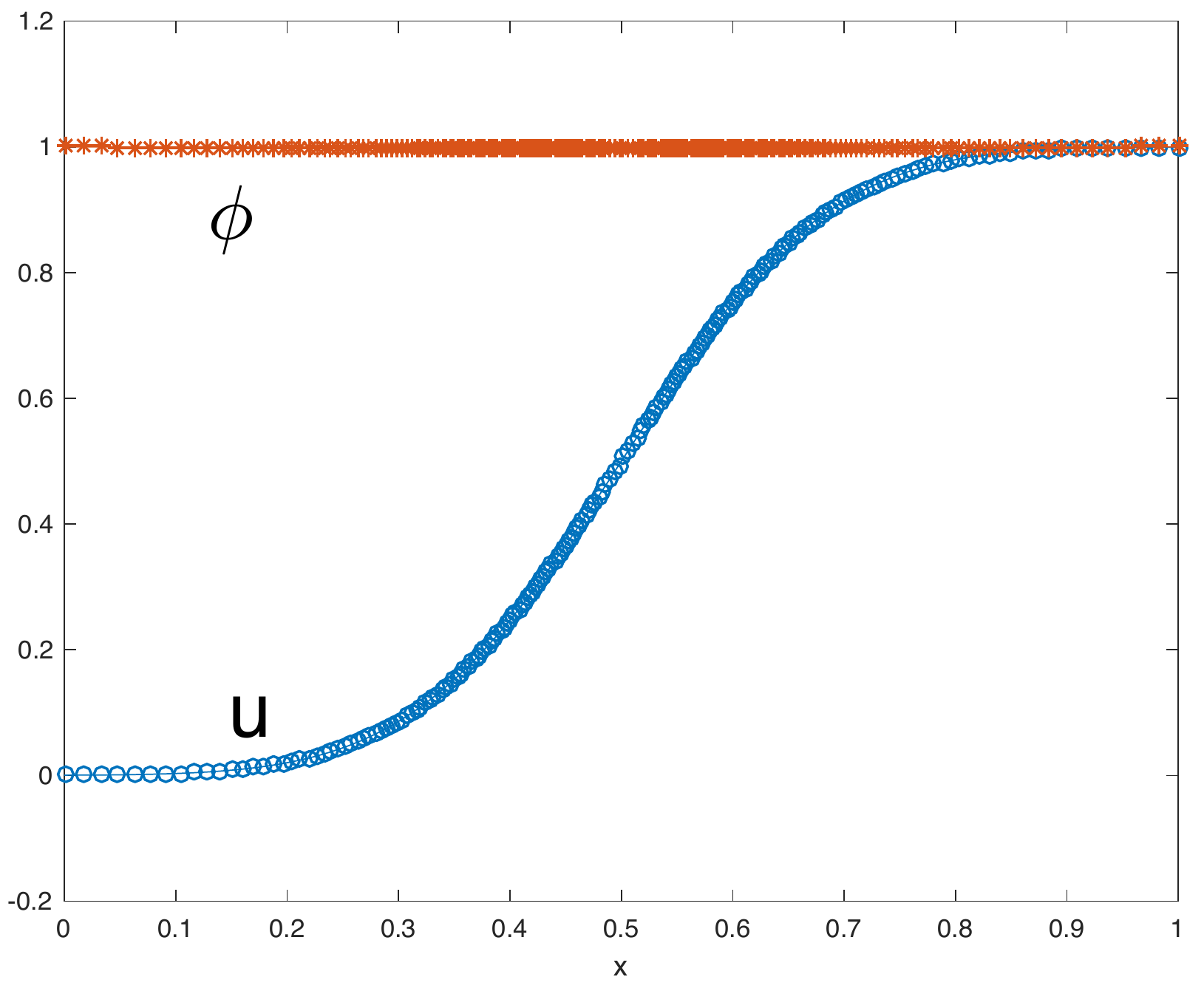}
\caption{ $t = 1$, $\epsilon = 10^{-5}$}
\end{subfigure}
 \begin{subfigure}{0.32\textwidth}
 \centering
\includegraphics[scale = 0.24]{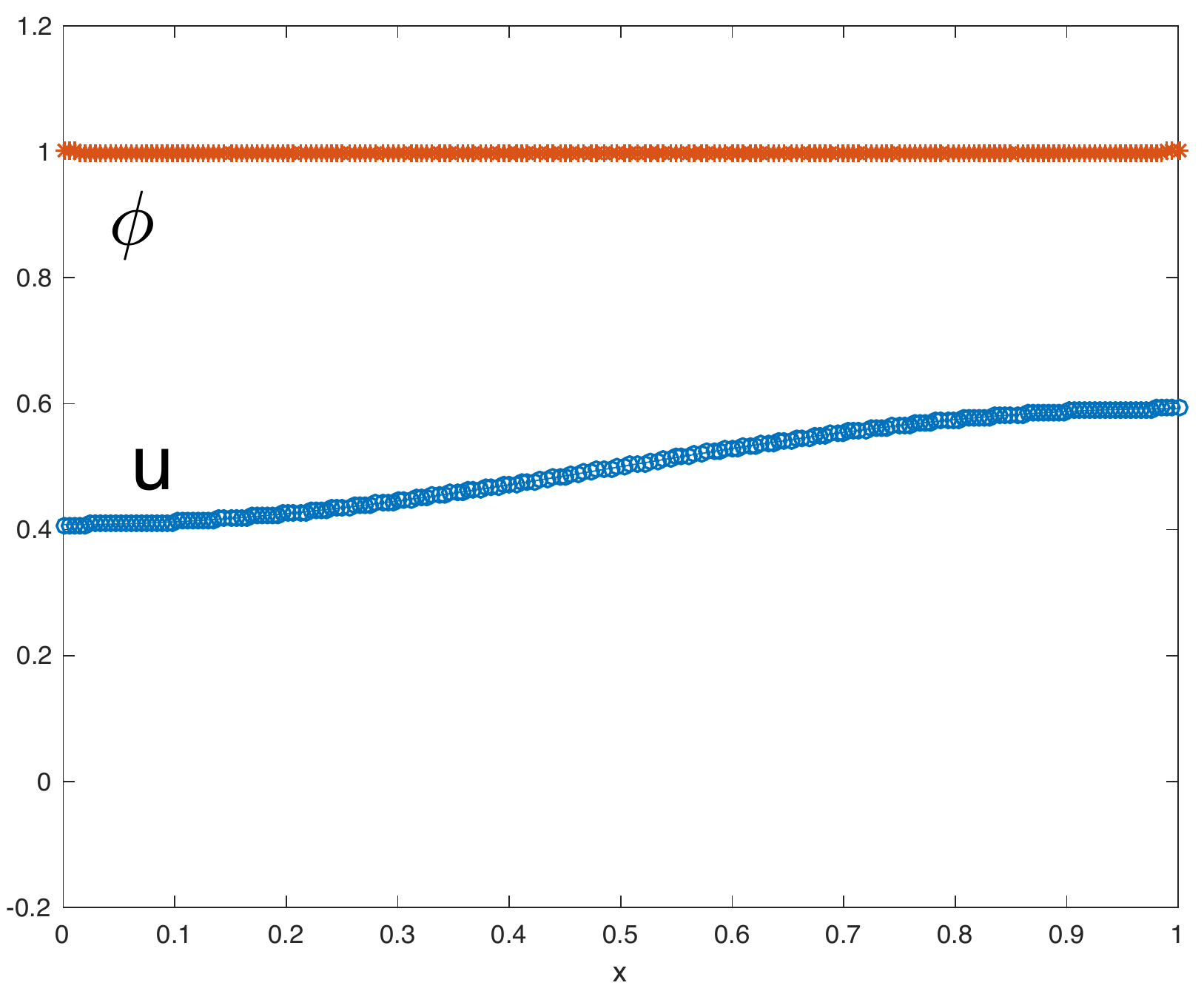}
\caption{ $t = 20$, $\epsilon = 10^{-5}$}
\end{subfigure}
\caption{Example~\ref{ex41}.
The computed solution $u_h$ and $\phi_h$ at three time instants for various values of $\epsilon$.
No scaling has been used on $g$ and $u$.}
\label{fig:3.11}
\end{figure}

\begin{exam}[1D hyperbolic tangent]
In this example, we take
\begin{equation}
g = 0.5(1+\tanh (100(x-0.5))), \quad x \in (0,1)
\label{ex41-1}
\end{equation}
which has a sharp jump at $x = 0.5$.
The initial conditions are taken as $u^0 = g$ and $ \phi^0 = 1$.
We take $N = 200$, $\alpha = 0.01$, $\beta = 10^{-3}$, $ \gamma = 10^{-3}$, and $k_{\epsilon} = 10^{-9}$.
The computed solution at three time instants for $\epsilon = 0.1$, $0.01$, and $10^{-5}$ is shown
in Fig. \ref{fig:3.11}. It can be seen that the mesh concentrates around and follows the sharp jumps
in the solution. This demonstrates the mesh adaptation ability of the MMPDE moving mesh strategy.

Recall that the jump in the solution simulates object edges in a real image and an ideal segmentation
should sharpen this jump while smoothing out the regions divided by the jump.
The first row of Fig.~\ref{fig:3.11} shows the evolution of $u$ and $\phi$ for $\epsilon = 0.1$.
One can see that the jump is not sharpened and $u$ is smoothed out on the whole domain as time evolves.
This indicates that the Ambrosio-Tortorelli functional with $\epsilon = 0.1$ does not provide a good segmentation.
The result is shown for a smaller $\epsilon = 0.01$ on the second row of the figure.
As time evolves, the jump is getting sharper and $u$ becomes piecewise constant essentially,
an indication for good image segmentation.
However, when $\epsilon$ continues to decrease, as shown on the last row ($\epsilon = 10^{-5}$)
of Fig.~\ref{fig:3.11},
the jump disappears for the time being, $\phi$ approaches to $1$, and $u$ becomes smooth
over the whole domain. This implies that the Ambrosio-Tortorelli functional loses its segmentation ability
for very small $\epsilon$, consistent with the analysis in Section~\ref{SEC:analysis}. 

It is interesting to see the transient behavior of $\phi$. From the simplified equation (\ref{simplephi}),
we have $\phi_t = -\alpha |\nabla u|^2 \phi$ initially due to the initial condition $\phi = 1$.
Thus, we expect that $\phi$ decreases initially and this decrease is more significant in the regions where
$\nabla u$ is larger. This is confirmed in the numerical results; see Fig.~\ref{fig:3.11}(a,d,g).
As time evolves, the system reaches an equilibrium state and $\phi$ is approximately given by (\ref{phi-2}).
When $\epsilon$ is not too small and $\nabla u$ is sufficiently large at some places, then $\phi$ can
become close to zero at the places and this yields a good segmentation; see the second row of Fig.~\ref{fig:3.11}.
However, when $\epsilon$ is too small, $\phi$ will essentially become 1 everywhere and the functional
loses its segmentation ability (cf. the third row of Fig. \ref{fig:3.11}).

\label{ex41}
\end{exam}

\vskip 5mm

\begin{figure}[htb]
\centering
\begin{subfigure}{0.32\textwidth}
\centering
\includegraphics[scale = 0.25]{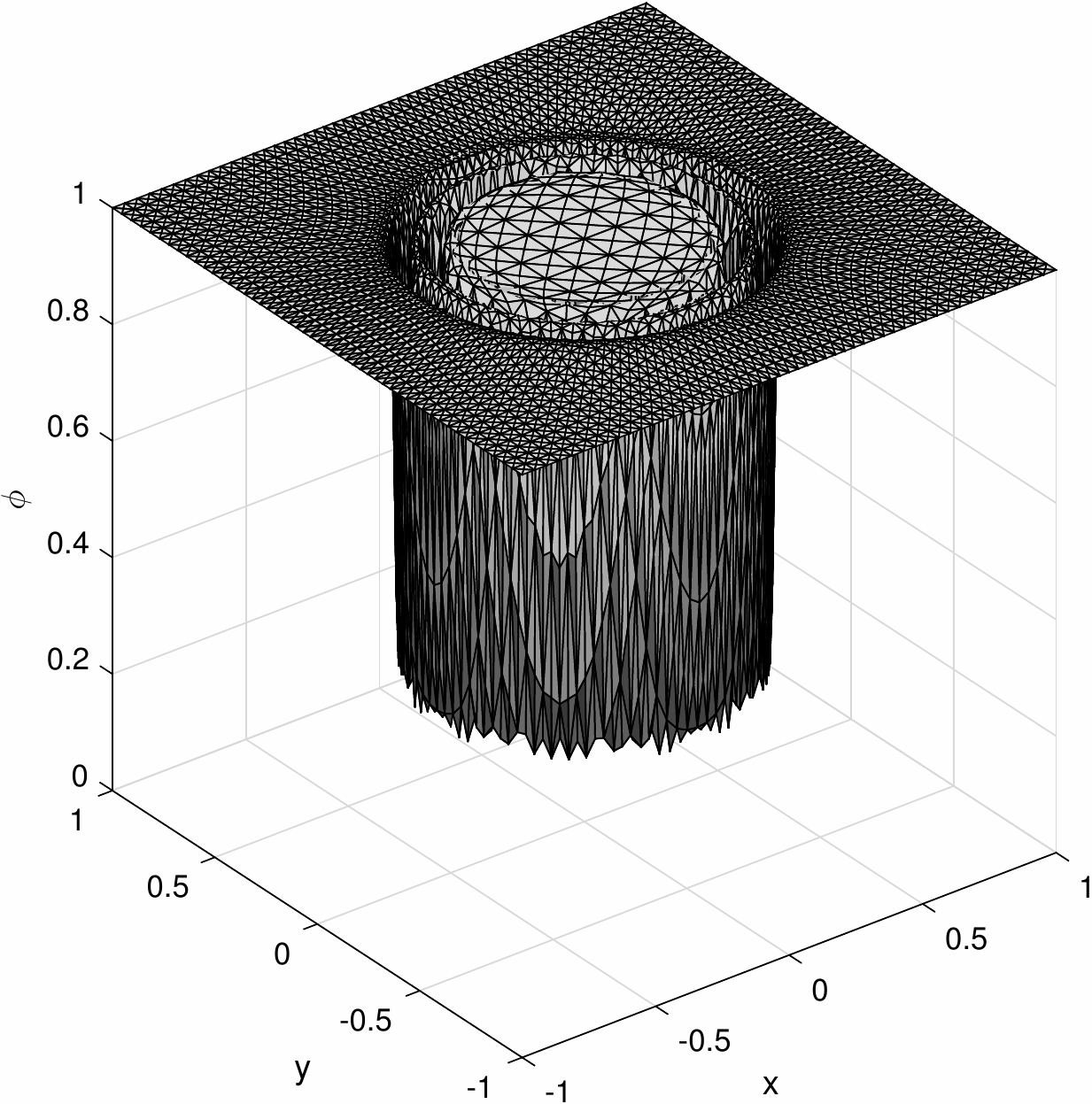}
\caption{t = 0.00015}
\end{subfigure}
\begin{subfigure}{0.32\textwidth}
\centering
\includegraphics[scale = 0.25]{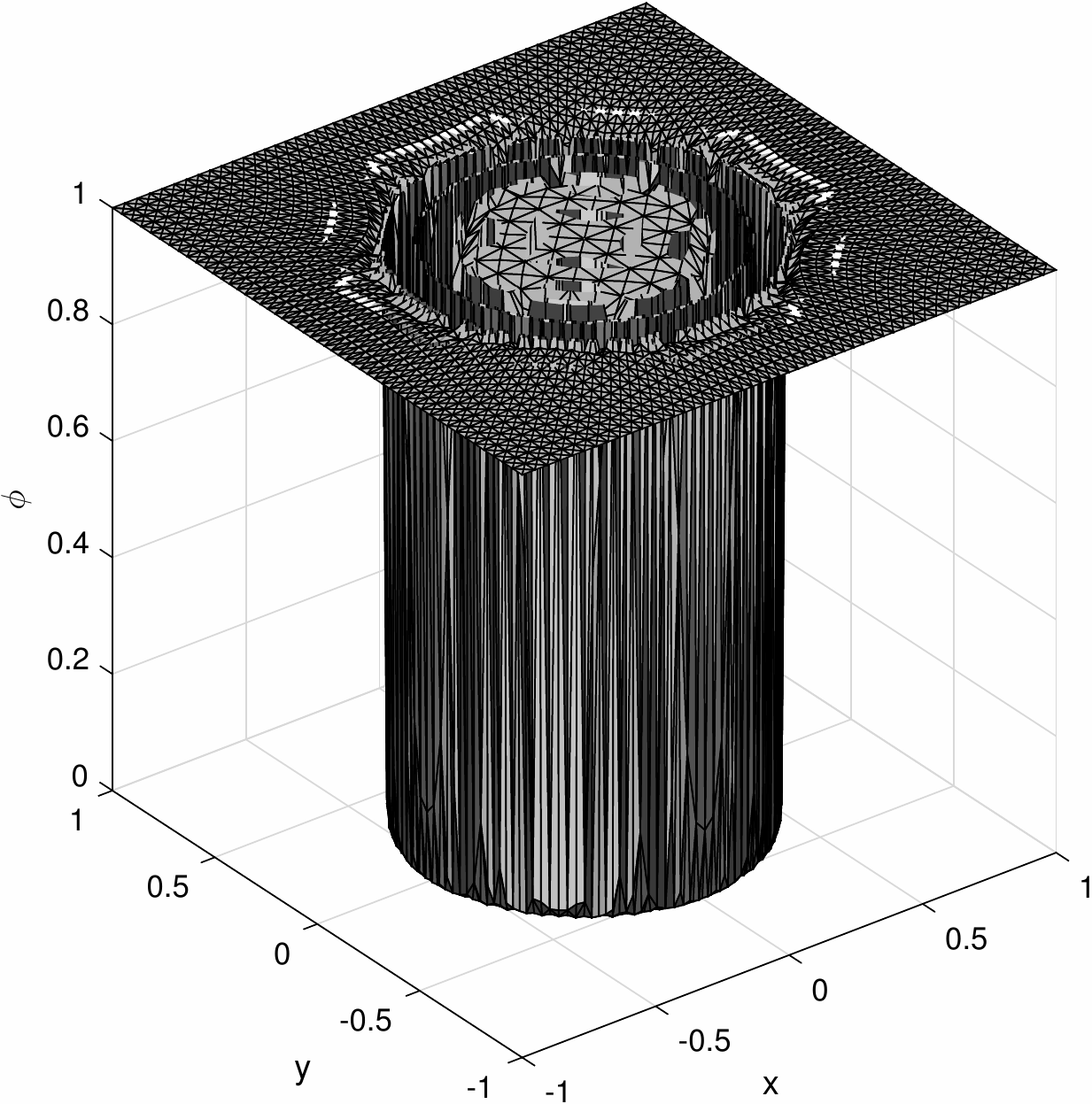}
\caption{t = 2}
\end{subfigure}
\begin{subfigure}{0.32\textwidth}
\centering
\includegraphics[scale = 0.25]{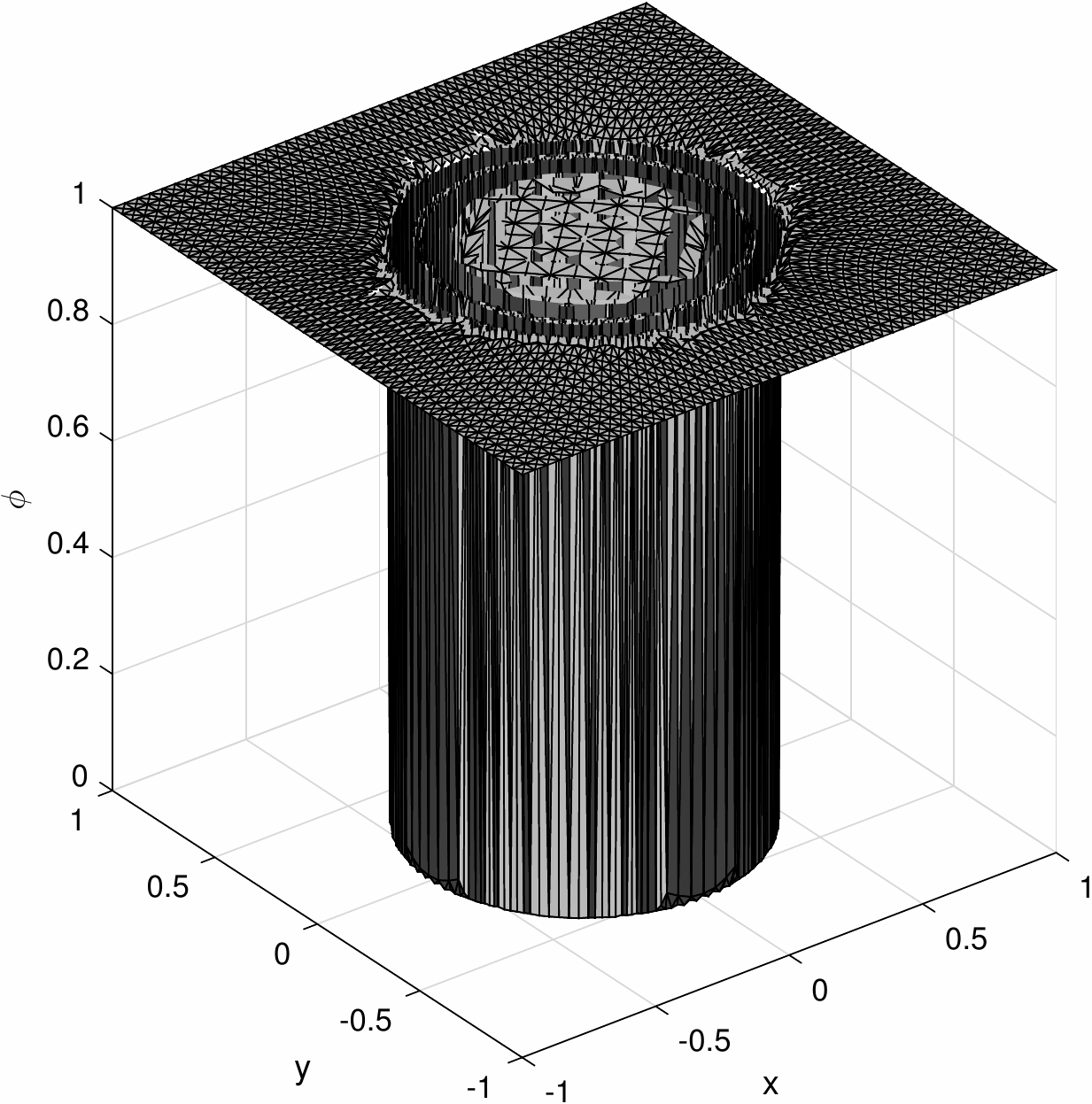}
\caption{t = 7}
\end{subfigure}
\begin{subfigure}{0.32\textwidth}
\centering
\includegraphics[scale = 0.24]{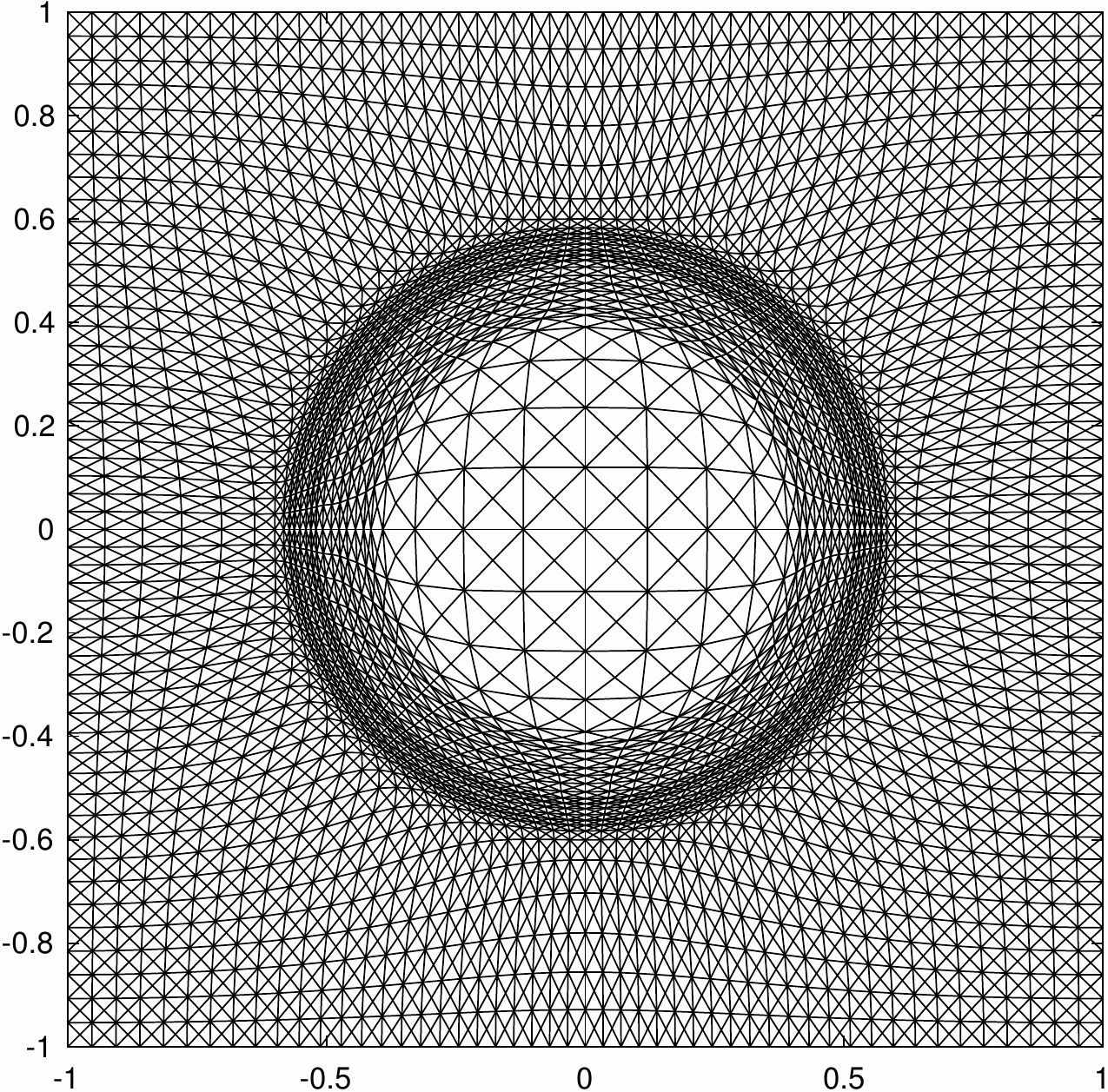}
\caption{t = 0.00015}
\end{subfigure}
\begin{subfigure}{0.32\textwidth}
\centering
\includegraphics[scale = 0.24]{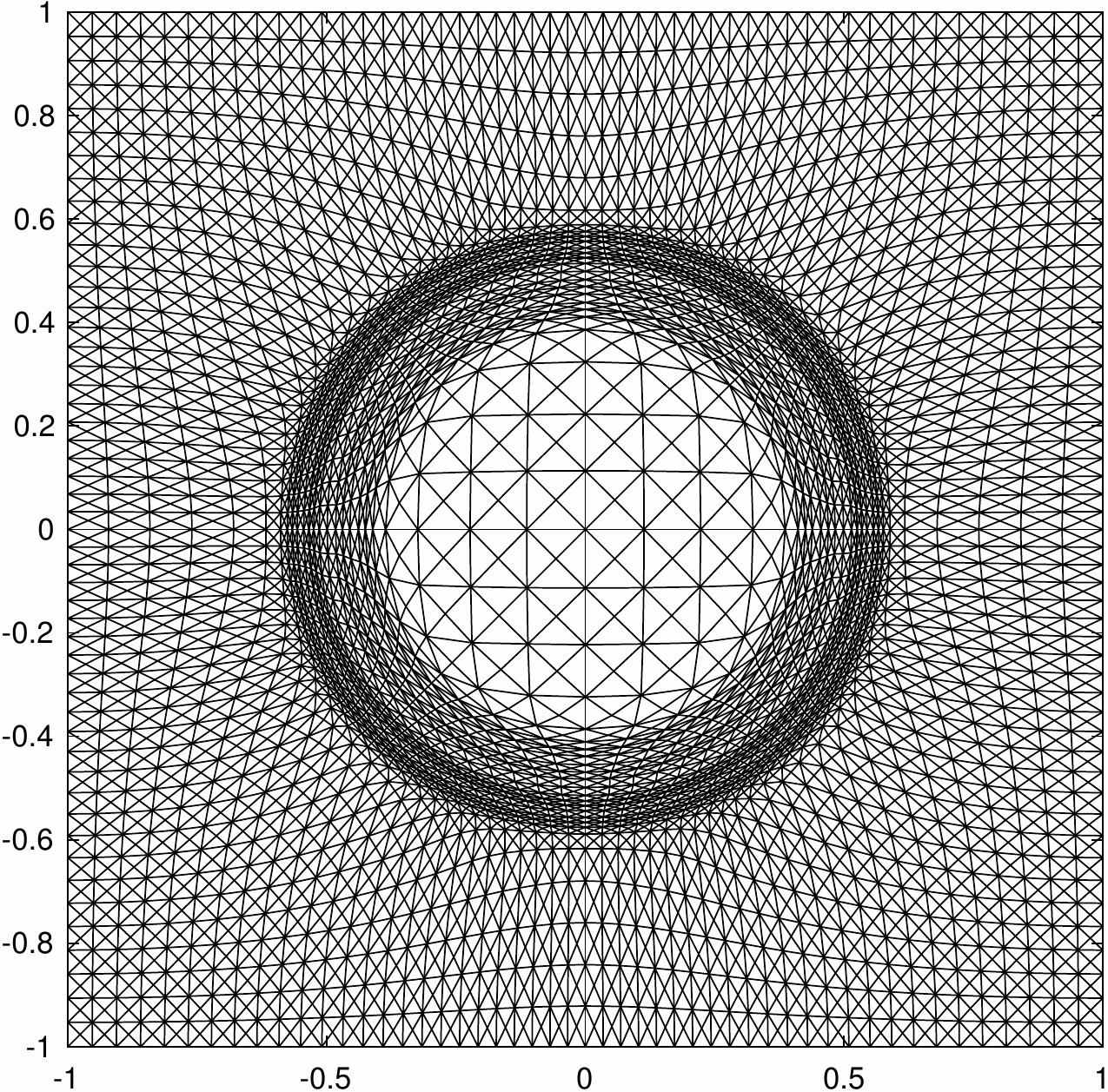}
\caption{t = 2}
\end{subfigure}
\begin{subfigure}{0.32\textwidth}
\centering
\includegraphics[scale = 0.24]{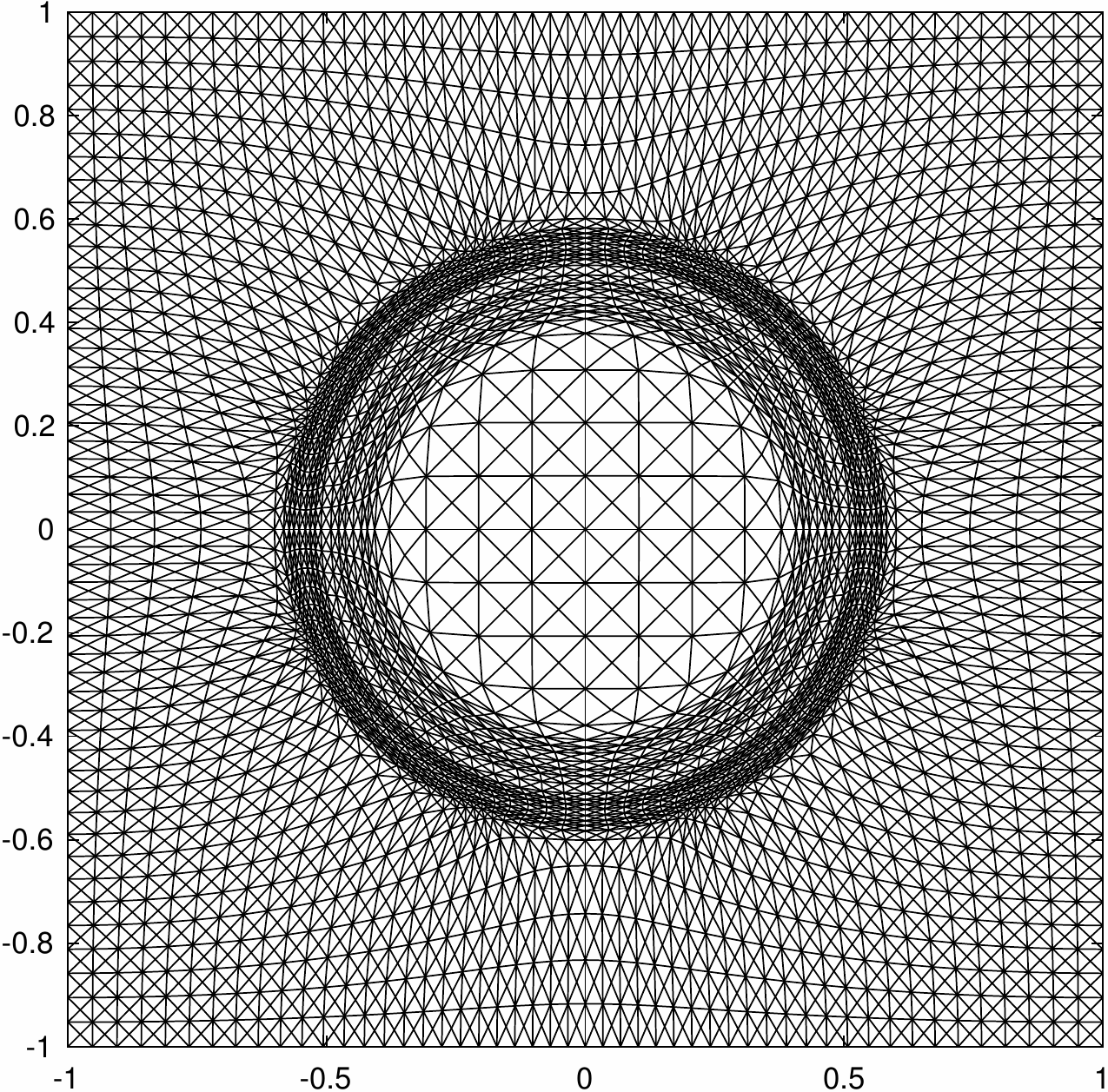}
\caption{t = 7}
\end{subfigure}
\begin{subfigure}{0.32\textwidth}
\centering
\includegraphics[scale = 0.33]{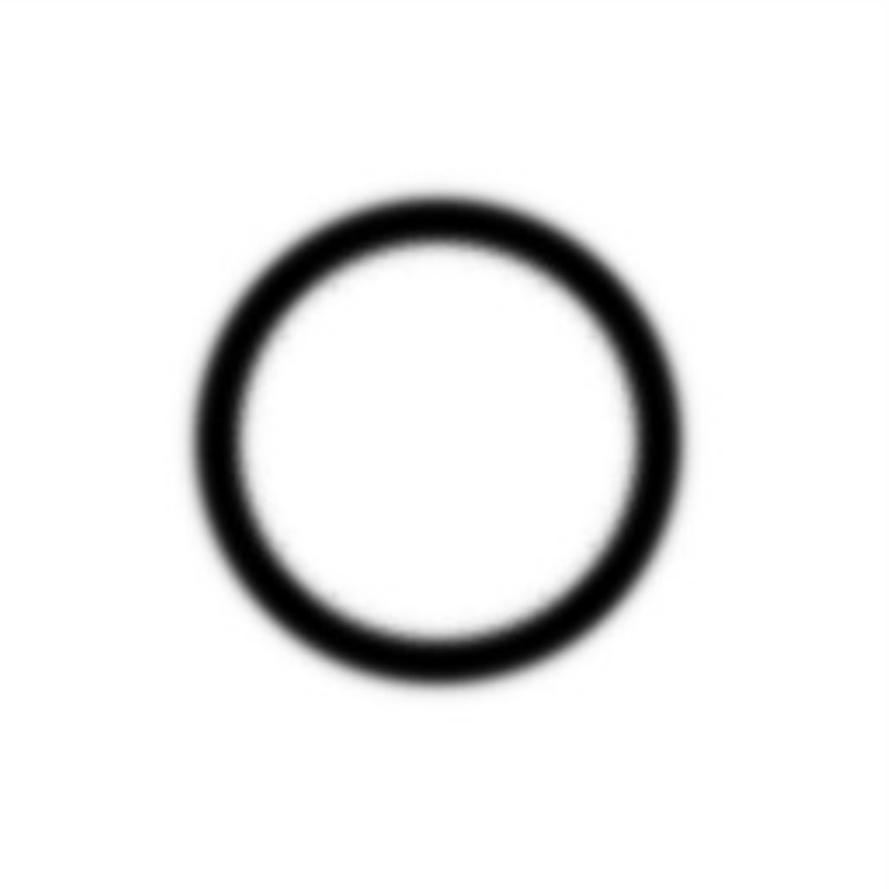}
\caption{t = 0.00015}
\end{subfigure}
\begin{subfigure}{0.32\textwidth}
\centering
\includegraphics[scale = 0.33]{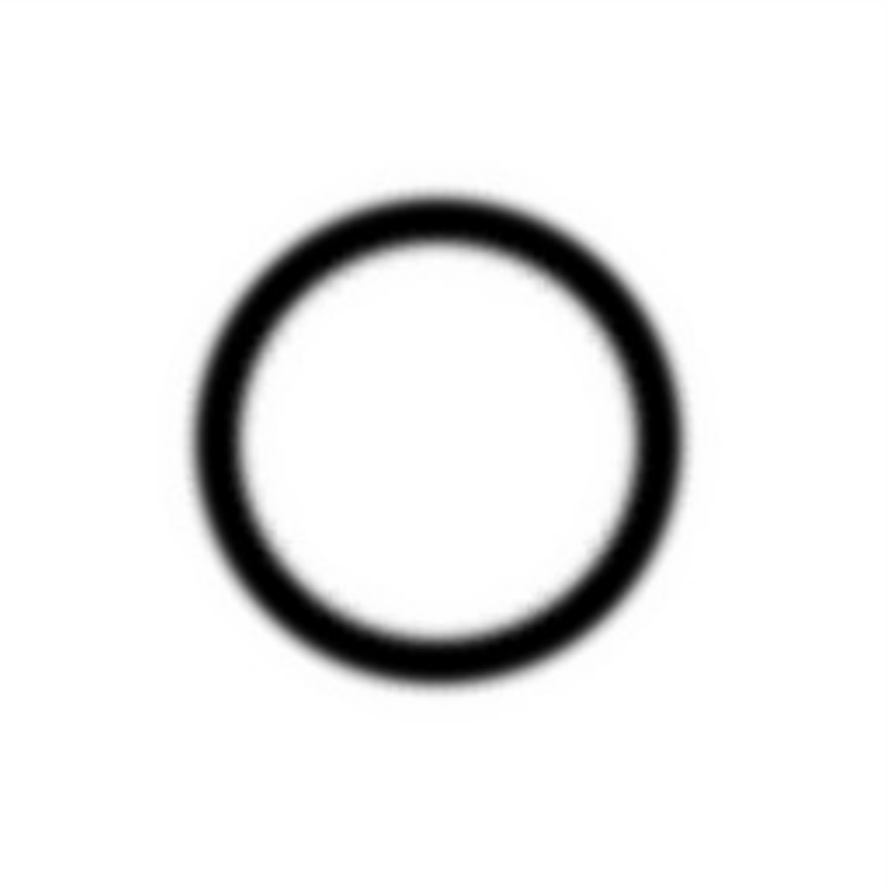}
\caption{t = 2}
\end{subfigure}
\begin{subfigure}{0.32\textwidth}
\centering
\includegraphics[scale = 0.33]{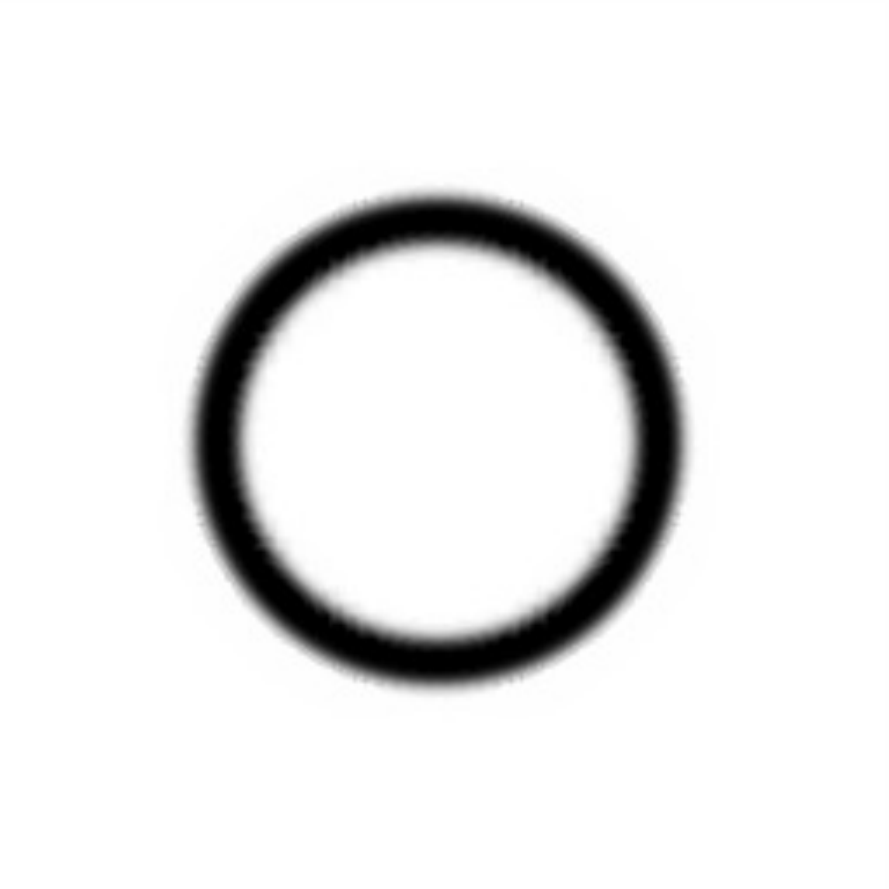}
\caption{t = 7}
\end{subfigure}
\caption{Example~\ref{ex42}. Evolution of the solution for $ \epsilon = 10^{-3}$.
The first, second, and third rows show the evolution of $\phi$, the moving mesh,  and the image of $u$, respectively.}
\label{fig:3.21}
\end{figure}

\begin{figure}[htb]
\centering
\begin{subfigure}{0.32\textwidth}
\centering
\includegraphics[scale = 0.25]{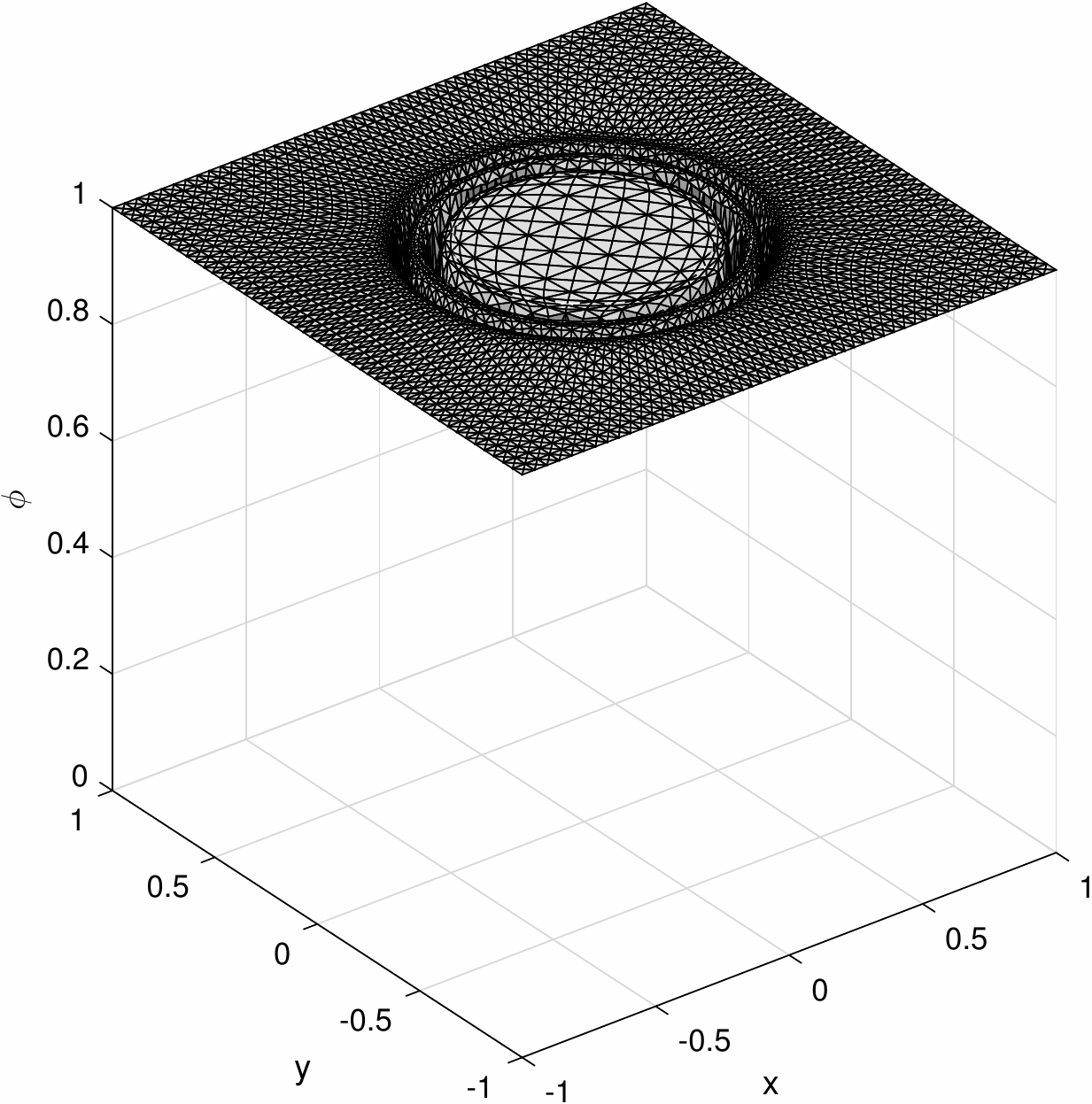}
\caption{t = 0.00015}
\end{subfigure}
\begin{subfigure}{0.32\textwidth}
\centering
\includegraphics[scale = 0.25]{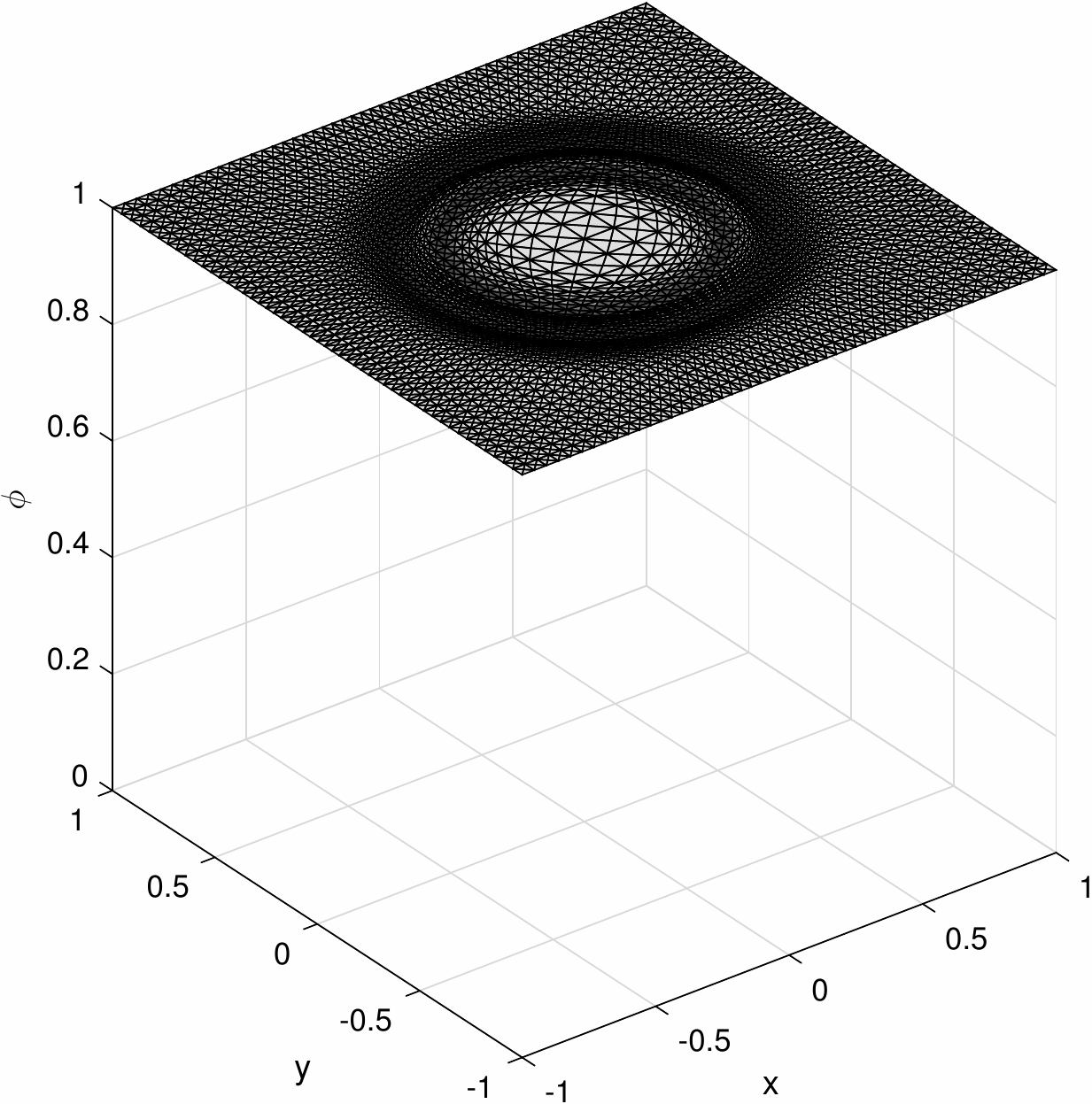}
\caption{t = 2}
\end{subfigure}
\begin{subfigure}{0.32\textwidth}
\centering
\includegraphics[scale = 0.25]{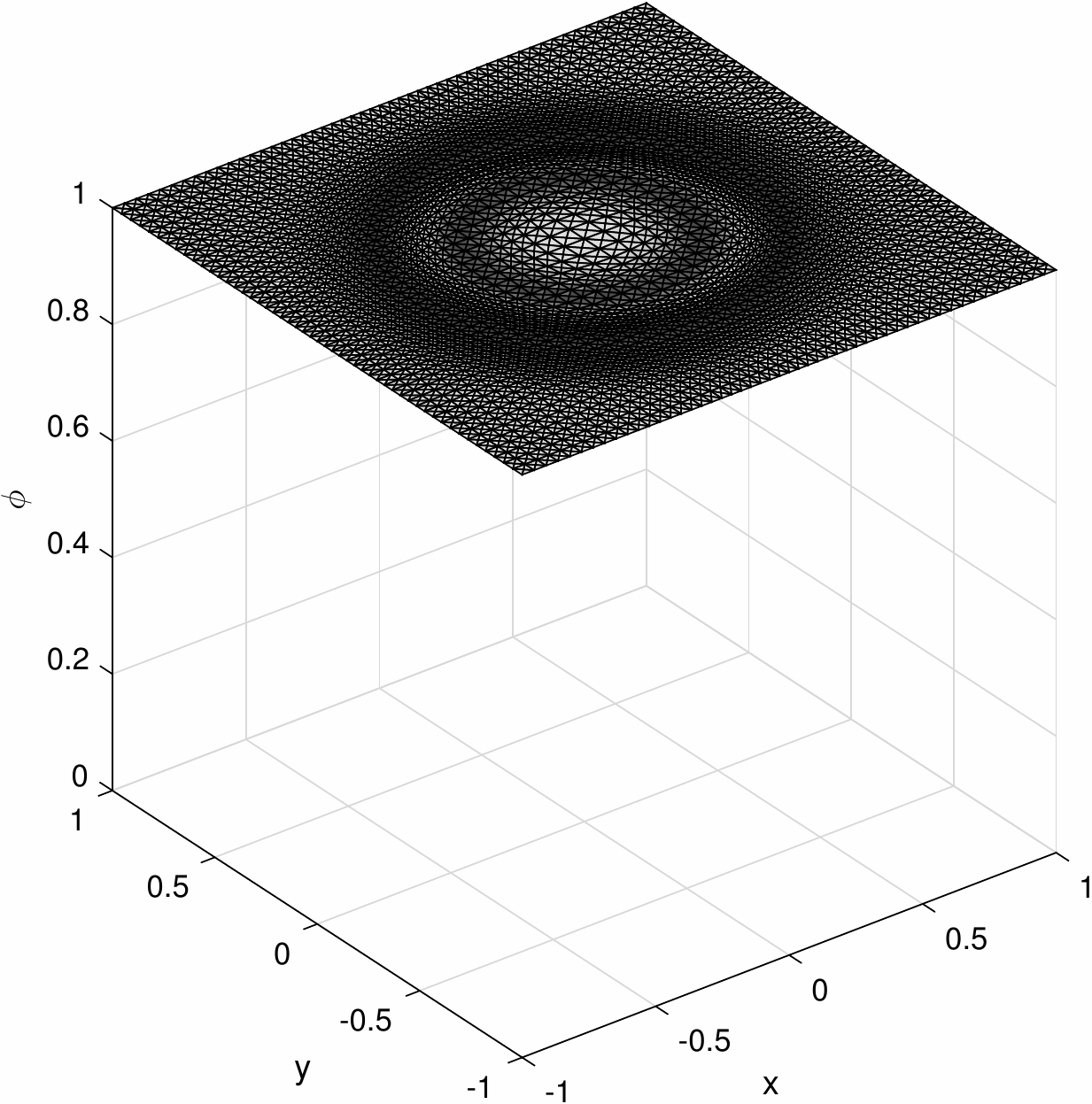}
\caption{t = 7}
\end{subfigure}
\begin{subfigure}{0.32\textwidth}
\centering
\includegraphics[scale = 0.24]{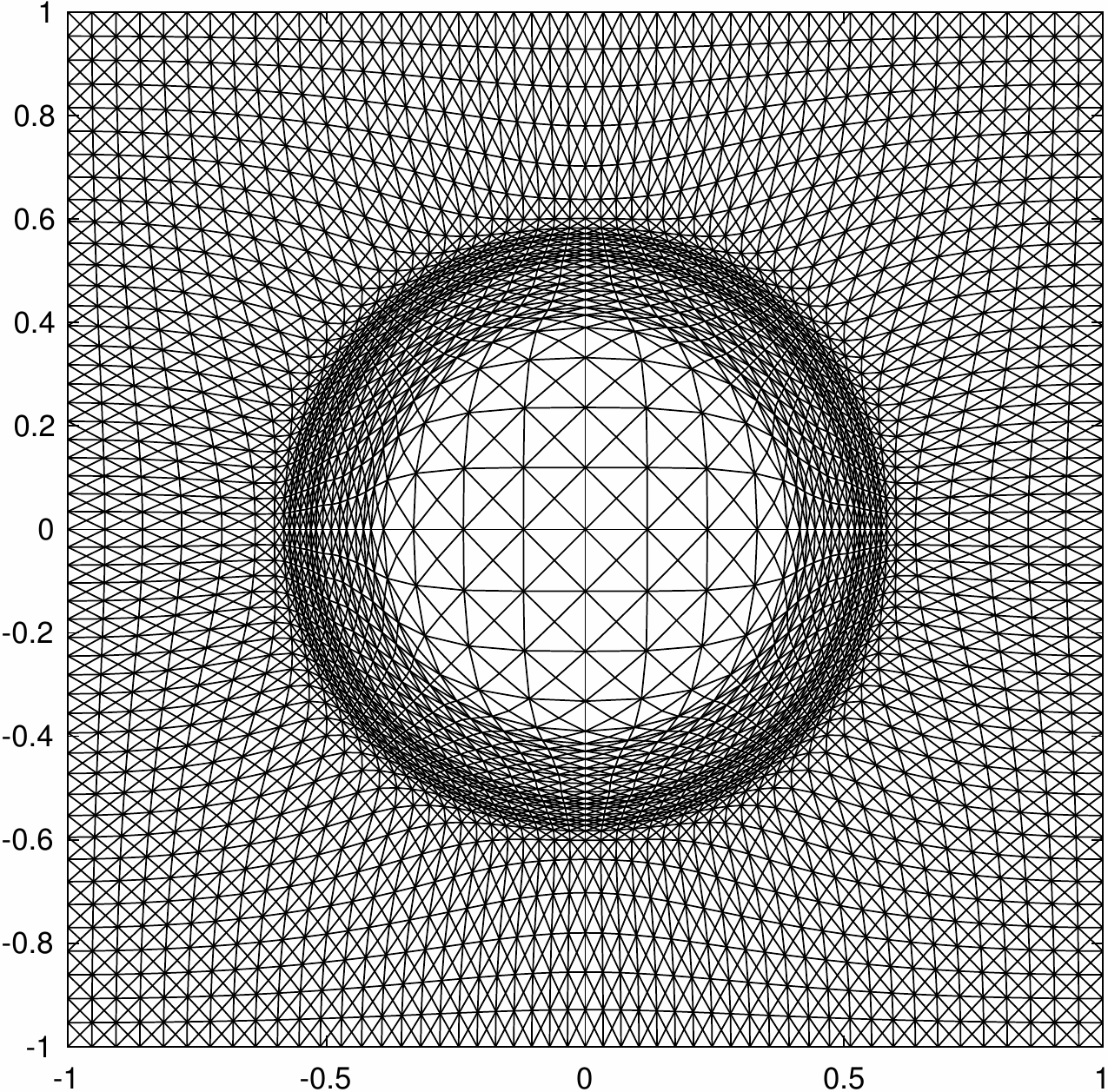}
\caption{t = 0.00015}
\end{subfigure}
\begin{subfigure}{0.32\textwidth}
\centering
\includegraphics[scale = 0.24]{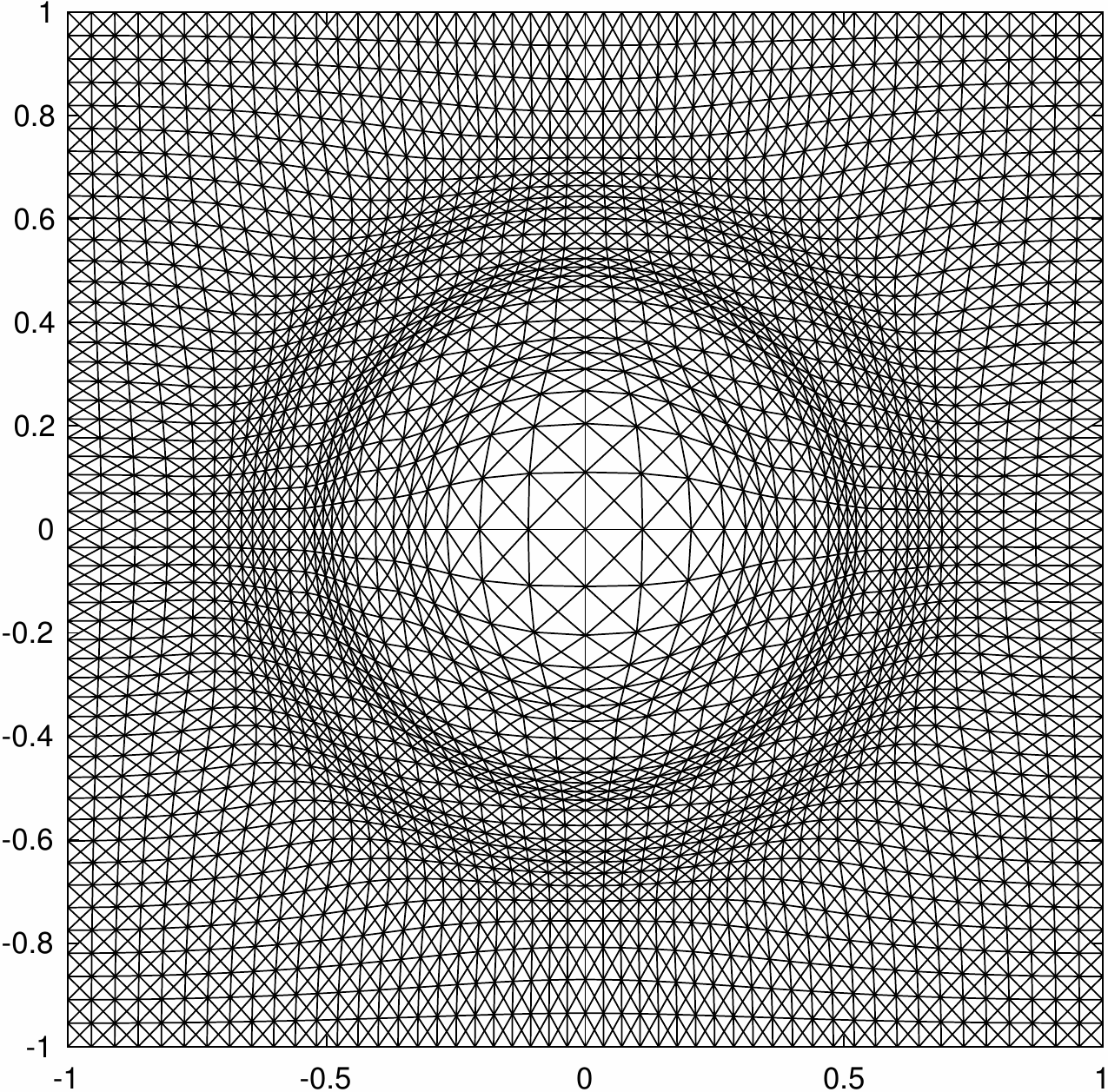}
\caption{t = 2}
\end{subfigure}
\begin{subfigure}{0.32\textwidth}
\centering
\includegraphics[scale = 0.24]{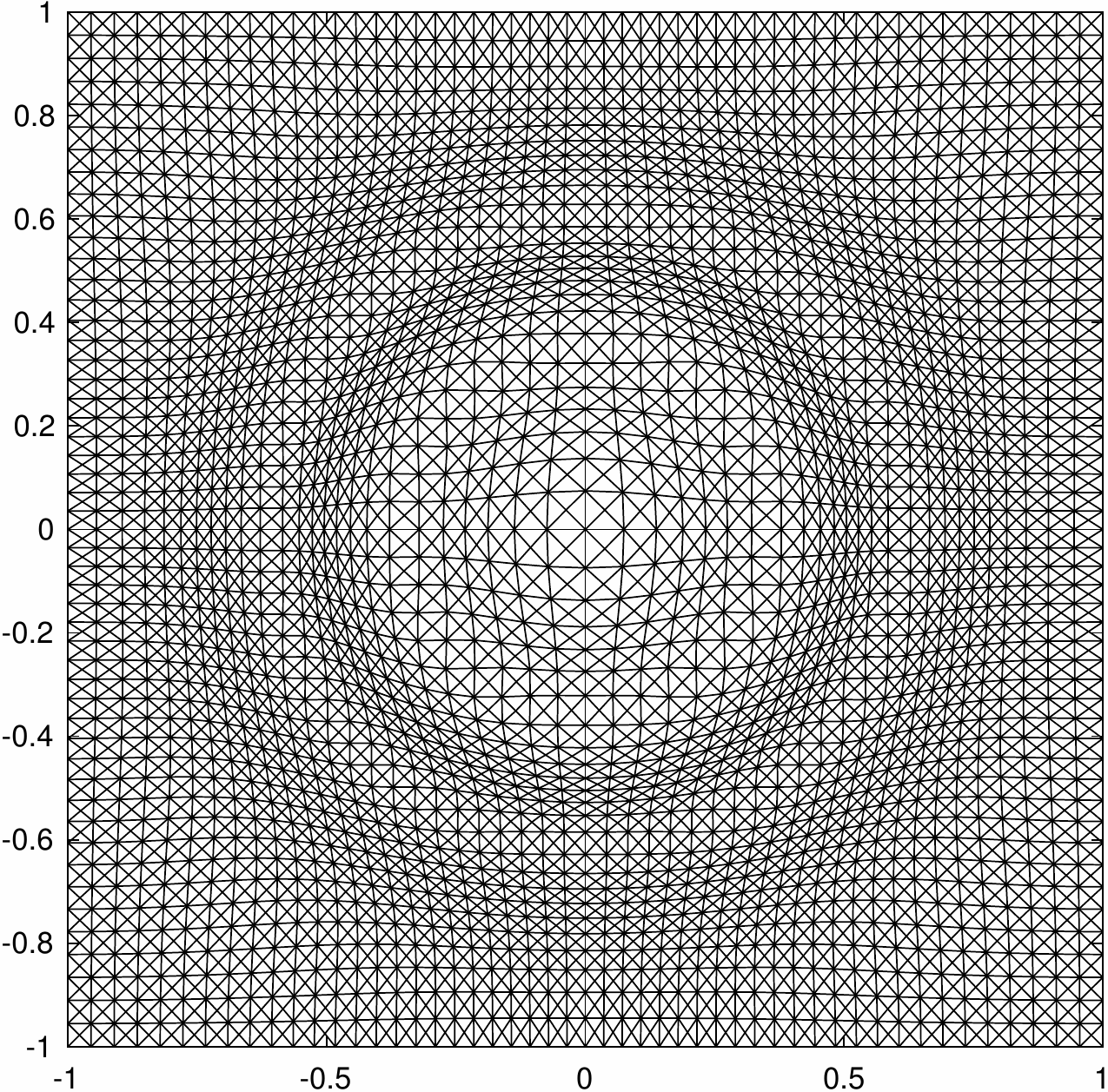}
\caption{t = 7}
\end{subfigure}
\begin{subfigure}{0.32\textwidth}
\centering
\includegraphics[scale = 0.33]{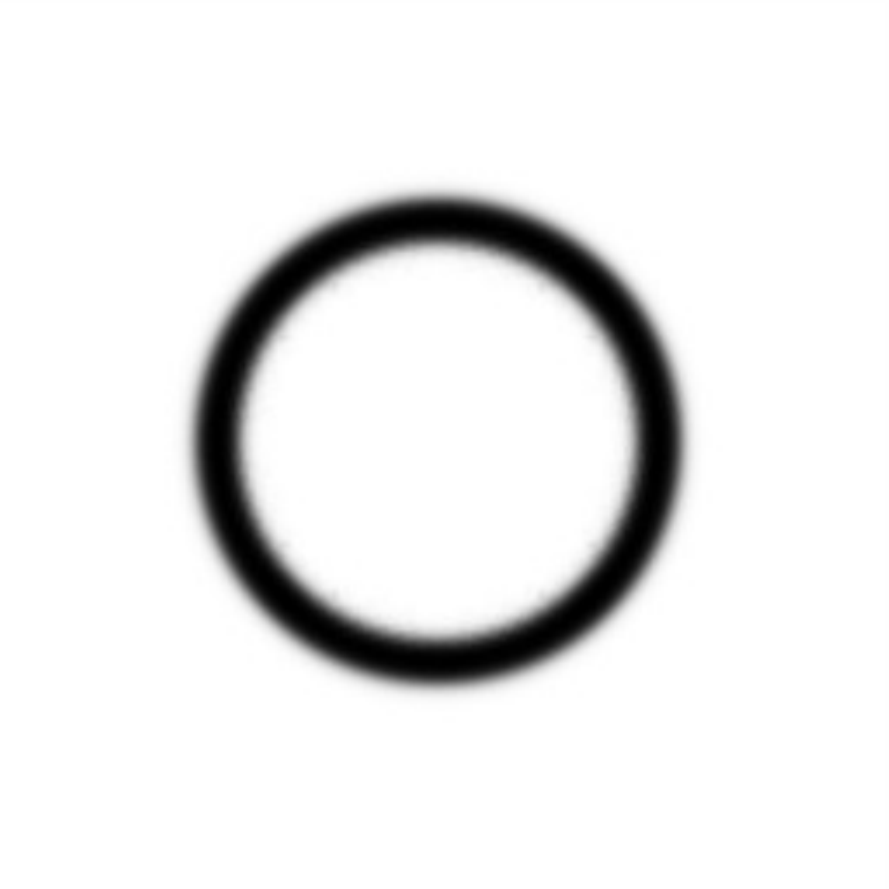}
\caption{t = 0.00015}
\end{subfigure}
\begin{subfigure}{0.32\textwidth}
\centering
\includegraphics[scale = 0.33]{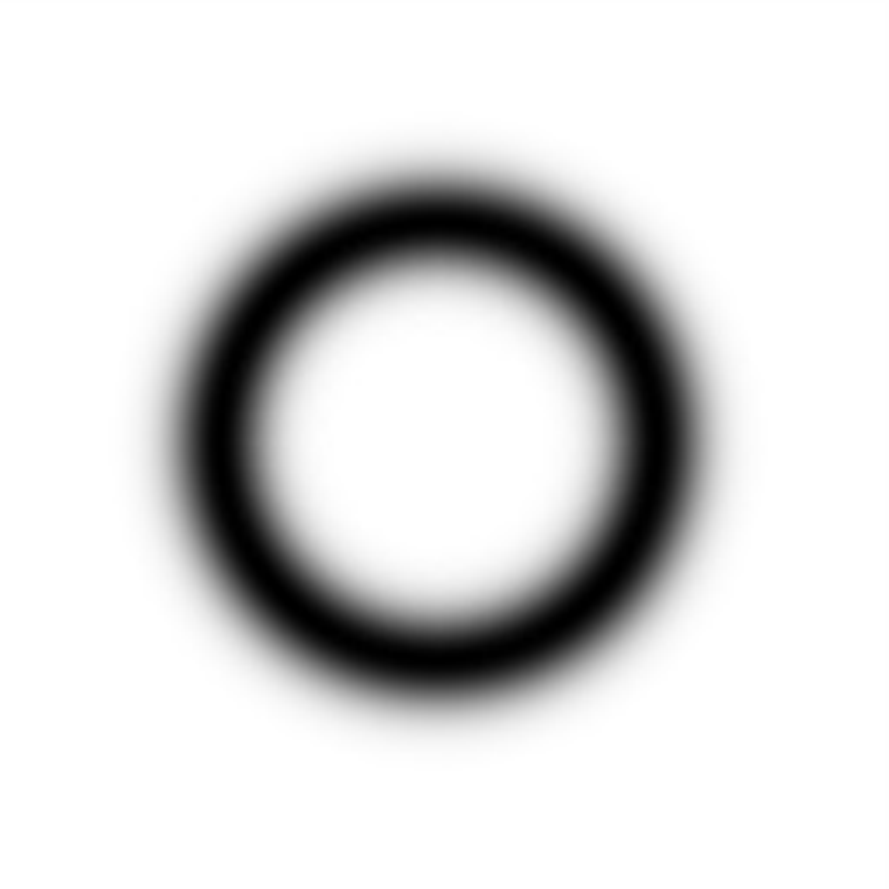}
\caption{t = 2}
\end{subfigure}
\begin{subfigure}{0.32\textwidth}
\centering
\includegraphics[scale = 0.33]{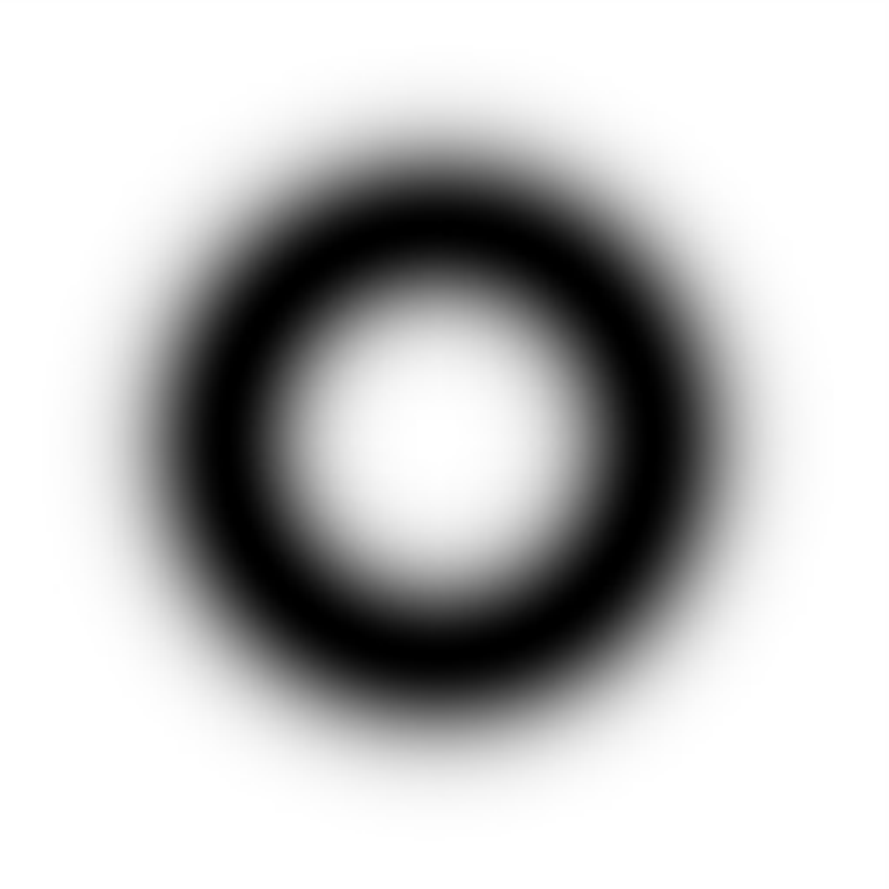}
\caption{t = 7}
\end{subfigure}
\caption{Example~\ref{ex42}. Evolution of the solution for $ \epsilon = 10^{-7}$.
The first, second, and third rows show the evolution of $\phi$, the moving mesh,  and the image of $u$, respectively.}
\label{fig:3.22}
\end{figure}

\begin{exam}[2D hyperbolic tangent]

In this example, we choose
\begin{align*}
g = &\, 0.49\left [ 2 + \tanh(50(\sqrt{(x-0.5)^2+(y-0.5)^2} - 0.05))\right .
\\
& \qquad \left. - \tanh(50(\sqrt{(x-0.5)^2+(y-0.5)^2} + 0.05))\right ],
\quad (x,y) \in (0,1) \times (0,1) 
\end{align*}
which models a circle, being close to 0 on the circle and approximately $1$ elsewhere.
For the reasons to be explained in Section~\ref{SEC:select}, $u$ and $g$ in the IBVP (\ref{grad-flow-1})
are scaled in this example according to (\ref{L-1}). 

We take $u^{0} = g$, $\phi^{0} = 1$, $N = 2\times 50 \times 50$, $\alpha = 10^{-3}$, $\gamma = 10^{-5}$,
$\beta = 10^{-2}$, and $k_{\epsilon} = 10^{-10}$. The numerical results obtained with $\epsilon = 10^{-3}$ and
$\epsilon = 10^{-7}$ are shown in Figs.~\ref{fig:3.21} and \ref{fig:3.22}, respectively.
They show that the mesh concentrates around the jump (the circle) very well, which, once again, demonstrates
the mesh adaptation ability of the MMPDE moving mesh method.

Fig. \ref{fig:3.21} shows that the Ambrosio-Tortorelli functional with $\epsilon = 10^{-3}$
makes a good segmentation. The evolution of $\phi$ is given on the first row, and $\phi$ deceases rapidly to $0$
along the circle at $t = 7$. The image of the circle is clear as shown on the third row. 
However, the situation changes when a smaller $\epsilon$ is used. 
As shown in Fig. \ref{fig:3.22} with $\epsilon = 10^{-7}$, the segmentation ability disappears. As $t$ increases,
$\phi$ becomes close to 1 in the whole domain, failing to identify the circle. In the same time, the image of $u$ blurs out.
As for Example~\ref{ex41}, the above observation is consistent with the analysis in Section~\ref{SEC:analysis},
that is, when $g$ is continuous, the segmentation ability of the Ambrosio-Tortorelli functional varies for small
but finite $\epsilon$ and disappears as $\epsilon \to 0$.
\label{ex42}
\end{exam}

\section{Selection of the regularization parameter and scaling of $g$ and $u$}
\label{SEC:select}

\subsection{Selection of the regularization parameter}

From the analysis in Section~\ref{SEC:analysis} and the examples in the previous section,
we have seen that it is crucial to choose a proper $\epsilon$ for the Ambrosio-Tortorelli functional to produce
a good segmentation when $g$ is continuous. To see how to choose $\epsilon$ properly,
we recall that $\phi$ is given in (\ref{phi-0})
for small $\epsilon$.  We want to have $\phi = 0 $ on object edges. Taking $\phi = 0 $ in (\ref{phi-0}) we get
\[
\epsilon = \frac{\beta}{2\alpha |\nabla u^{(0)}|^2},
\]
where $u^{(0)}$ is the solution of (\ref{u-3}) subject to a homogeneous Neumann boundary condition.
Since $u^{(0)}$ is completely determined by its initial value $g$ and an objective of the Ambrosio-Tortorelli functional
is to make $u$ (and thus $u^{(0)}$) close to $g$, it is reasonable to replace $u^{(0)}$
by $g$ in the above formula, i.e.,
\[
\epsilon = \frac{\beta}{2\alpha |\nabla g|^2} .
\]
Since $|\nabla g|$ varies from place to place and $\epsilon$ is a constant, in our computation we replace
the former with $( |\nabla g|_{\text{max}} + |\nabla g|_{\text{min}})/2$ and have 
\begin{equation}
\epsilon = \frac{\beta}{2\alpha \left (( |\nabla g|_{\text{max}} + |\nabla g|_{\text{min}})/2\right )^2} .
\label{epsilon-1}
\end{equation}

\begin{figure}[htb]
\centering
\begin{subfigure}{0.32\textwidth}
\centering
\includegraphics[scale = 0.24]{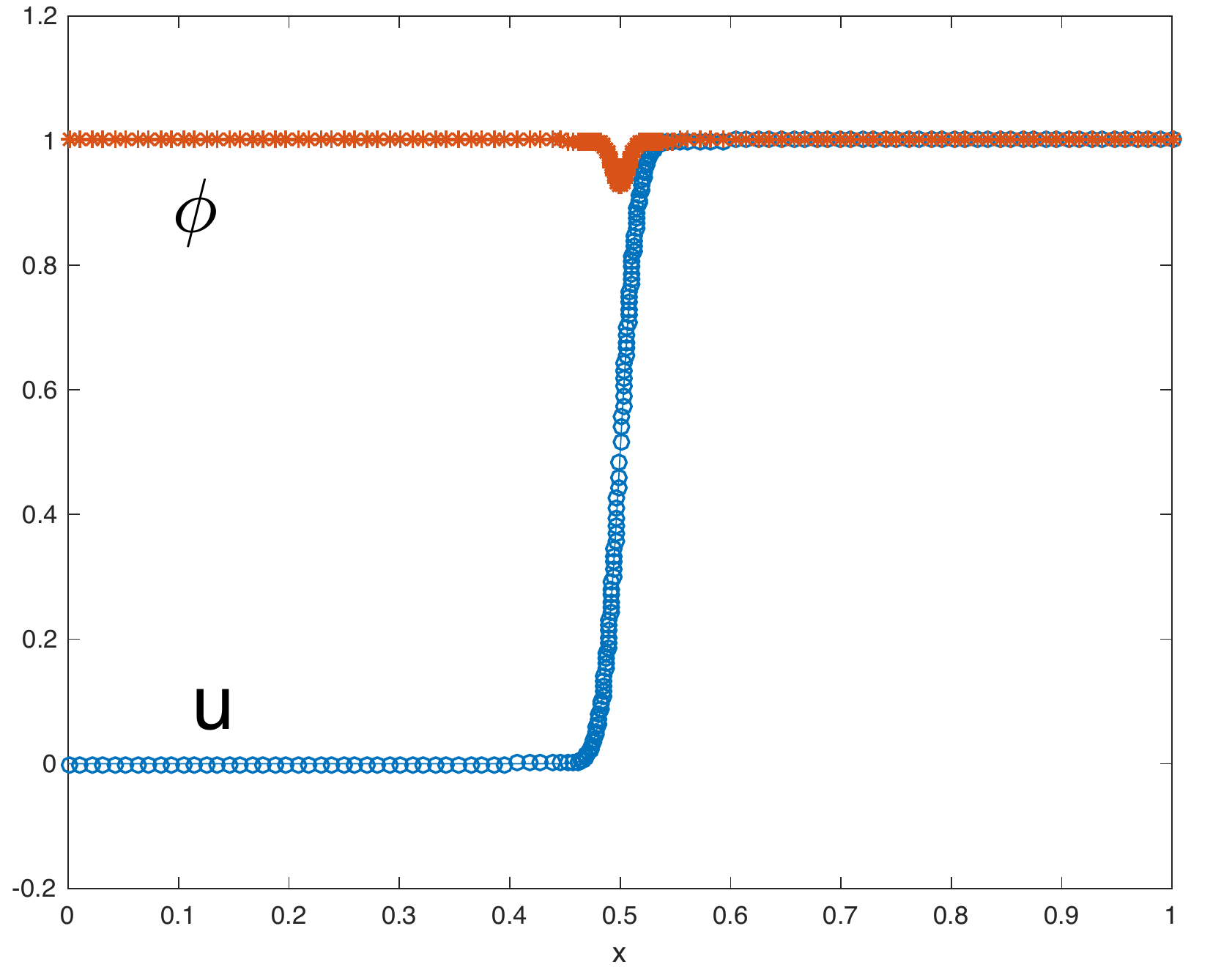}
\caption{t = 0.005}
\end{subfigure}
\begin{subfigure}{0.32\textwidth}
\centering
\includegraphics[scale = 0.24]{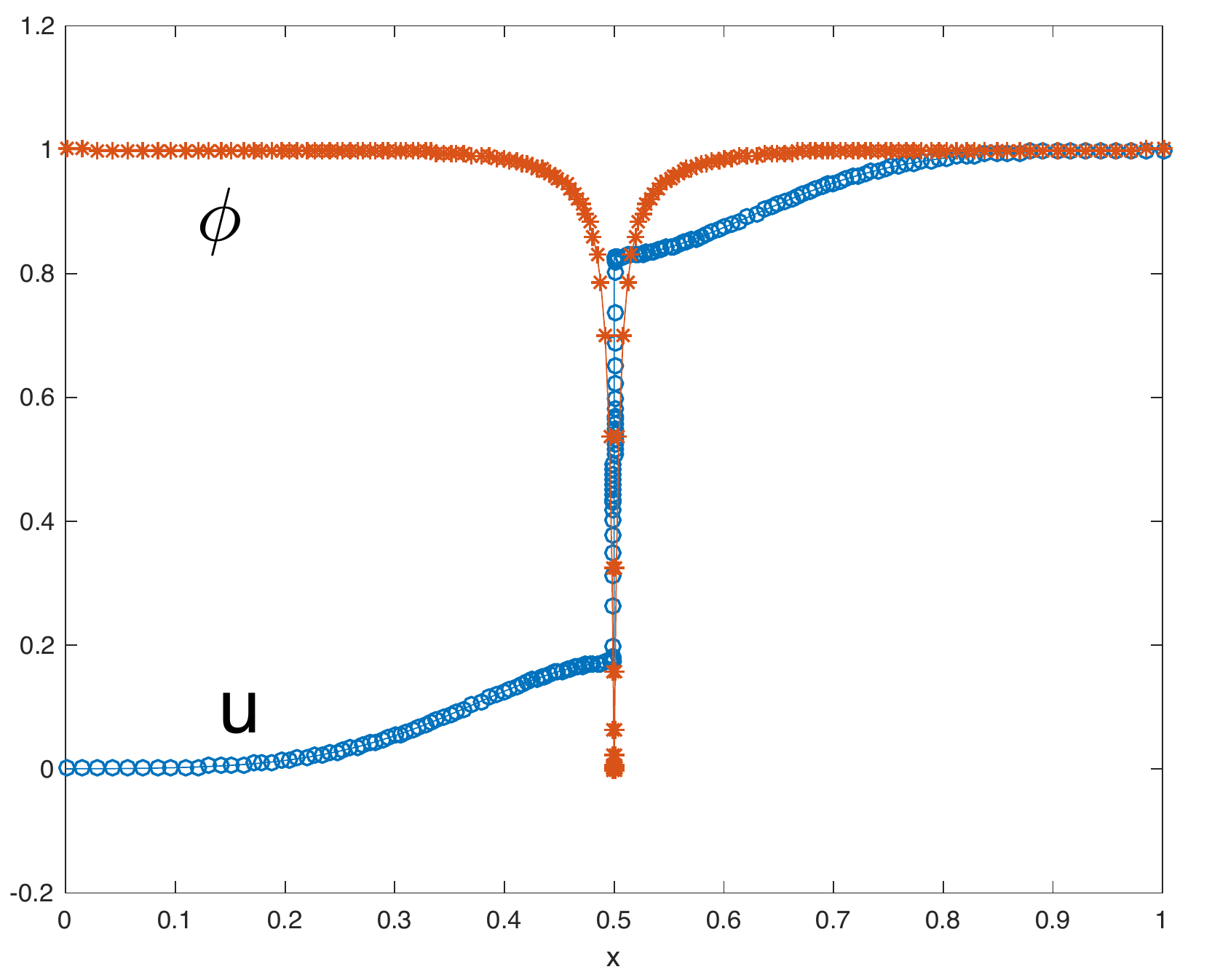}
\caption{t = 1}
\end{subfigure}
 \begin{subfigure}{0.32\textwidth}
 \centering
\includegraphics[scale = 0.24]{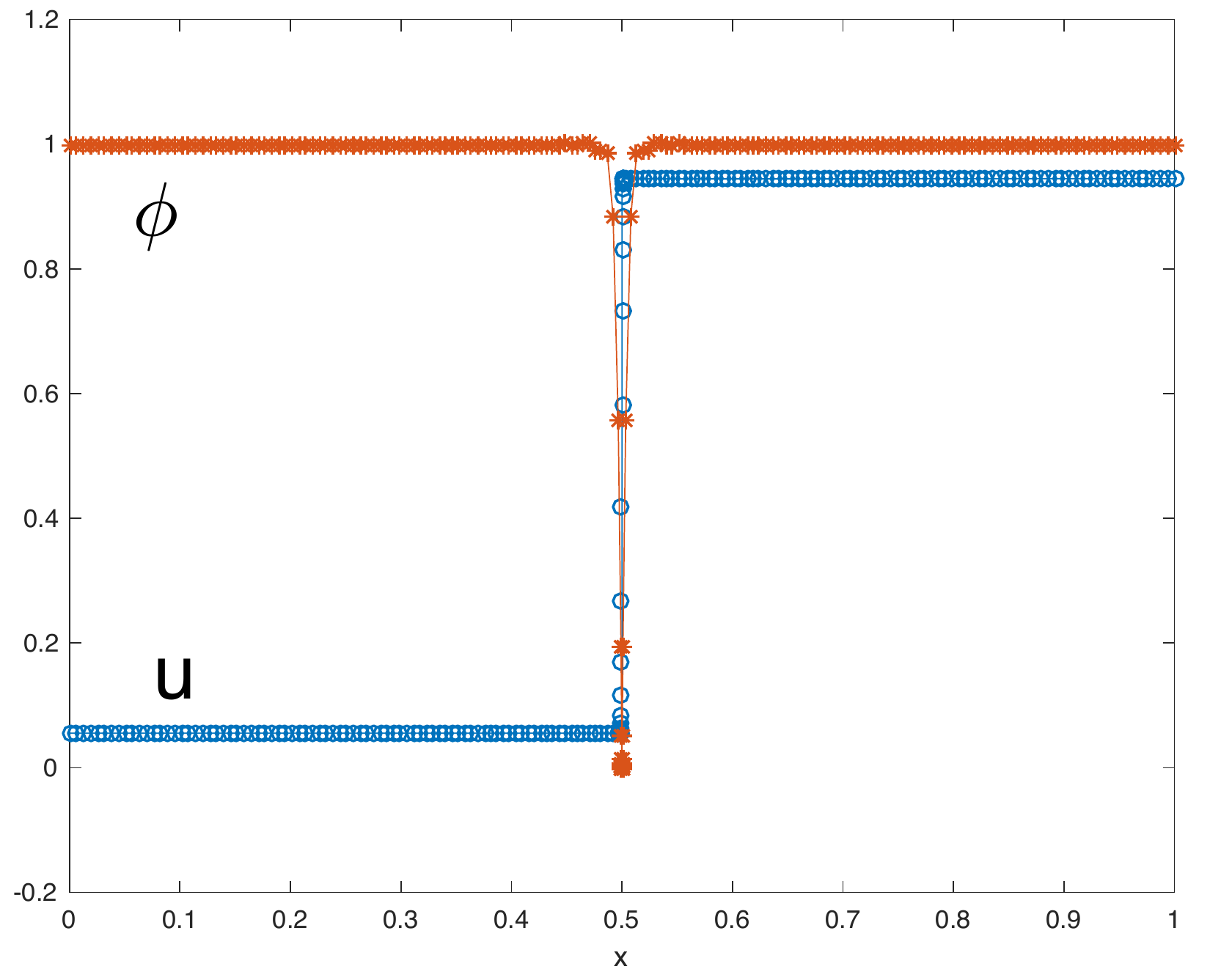}
\caption{t = 20}
\end{subfigure}
\caption{Example~\ref{ex41}. The evolution of $u$ and $\phi$ for $\epsilon = 0.008$ (determined by (\ref{epsilon-1})).
No scaling has been used on $u$ and $g$.}
\label{fig4c}
\end{figure}

To demonstrate this choice of $\epsilon$, we apply it to Example~\ref{ex41} and obtain $\epsilon = 0.008$.
The numerical result obtained with the same initial condition and parameters (other than $\epsilon$)
is shown in Fig. \ref{fig4c}. One can see that this value of $\epsilon$ leads to a good segmentation of the image.

\subsection{Scaling of $g$ and $u$}
\label{SEC:scaling}

Our experience shows that (\ref{epsilon-1}) works well when the difference in $\nabla g$ between the
objects and their edges is sufficiently large. However, when the change of $\nabla g$ is small,
the Ambrosio-Tortorelli functional can still fail to produce a segmentation of good quality.
To avoid this difficulty, we propose to scale $u$ and $g$ in (\ref{grad-flow-1}), i.e., $u \to L u$ and
$g \to L g$ for some parameter $L \ge 1$. This will make the change of $\nabla g$ from place to place
more significant. Moreover, the first equation of (\ref{grad-flow-1}) will stay invariant.
The second equation becomes
\[
\phi_{t} = 2\beta \epsilon \Delta \phi - L^2 \alpha | \nabla u|^{2}\phi + \frac{\beta}{2 \epsilon}(1 - \phi),
\]
where the second term on the right-hand side is made larger, helping decrease $\phi$. We choose
\begin{equation}
L =  \max \left \{1, \frac{|\nabla g|_{\text{cr}}}{|\nabla g|_{\max}}\right \},
\label{L-1}
\end{equation}
where $|\nabla g|_{\text{cr}}$ is a parameter. Generally speaking,
the larger $|\nabla g|_{\text{cr}}$ (and $L$) is, the more likely the segmentation works, but this will also make (\ref{grad-flow-1}) harder to integrate.
We take $|\nabla g|_{\text{cr}} = 3\times 10^3$ (by trial and error) in our computation, unless stated otherwise.

\begin{figure}[htb]
\centering
\begin{subfigure}{0.32\textwidth}
\centering
\includegraphics[scale = 0.24]{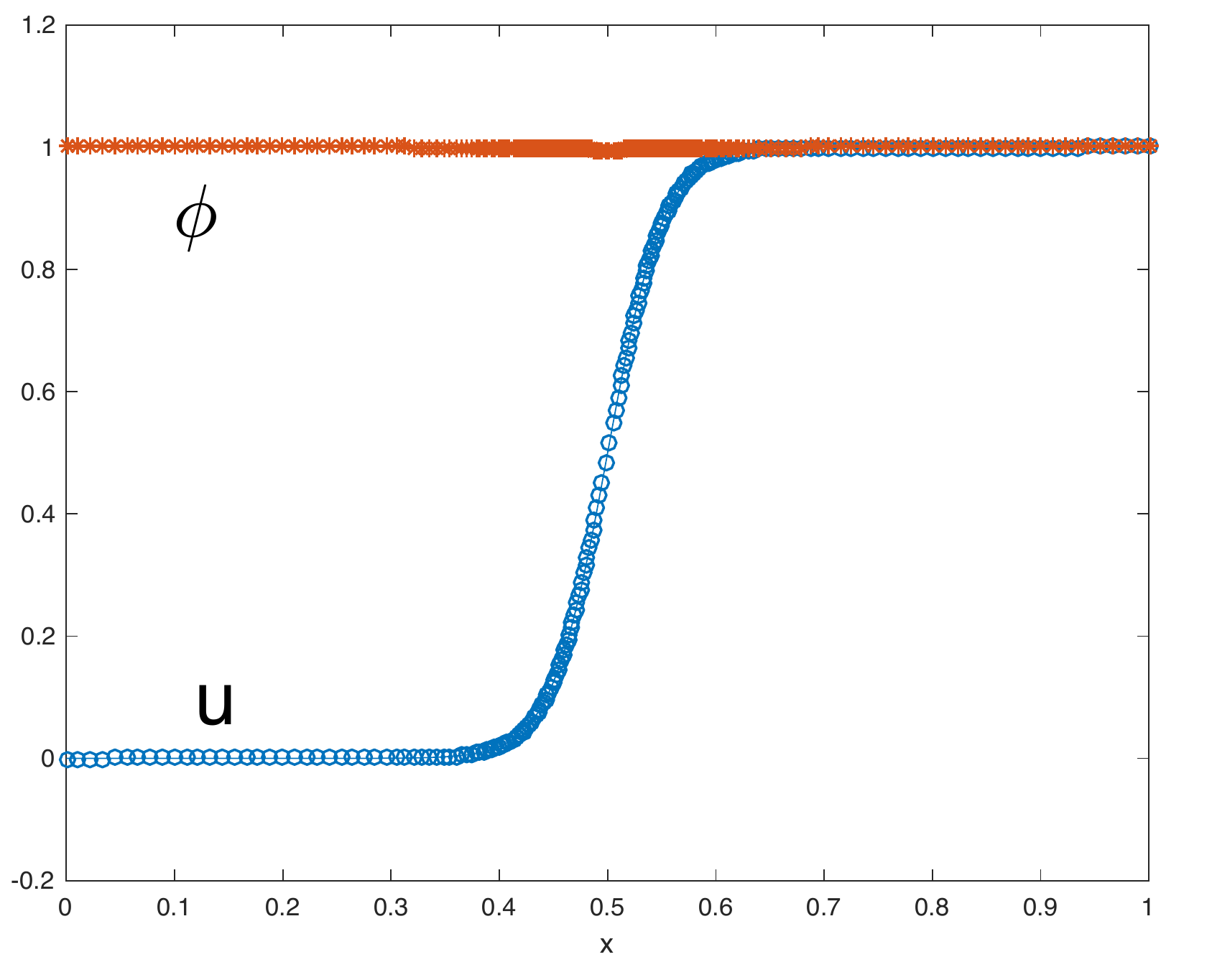}
\caption{ $t = 0.005$}
\end{subfigure}
\begin{subfigure}{0.32\textwidth}
\centering
\includegraphics[scale = 0.24]{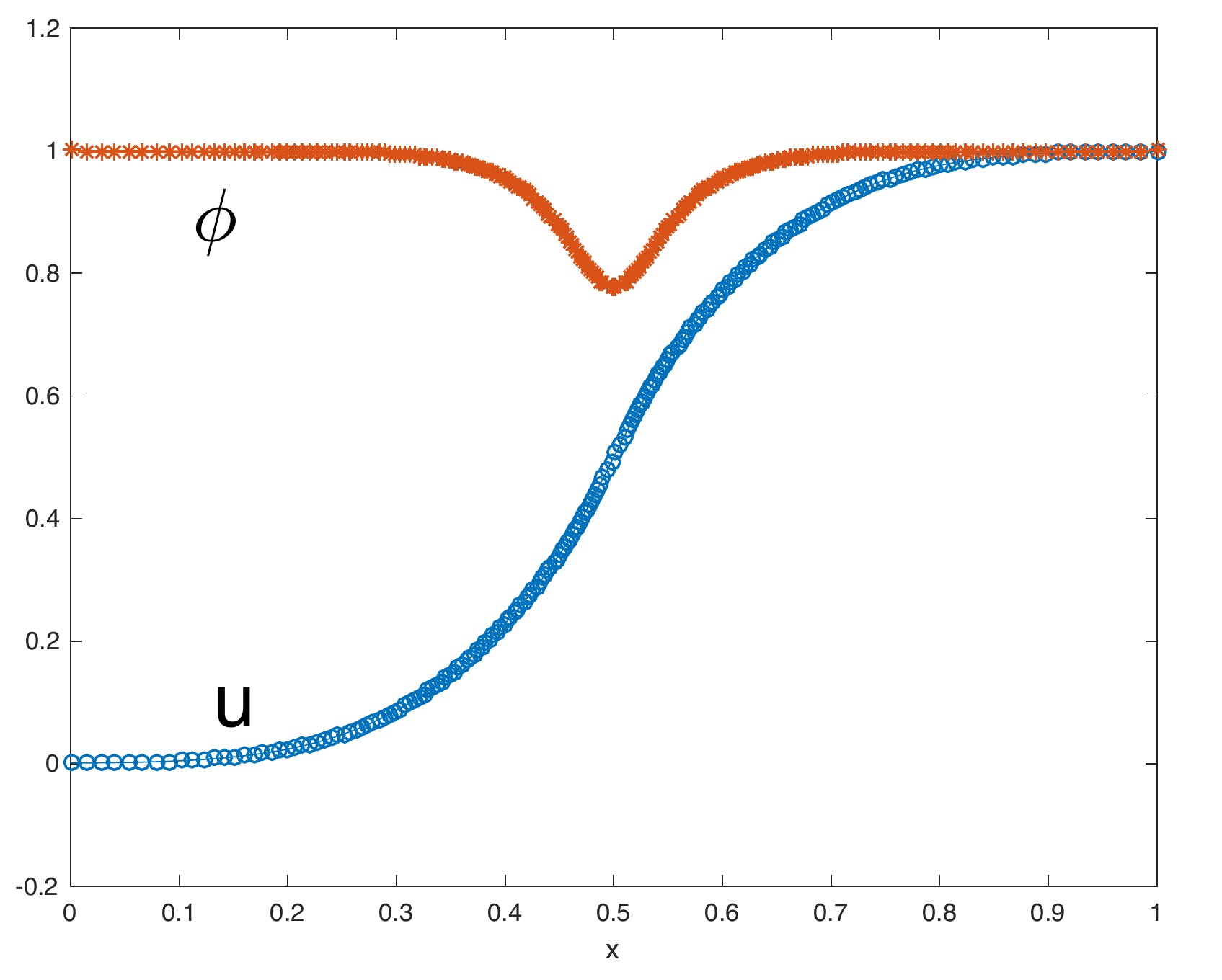}
\caption{ $t = 1$}
\end{subfigure}
\begin{subfigure}{0.32\textwidth}
\centering
\includegraphics[scale = 0.24]{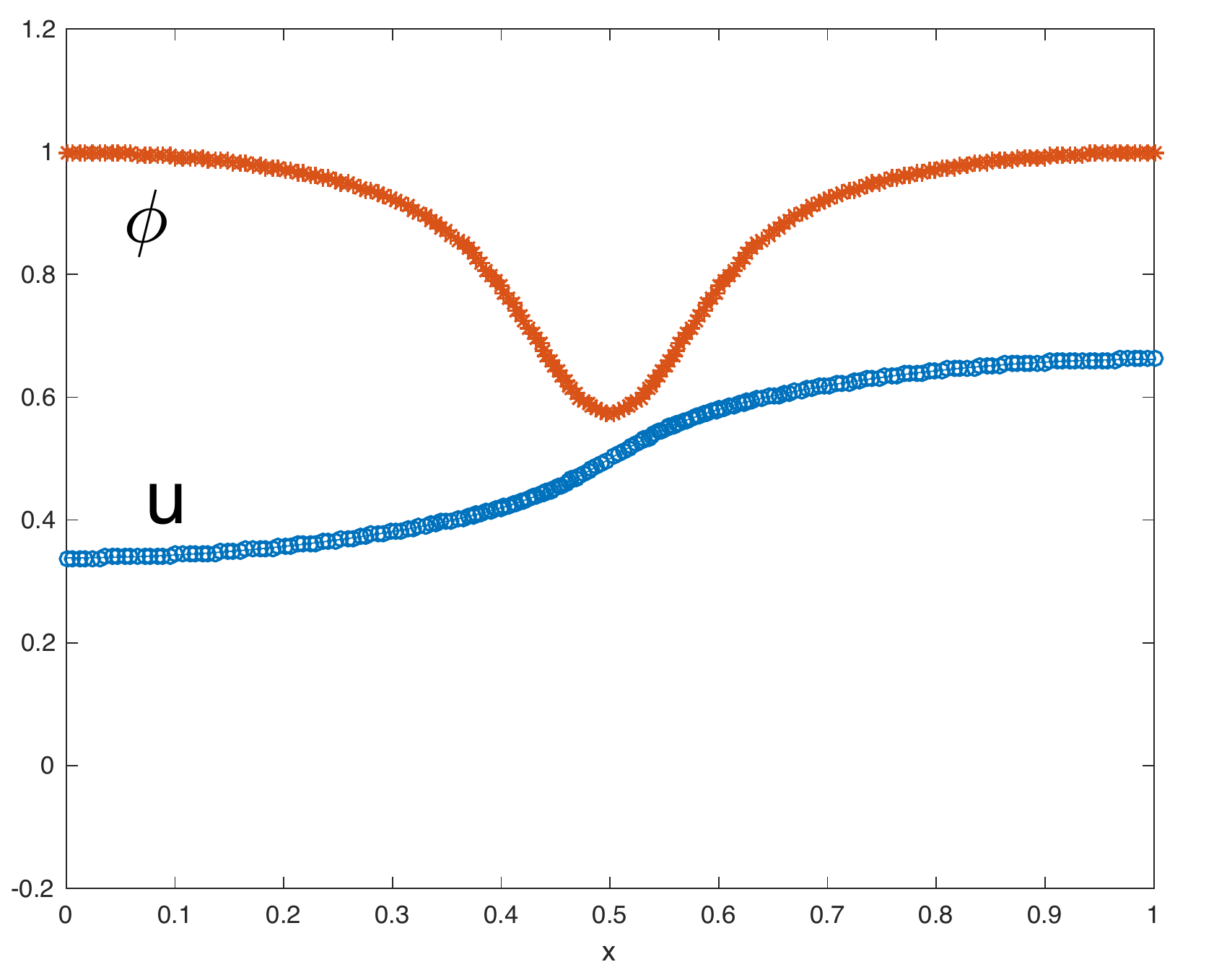}
\caption{ $t = 20$}
\end{subfigure}
\begin{subfigure}{0.32\textwidth}
\centering
\includegraphics[scale = 0.24]{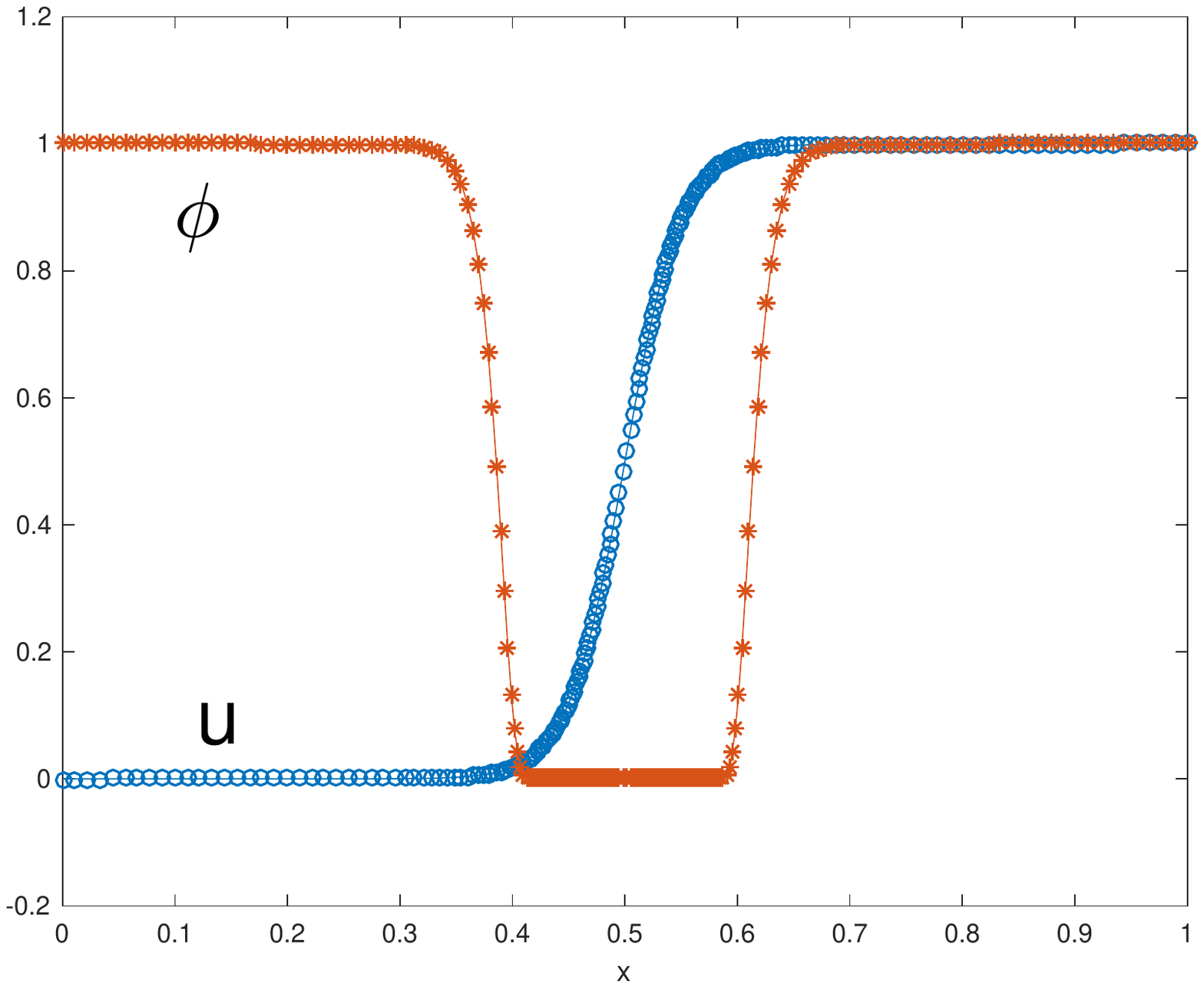}
\caption{ $t = 0.005$}
\end{subfigure}
\begin{subfigure}{0.32\textwidth}
\centering
\includegraphics[scale = 0.24]{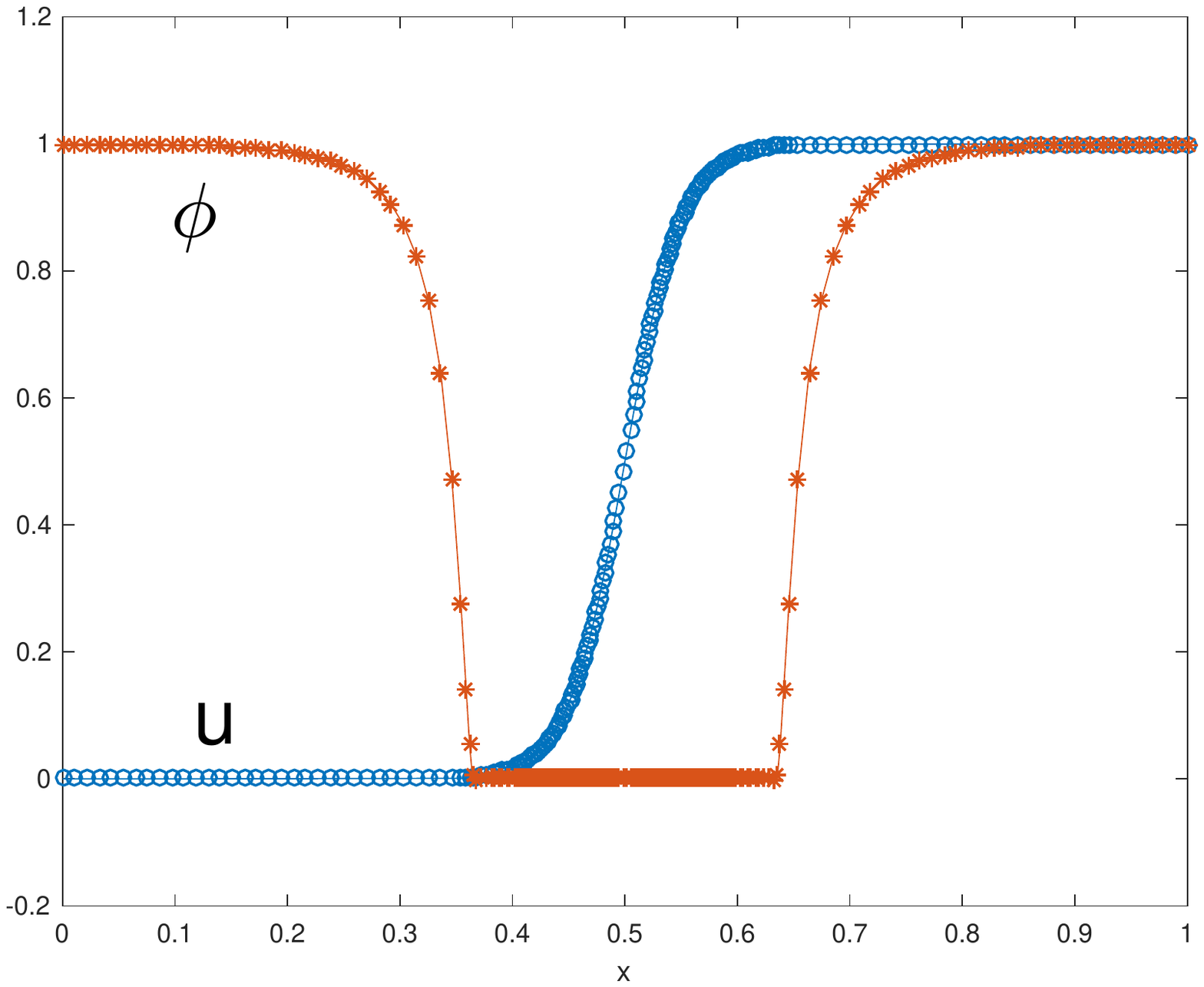}
\caption{ $t = 1$}
\end{subfigure}
 \begin{subfigure}{0.32\textwidth}
 \centering
\includegraphics[scale = 0.24]{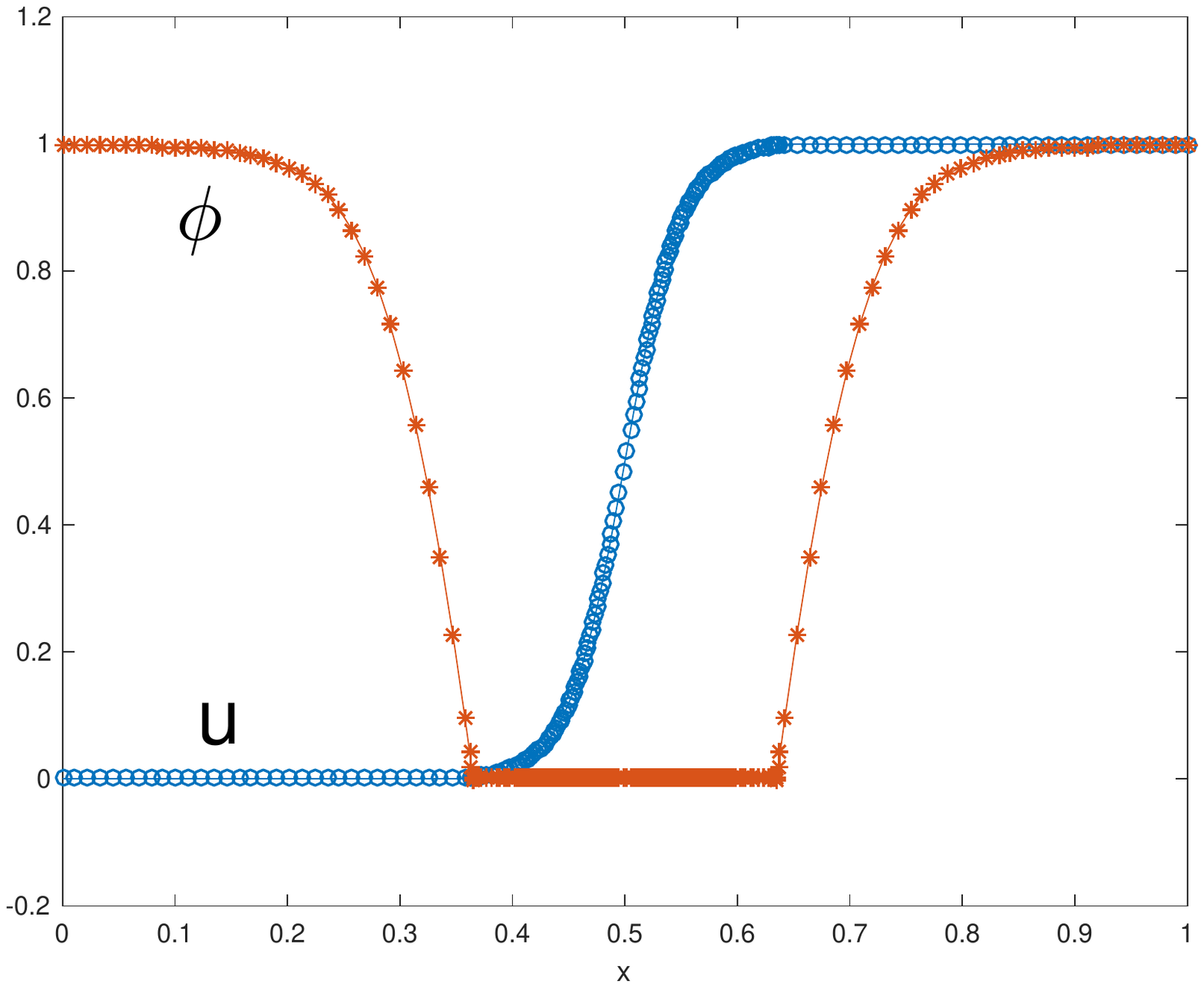}
\caption{ $t = 20$}
\end{subfigure}
\caption{Example~\ref{ex41} with $g = 0.5(1+\tanh (20(x-0.5)))$ and $u^{0} = g$. $\epsilon$ is chosen as in
(\ref{epsilon-1}) and other parameters are the same as in Example~\ref{ex41}. No scaling is used for the top row
while the scaling with (\ref{L-1}) for $u$ and $g$ is used for the bottom row.}
\label{rescaling}
\end{figure}

To demonstrate the effects of the scaling, we recompute Example~\ref{ex41} with
$u^0 = g = 0.5( 1 + \tanh(20(x-0.5)))$, which has a less steep jump at $x=0.5$ than the function (\ref{ex41-1}).
Results with and without scaling are shown in Fig.~\ref{rescaling}. It can be seen that
scaling improves the segmentation ability of the Ambrosio-Tortorelli functional.

\subsection{Segmentation for real images}

\begin{figure}[htb]
\centering
\begin{subfigure}{0.32\textwidth}
\centering
\includegraphics[scale = 0.45]{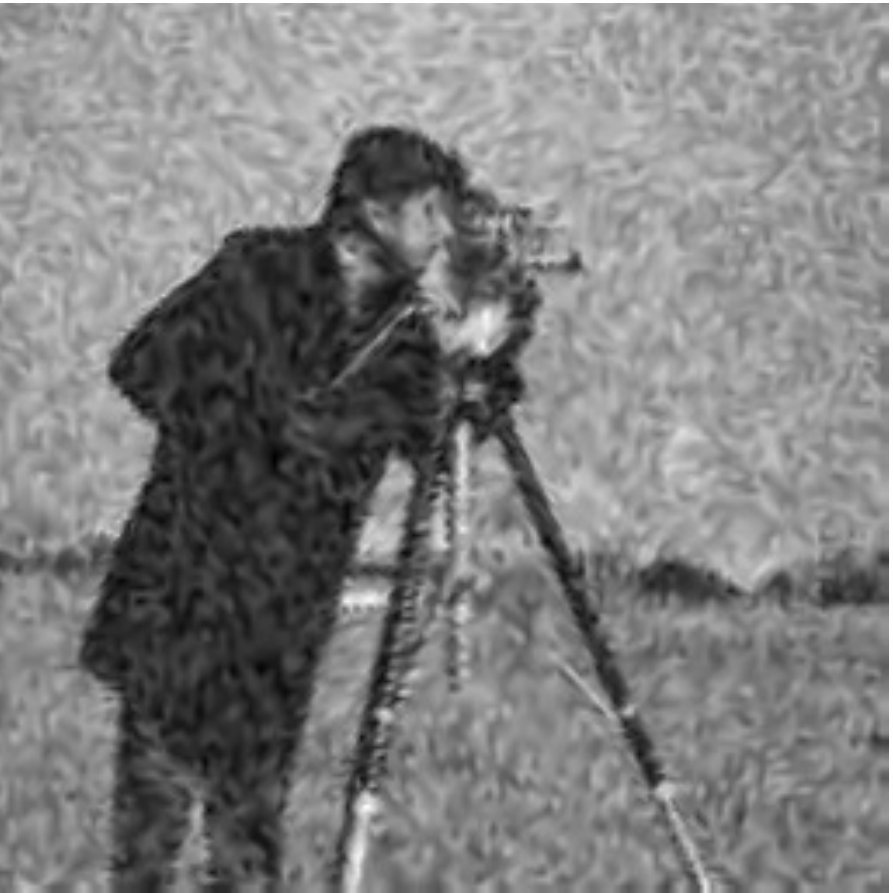}
\caption{ $t = 0.005$, $\epsilon = 10^{-5}$}
\end{subfigure}
\begin{subfigure}{0.32\textwidth}
\centering
\includegraphics[scale = 0.45]{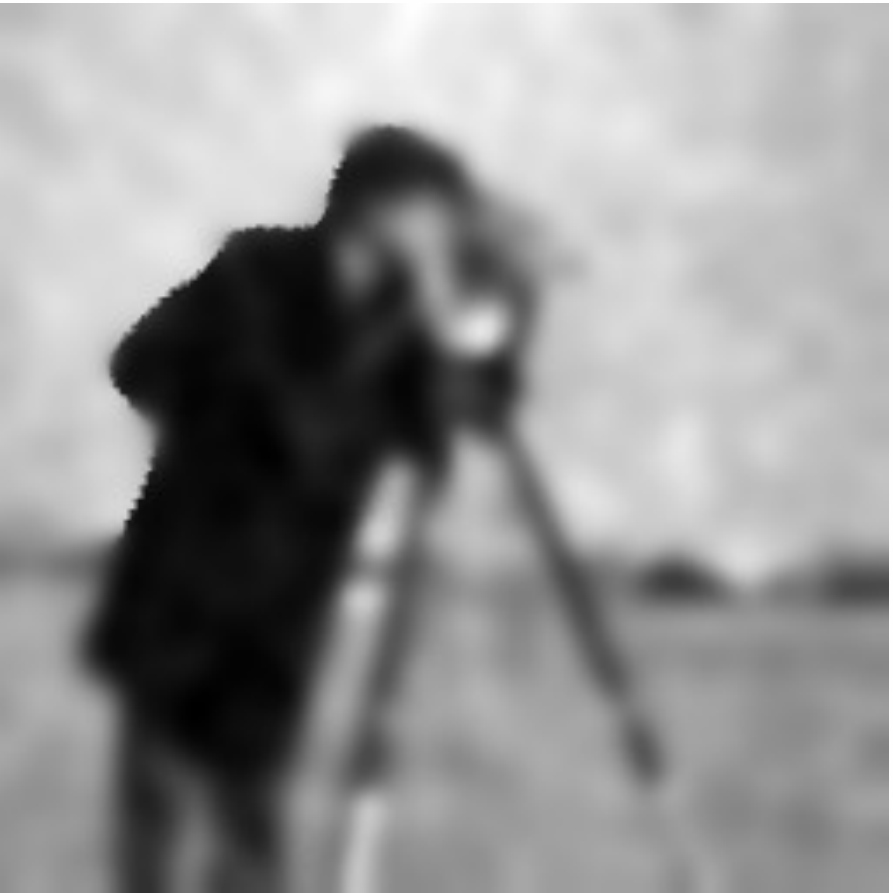}
\caption{ $t = 0.16$, $\epsilon = 10^{-5}$}
\end{subfigure}
 \begin{subfigure}{0.32\textwidth}
 \centering
\includegraphics[scale = 0.45]{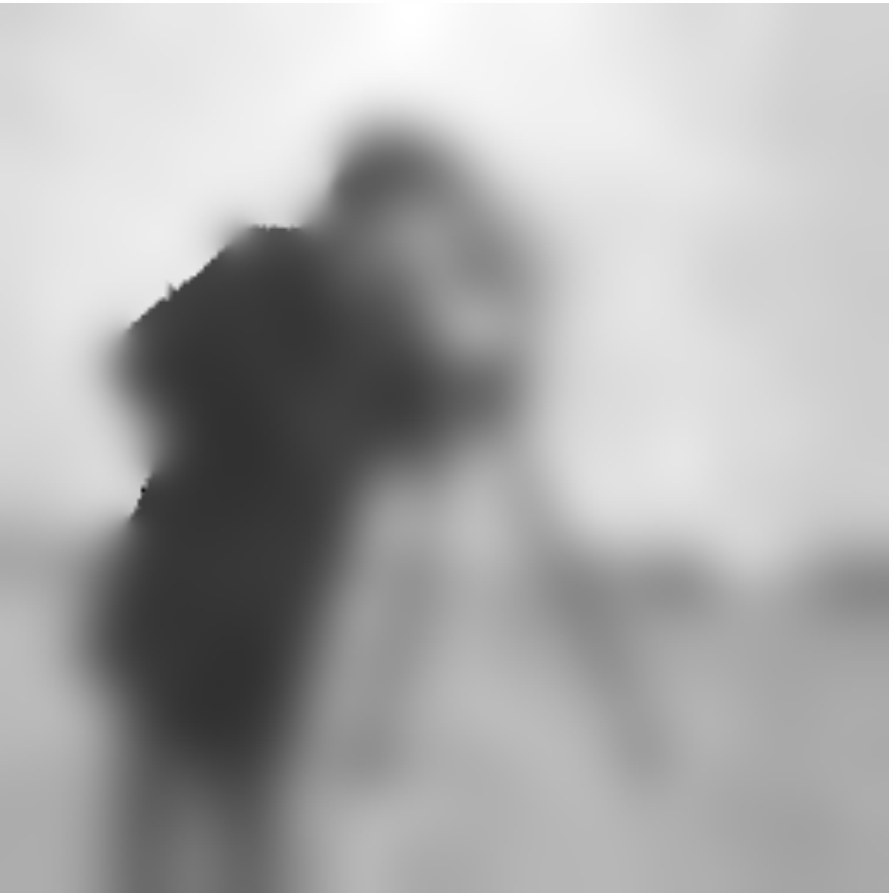}
\caption{$t = 0.6$, $\epsilon = 10^{-5}$}
\end{subfigure}
\begin{subfigure}{0.32\textwidth}
\centering
\includegraphics[scale = 0.45]{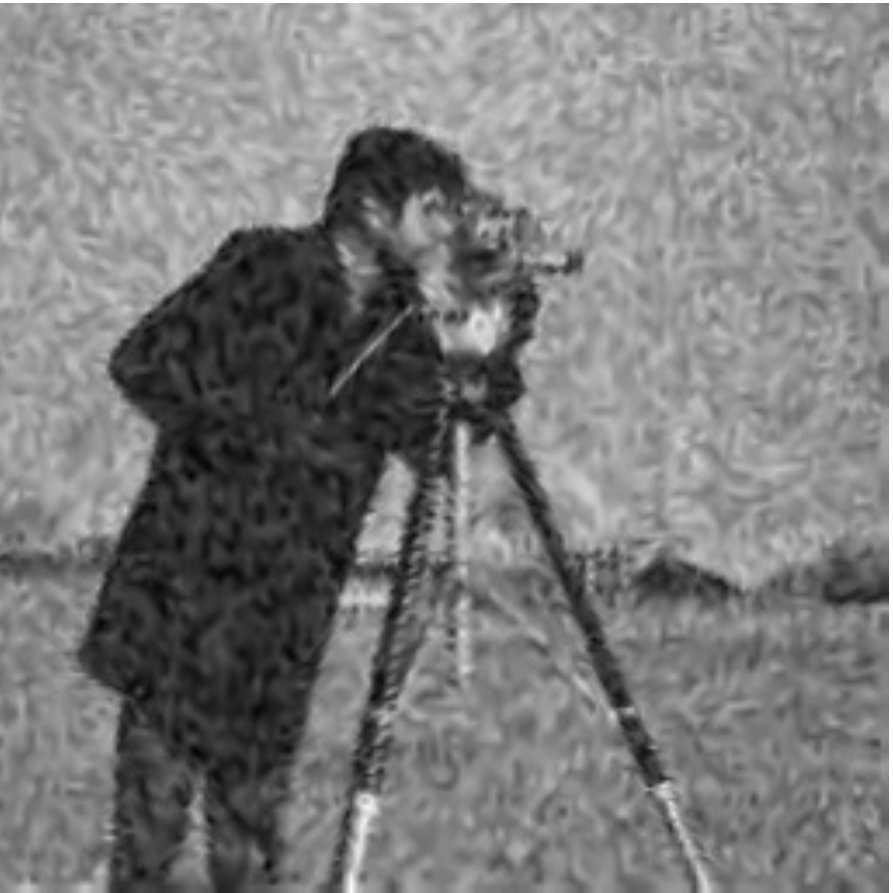}
\caption{$t = 0.005$, $\epsilon$ is chosen by (\ref{epsilon-1})}
\end{subfigure}
\begin{subfigure}{0.32\textwidth}
\centering
\includegraphics[scale = 0.45]{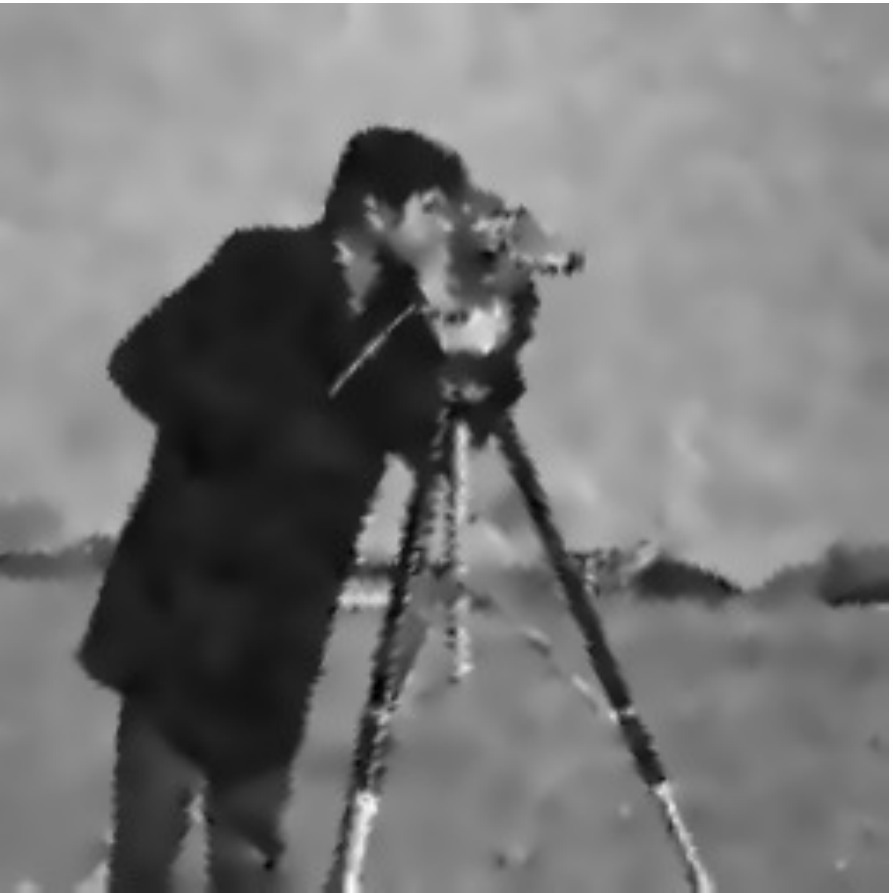}
\caption{ $t = 0.16$, $\epsilon$ is chosen by (\ref{epsilon-1})}
\end{subfigure}
 \begin{subfigure}{0.32\textwidth}
 \centering
\includegraphics[scale = 0.45]{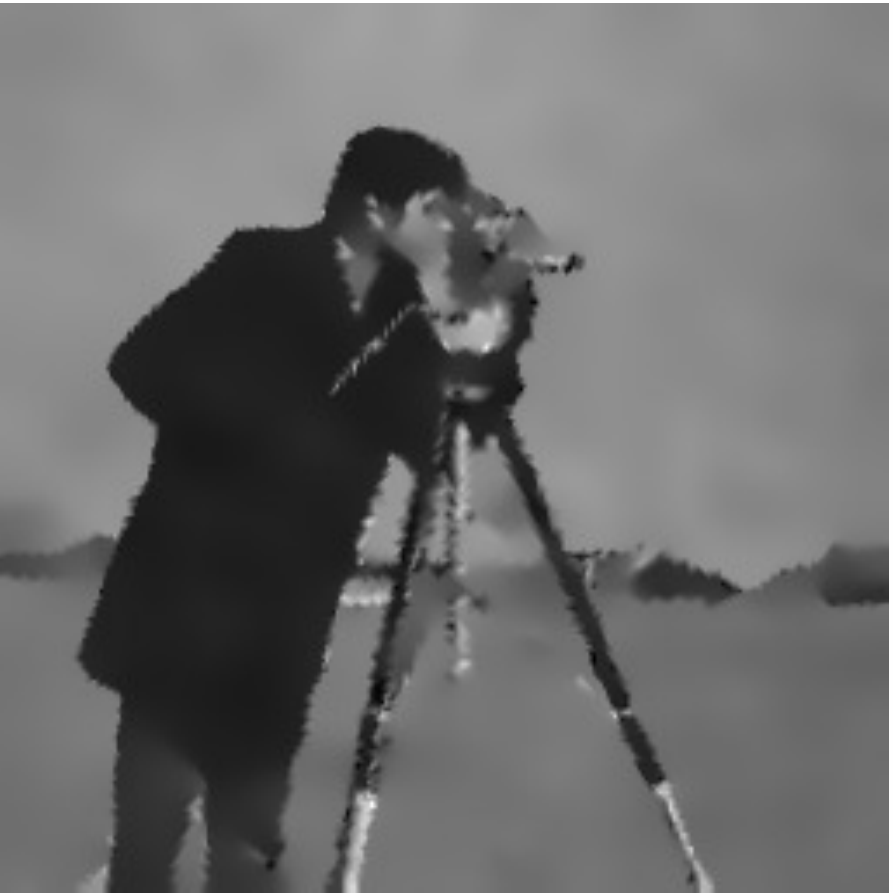}
\caption{$t = 0.6$, $\epsilon$ is chosen by (\ref{epsilon-1})}
\end{subfigure}
\caption{A comparison of the image segmentation with different $\epsilon$ values.}
\label{cameraman}
\end{figure}

\begin{figure}[htb]
\centering
\begin{subfigure}{0.32\textwidth}
\centering
\includegraphics[scale = 0.3]{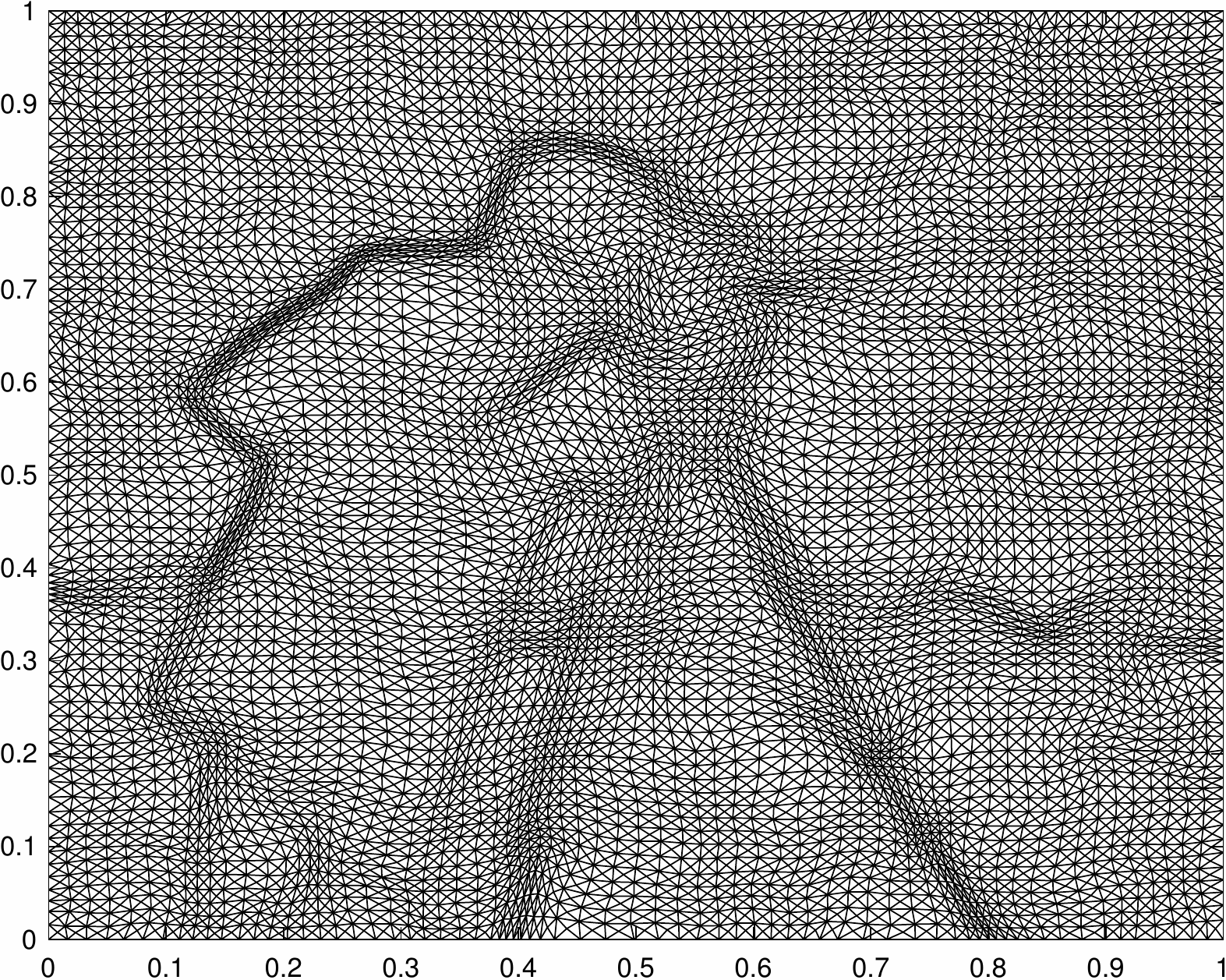}
\caption{$t = 0.005$, $\epsilon = 10^{-5}$}
\end{subfigure}
\begin{subfigure}{0.32\textwidth}
\centering
\includegraphics[scale = 0.3]{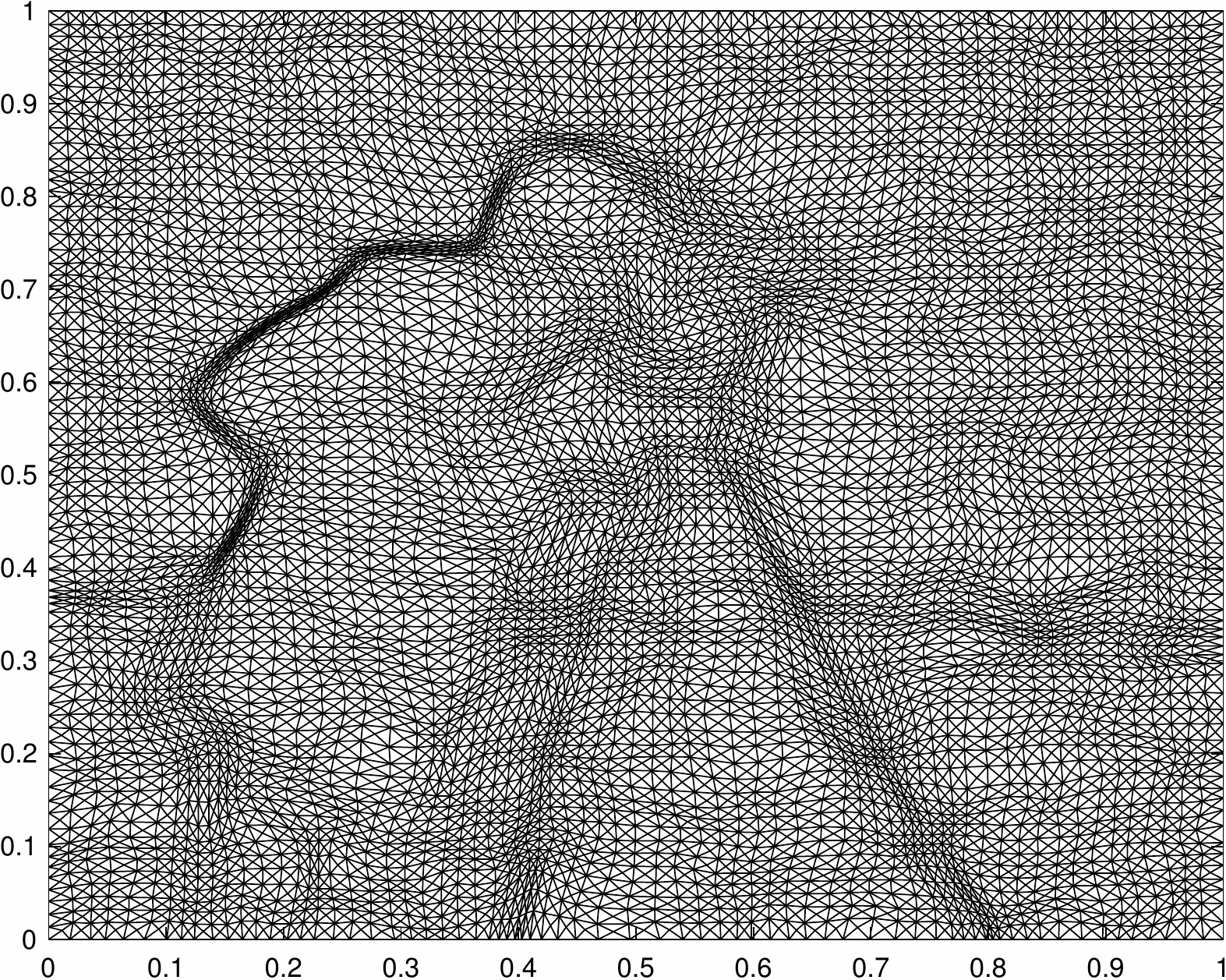}
\caption{$t = 0.6$, $\epsilon = 10^{-5}$}
\end{subfigure}
 \begin{subfigure}{0.32\textwidth}
 \centering
\includegraphics[scale = 0.3]{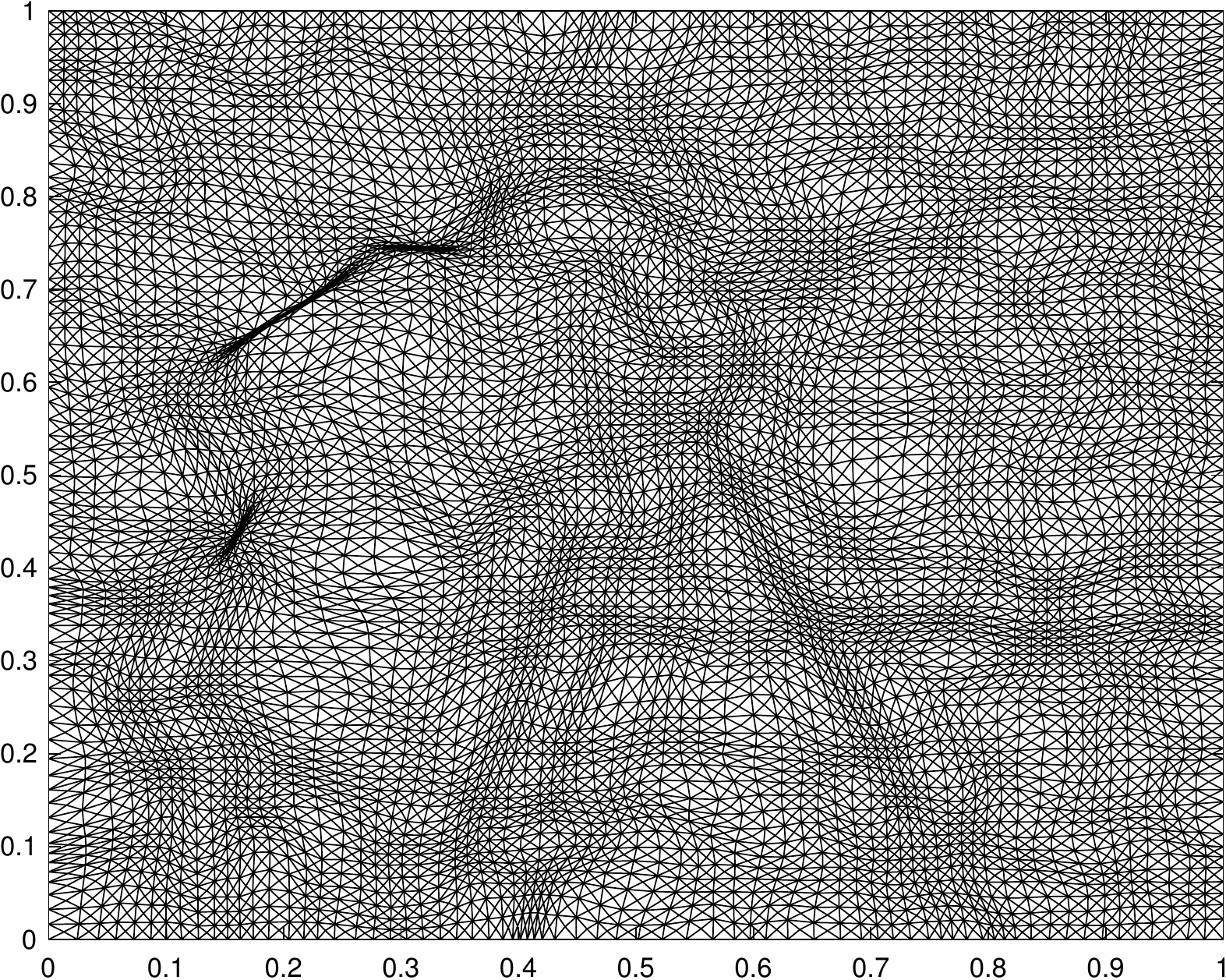}
\caption{$t = 0.6$, $\epsilon = 10^{-5}$}
\end{subfigure}
\begin{subfigure}{0.32\textwidth}
\centering
\includegraphics[scale = 0.3]{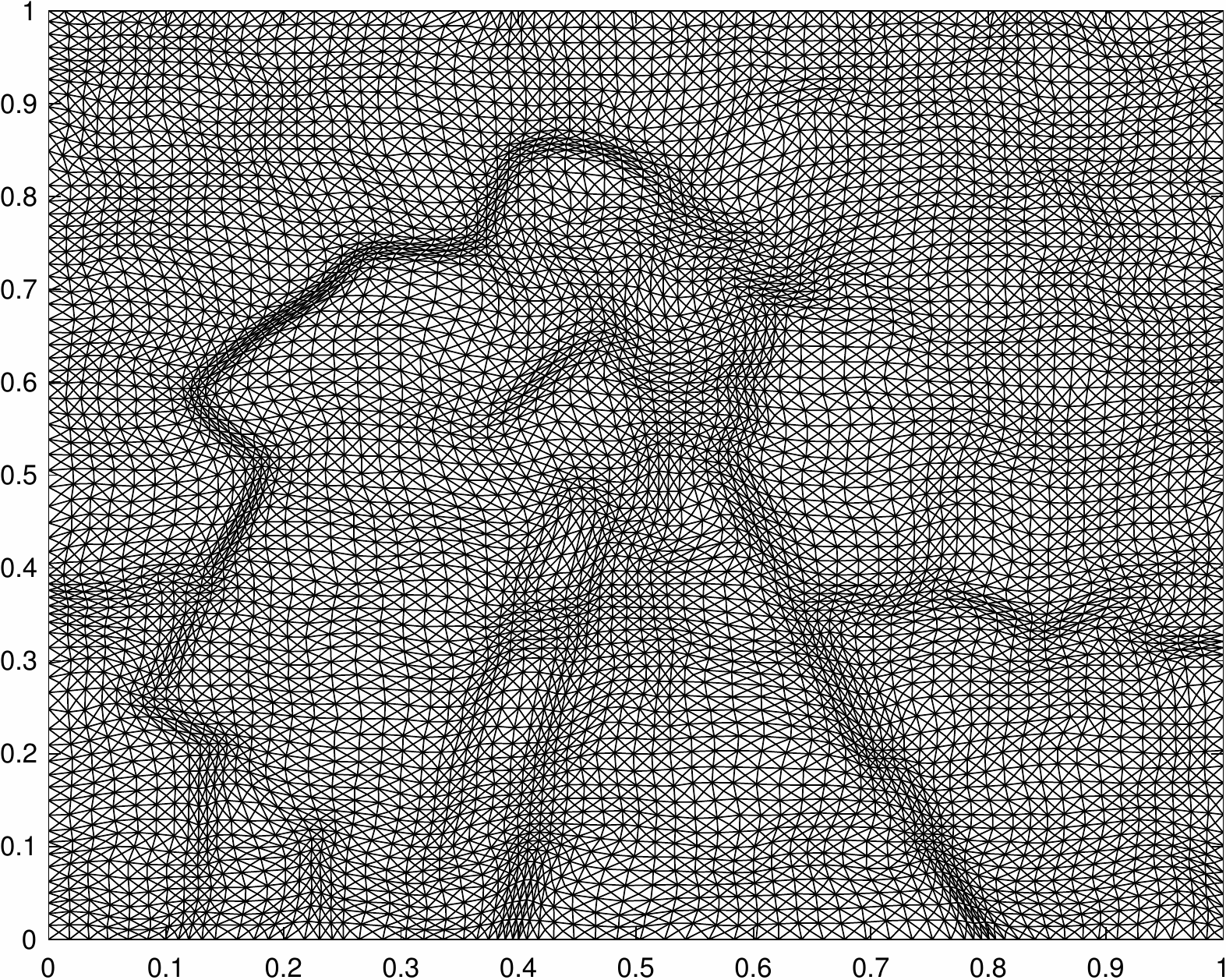}
\caption{$t = 0.005$, $\epsilon$ is chosen by (\ref{epsilon-1})}
\end{subfigure}
\begin{subfigure}{0.32\textwidth}
\centering
\includegraphics[scale = 0.3]{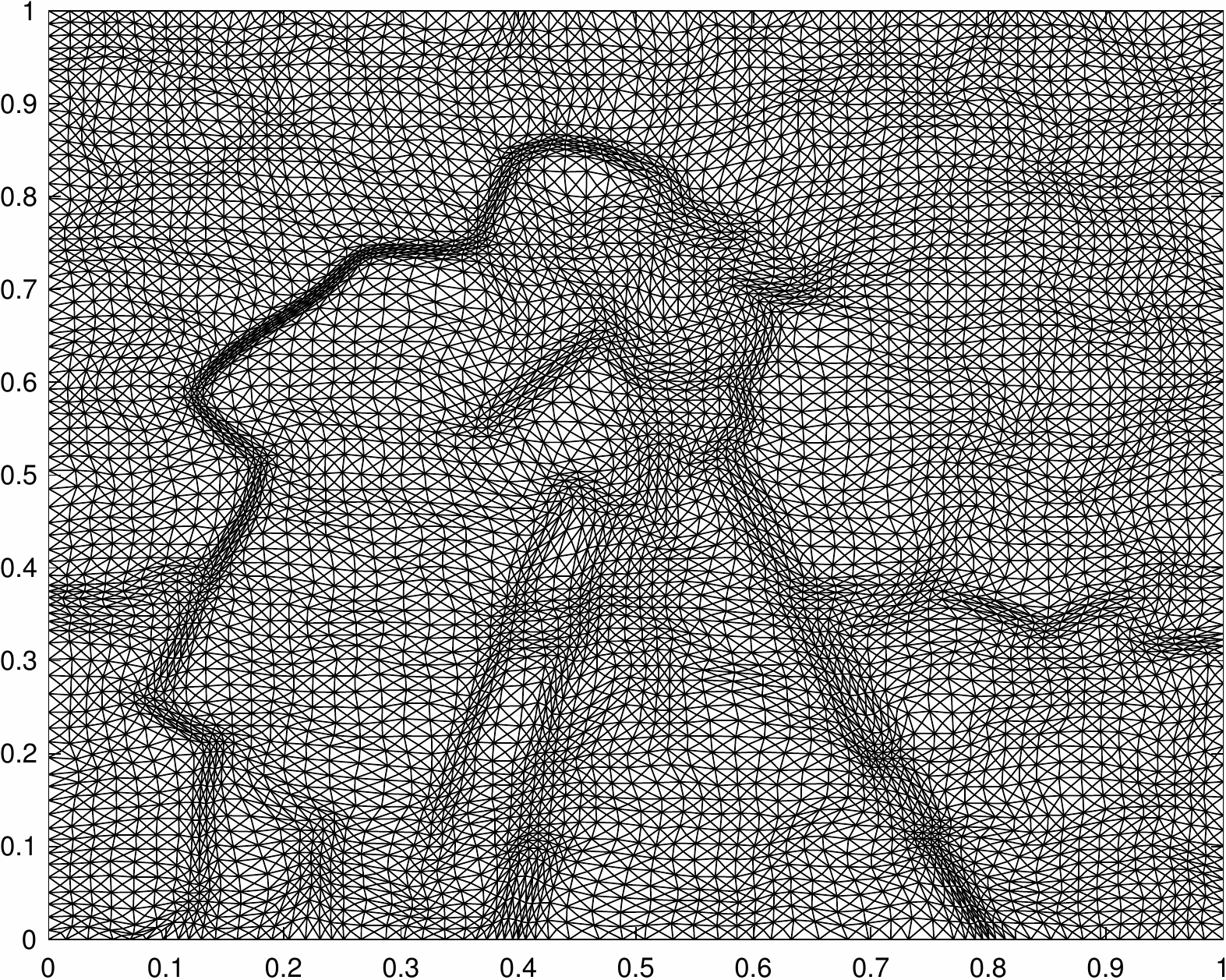}
\caption{$t = 0.16$, $\epsilon$ is chosen by (\ref{epsilon-1})}
\end{subfigure}
 \begin{subfigure}{0.32\textwidth}
 \centering
\includegraphics[scale = 0.3]{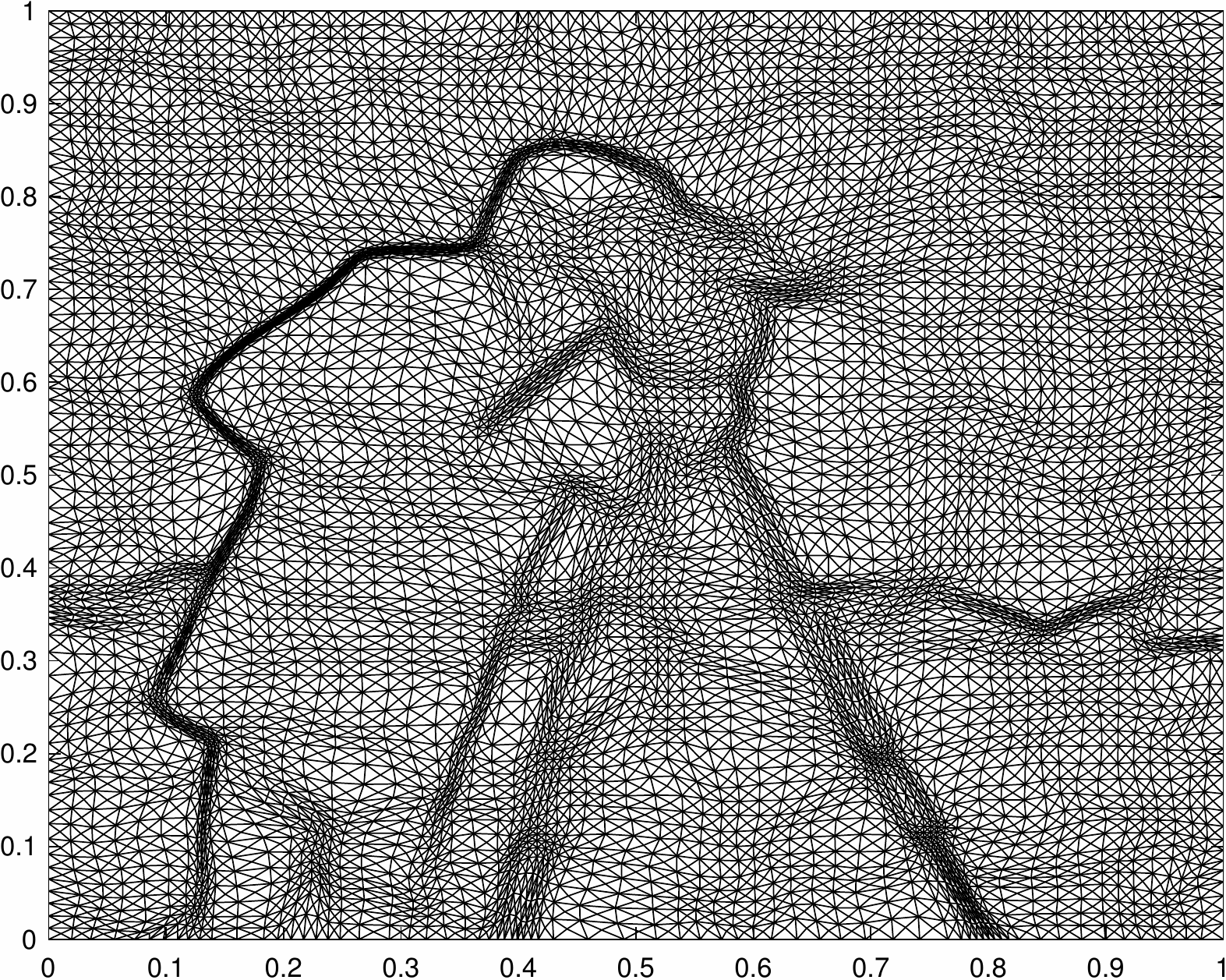}
\caption{$t = 0.6$, $\epsilon$ is chosen by (\ref{epsilon-1})}
\end{subfigure}
\caption{The meshes corresponding to Fig.~\ref{cameraman}.}
\label{cameramanmesh}
\end{figure}

\begin{figure}[htb]
\centering
\begin{subfigure}{0.34\textwidth}
\centering
\includegraphics[scale = 0.24]{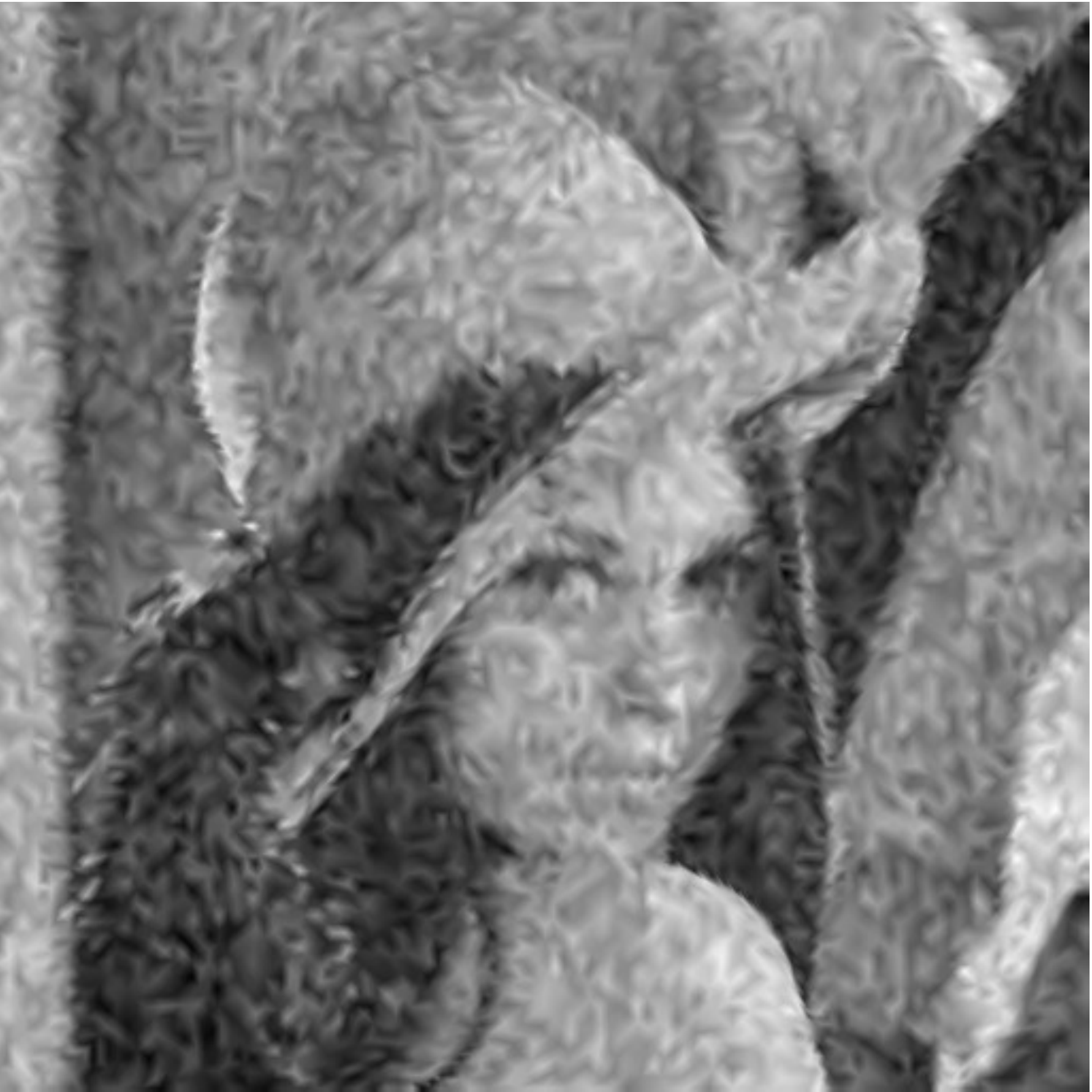}
\caption{$t = 6 \times 10^{-6}$, $\epsilon = 10^{-6}$}
\end{subfigure}
\begin{subfigure}{0.32\textwidth}
\centering
\includegraphics[scale = 0.24]{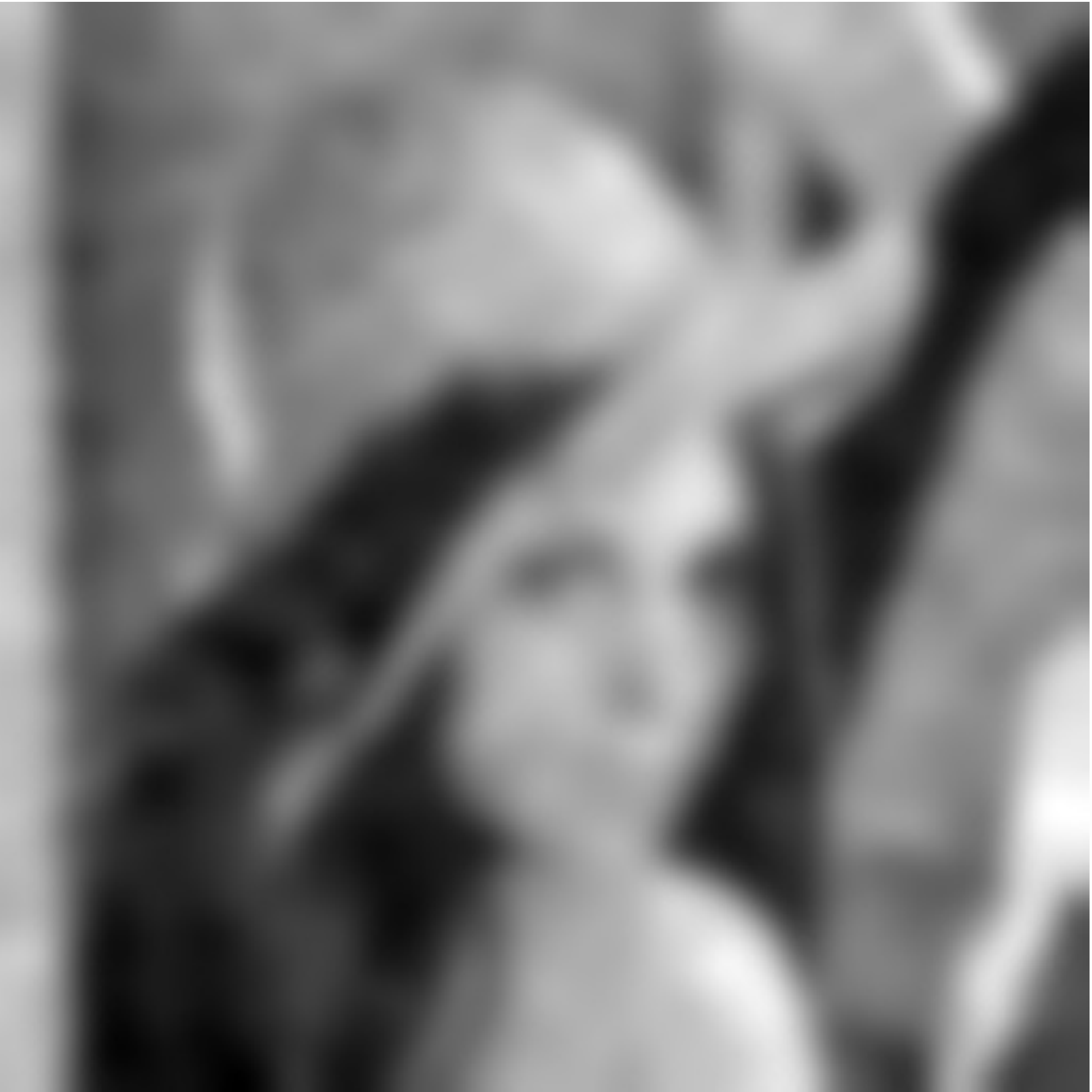}
\caption{ $t = 0.15$, $\epsilon = 10^{-6}$}
\end{subfigure}
 \begin{subfigure}{0.32\textwidth}
 \centering
\includegraphics[scale = 0.24]{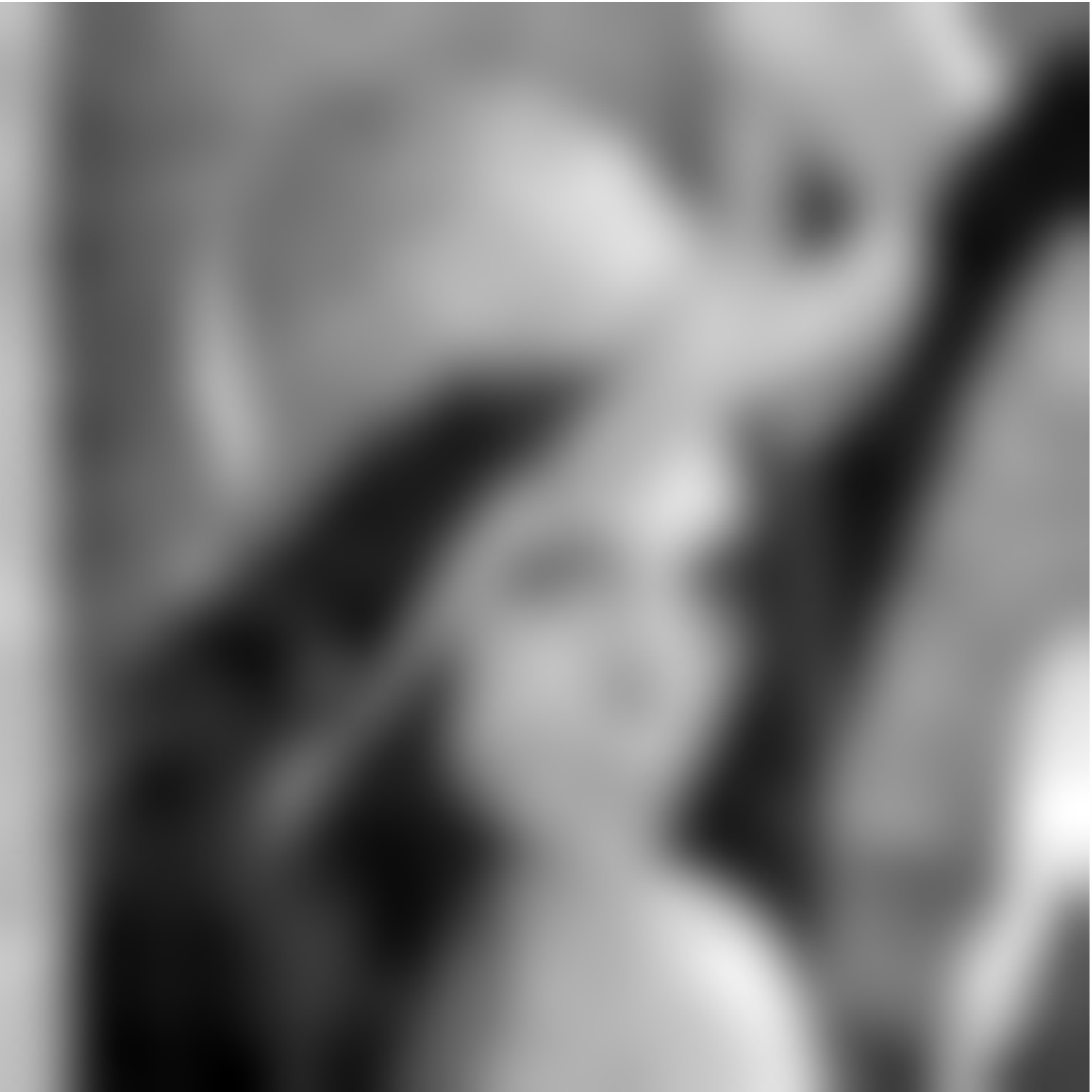}
\caption{ $t = 0.3$, $\epsilon = 10^{-6}$}
\end{subfigure}
\begin{subfigure}{0.34\textwidth}
\centering
\includegraphics[scale = 0.24]{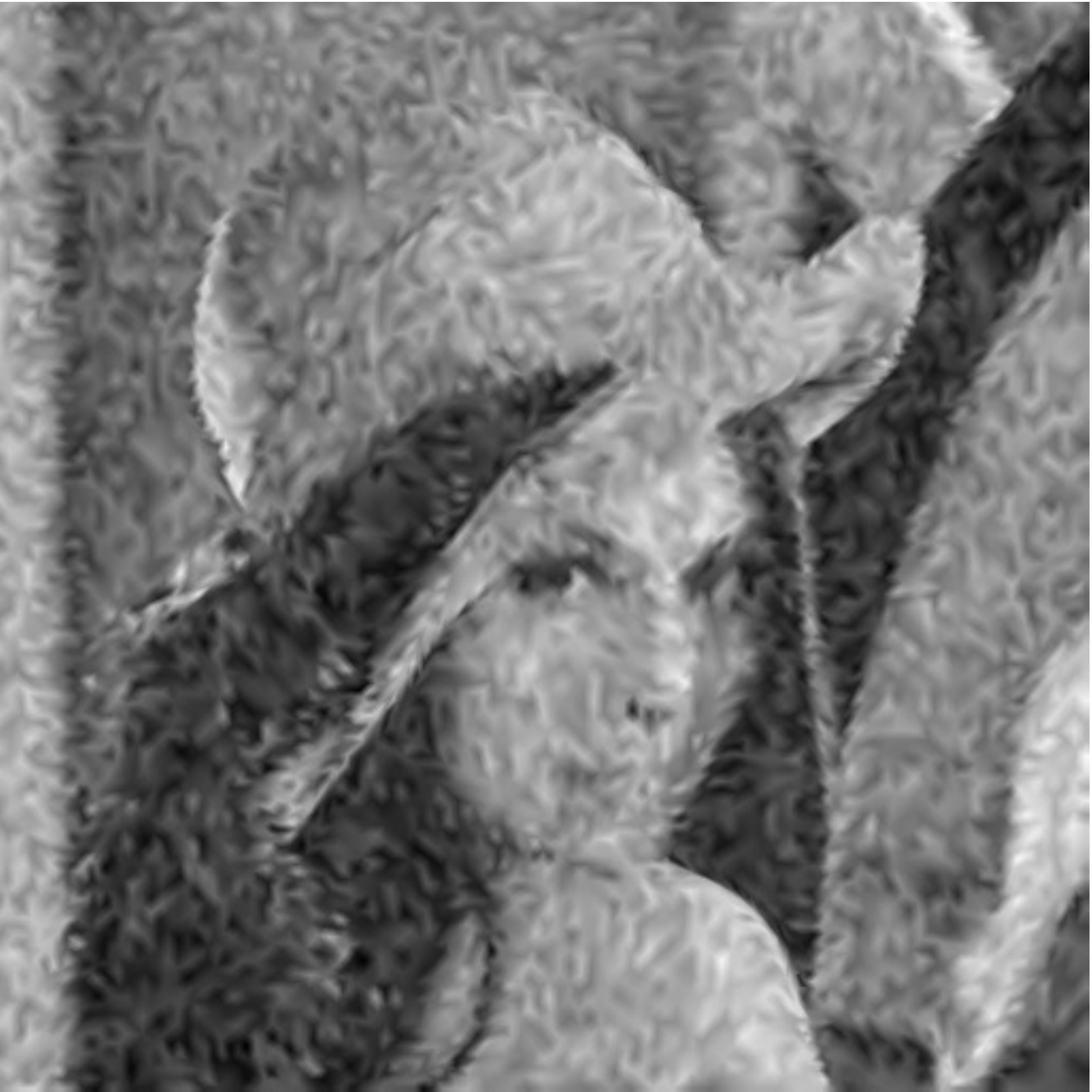}
\caption{$ t = 6 \times 10^{-6}$, $\epsilon$ is chosen by (\ref{epsilon-1})}
\end{subfigure}
\begin{subfigure}{0.32\textwidth}
\centering
\includegraphics[scale = 0.24]{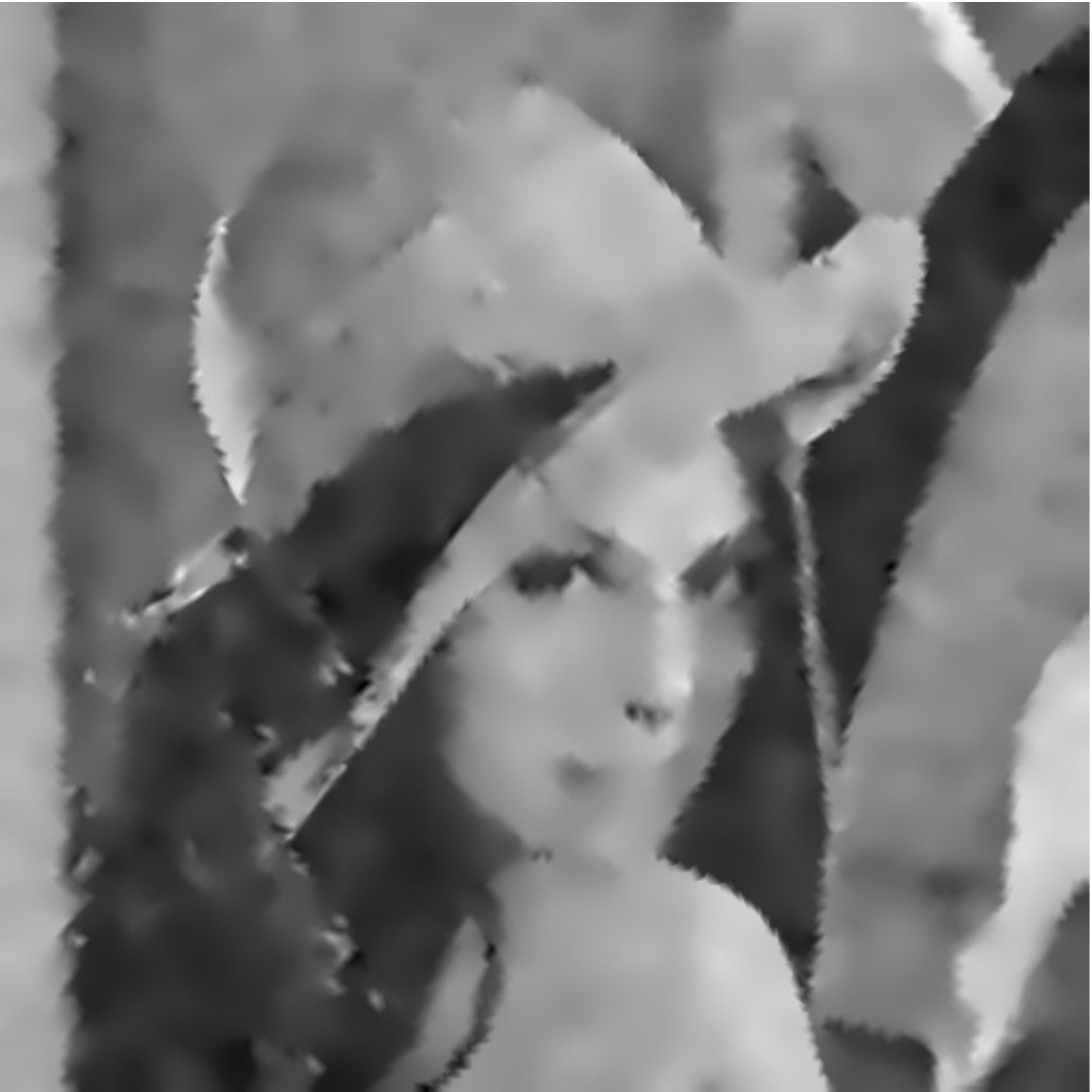}
\caption{$t = 0.15$, $\epsilon$ is chosen by (\ref{epsilon-1})}
\end{subfigure}
\begin{subfigure}{0.32\textwidth}
\centering
\includegraphics[scale = 0.24]{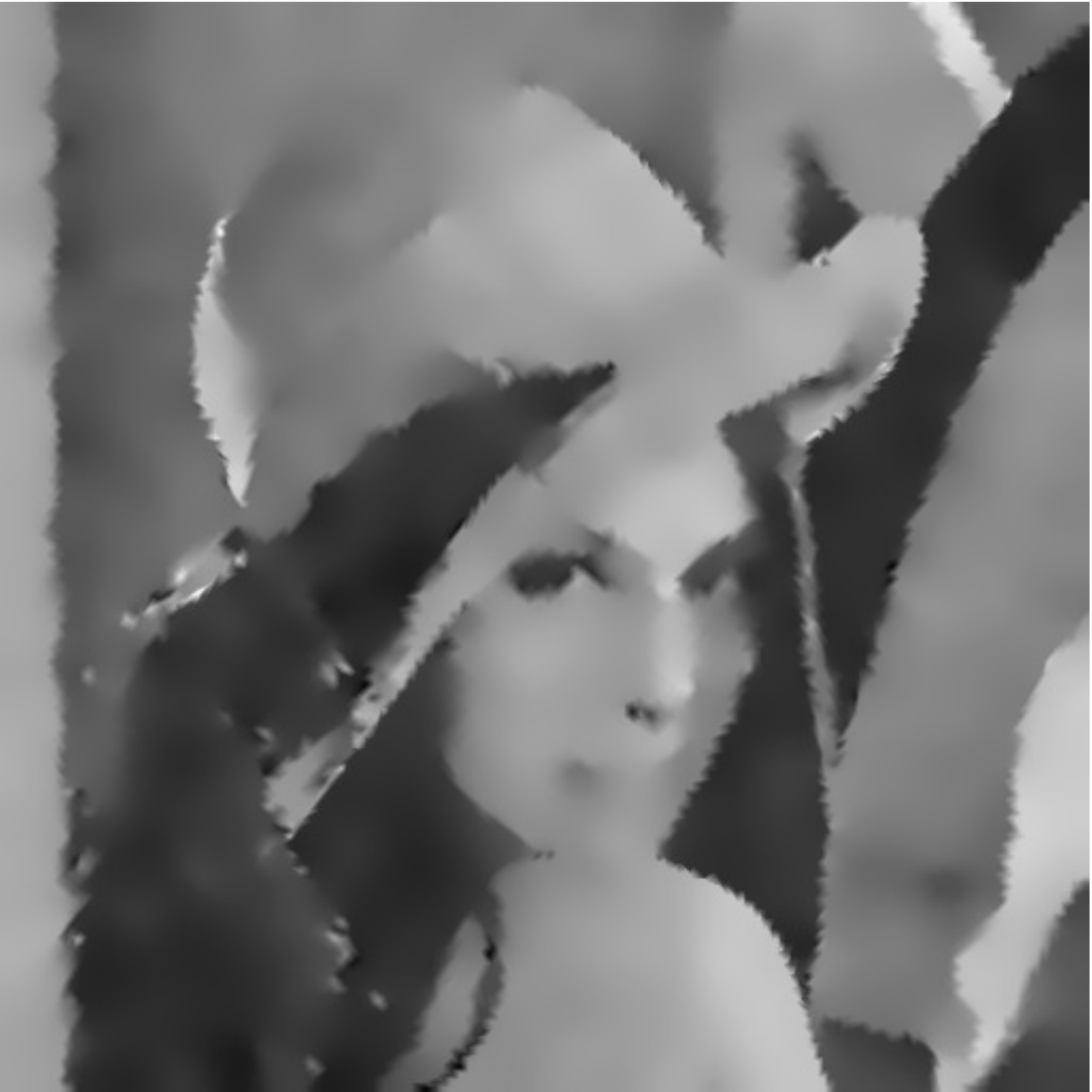}
\caption{ $t = 0.3$, $\epsilon$ is chosen by (\ref{epsilon-1})}
\end{subfigure}
\caption{Evolution of the image.}
\label{lenaimage}
\end{figure}

\begin{figure}[htb]
\centering
\begin{subfigure}{0.34\textwidth}
\centering
\includegraphics[scale = 0.3]{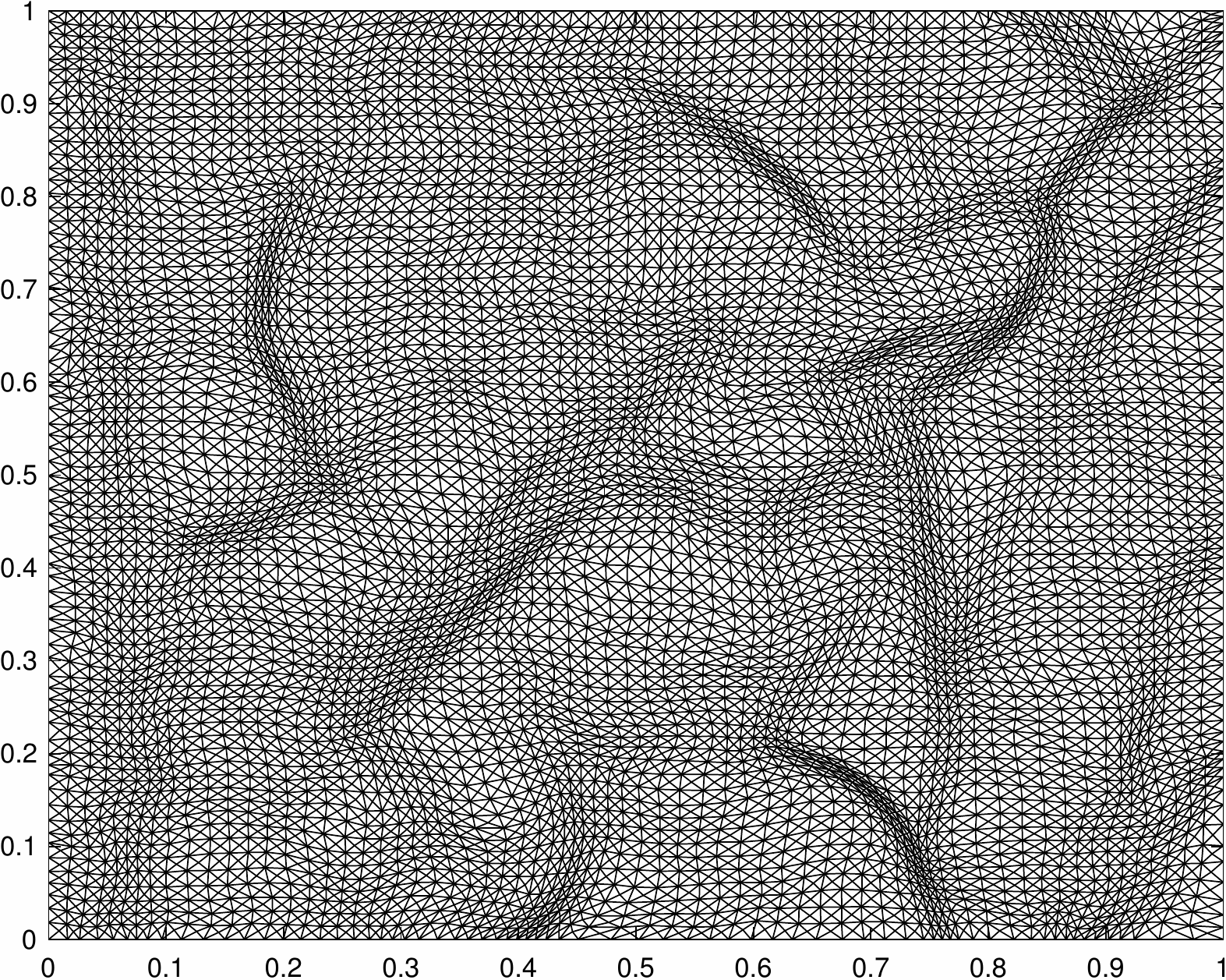}
\caption{$t = 6 \times 10^{-6}$, $\epsilon = 10^{-6}$}
\end{subfigure}
\begin{subfigure}{0.32\textwidth}
\centering
\includegraphics[scale = 0.3]{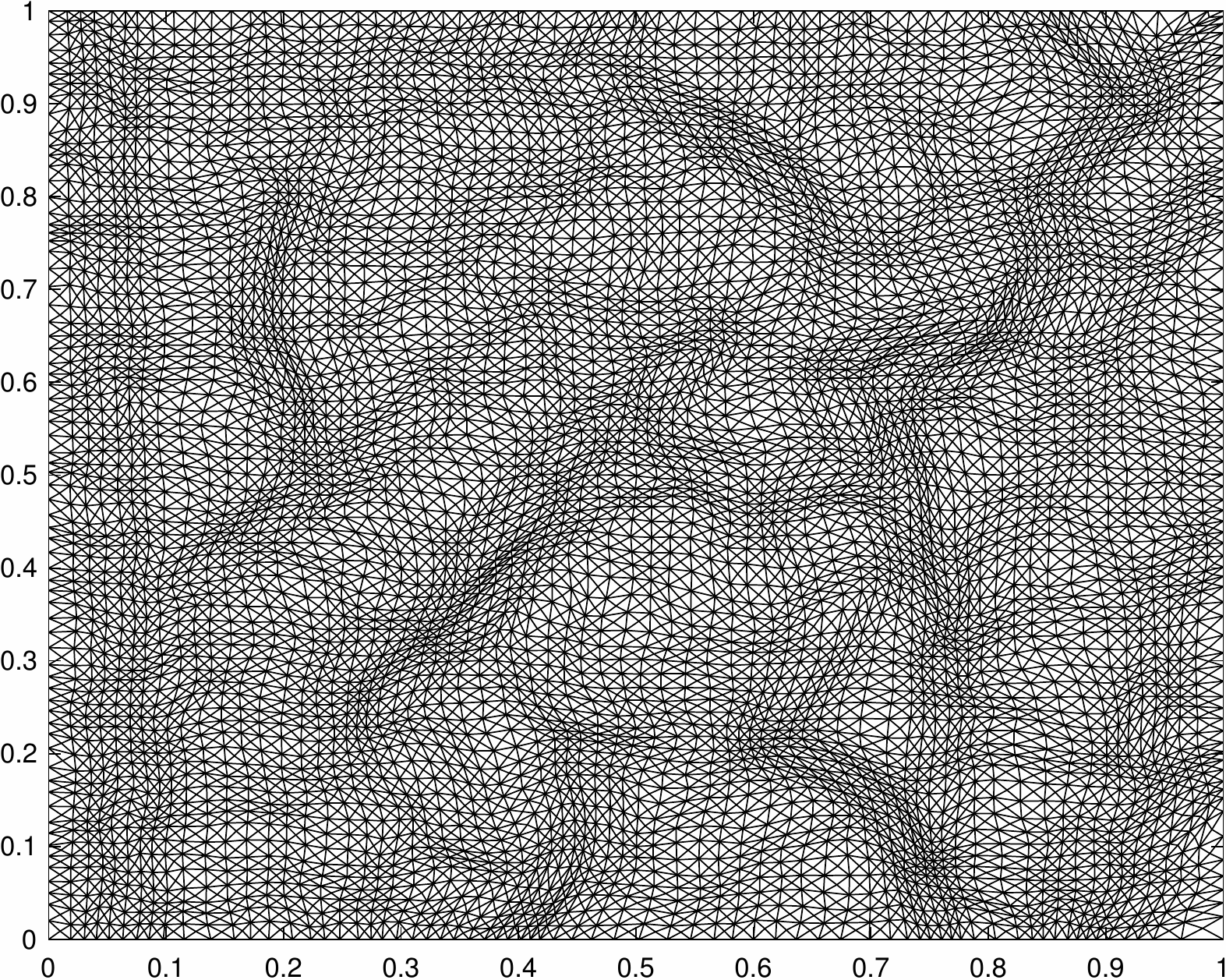}
\caption{$t = 0.15$, $\epsilon = 10^{-6}$}
\end{subfigure}
 \begin{subfigure}{0.32\textwidth}
 \centering
\includegraphics[scale = 0.3]{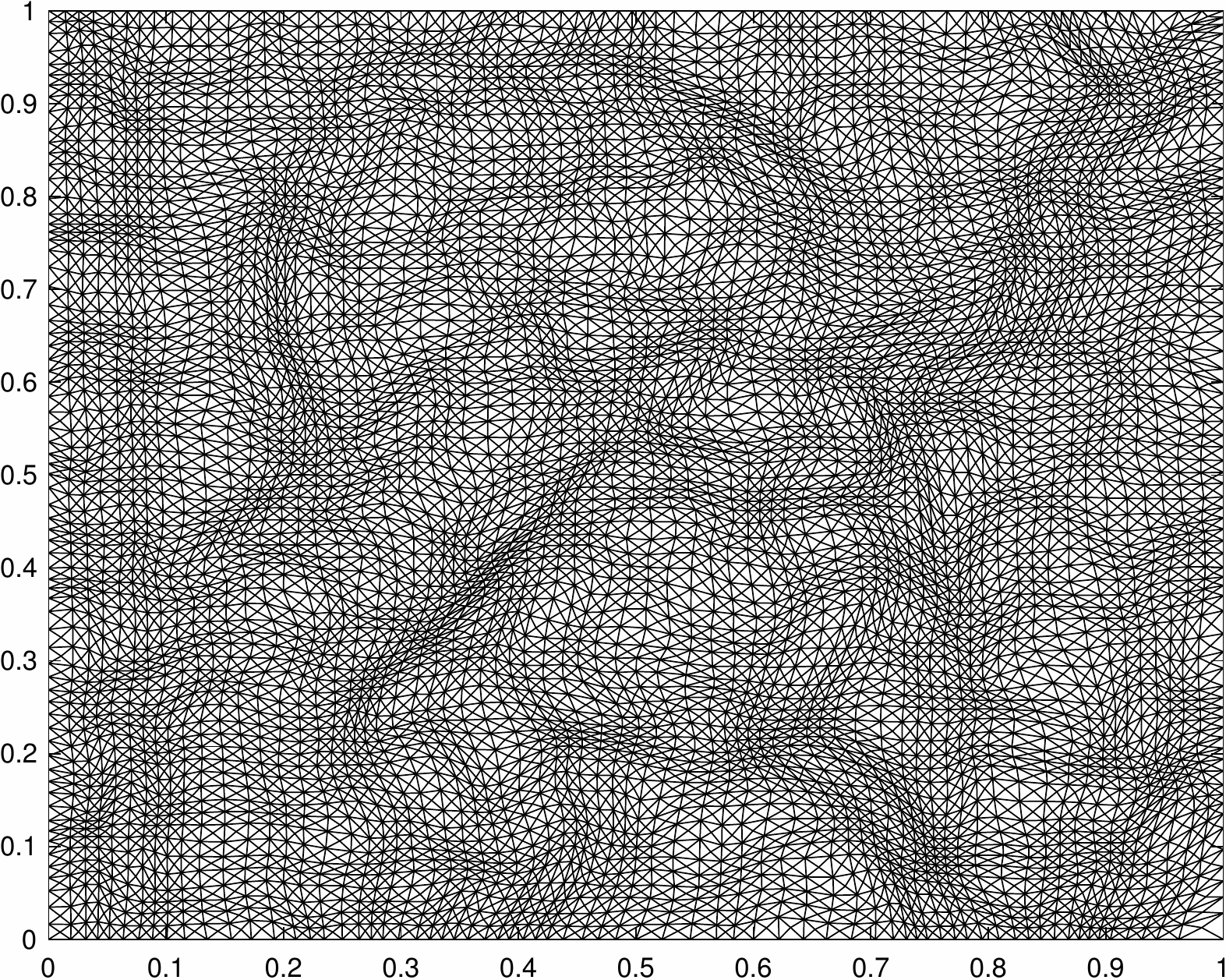}
\caption{$t = 0.3$, $\epsilon = 10^{-6}$}
\end{subfigure}
\begin{subfigure}{0.34\textwidth}
\centering
\includegraphics[scale = 0.3]{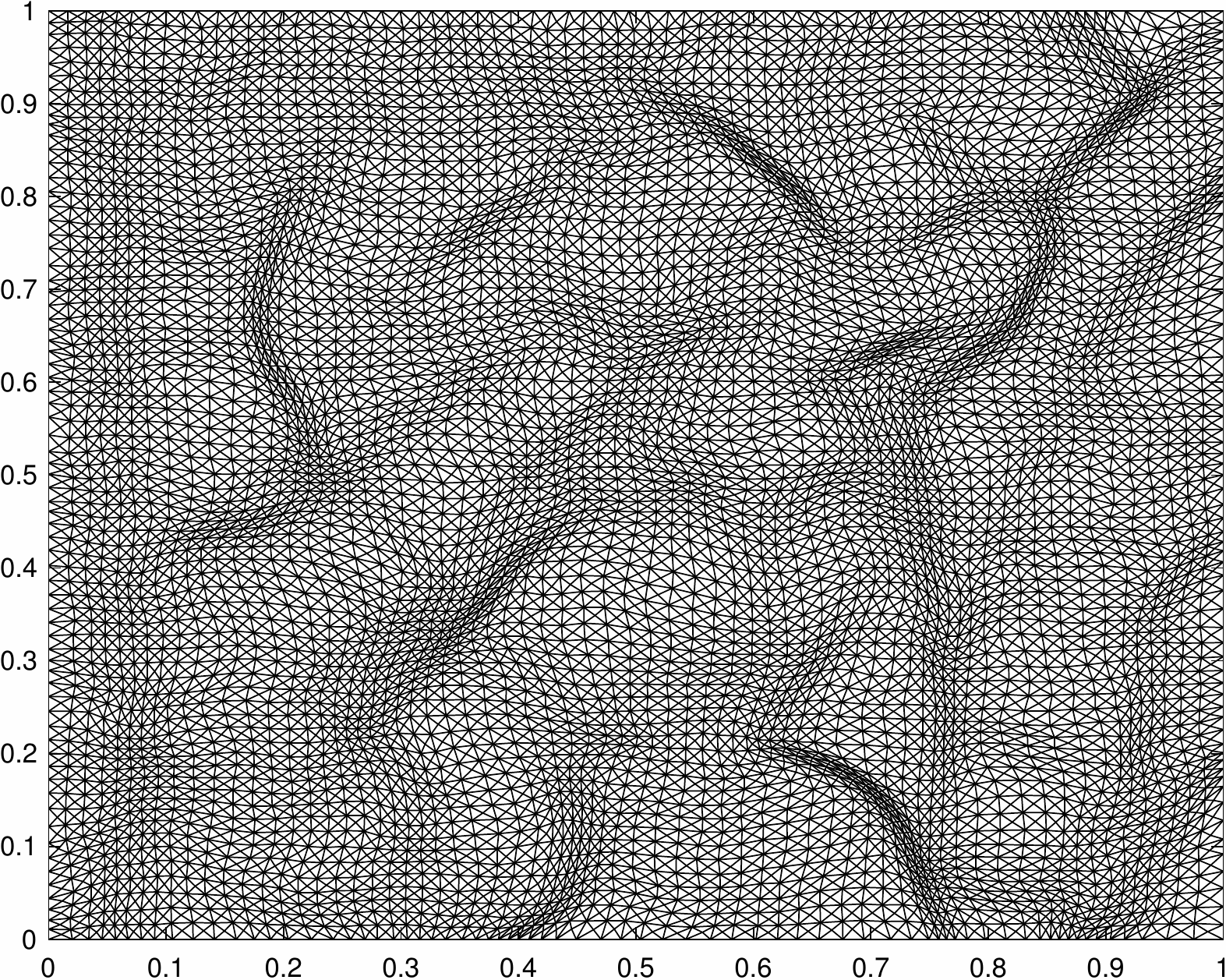}
\caption{$ t = 6 \times 10^{-6}$, $\epsilon$ is chosen by (\ref{epsilon-1})}
\end{subfigure}
\begin{subfigure}{0.32\textwidth}
\centering
\includegraphics[scale = 0.3]{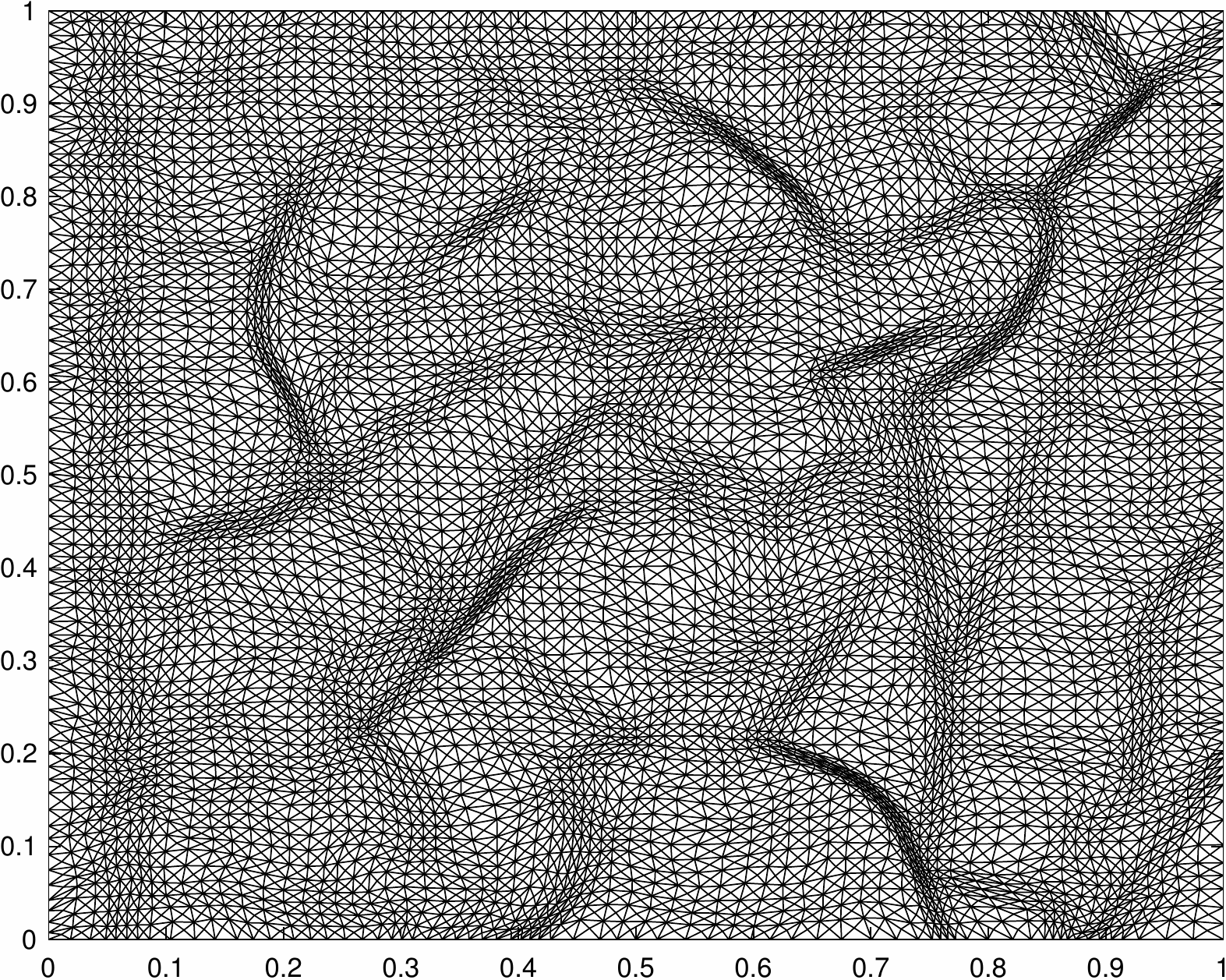}
\caption{ $t = 0.15$, $\epsilon$ is chosen by (\ref{epsilon-1})}
\end{subfigure}
 \begin{subfigure}{0.32\textwidth}
 \centering
\includegraphics[scale = 0.3]{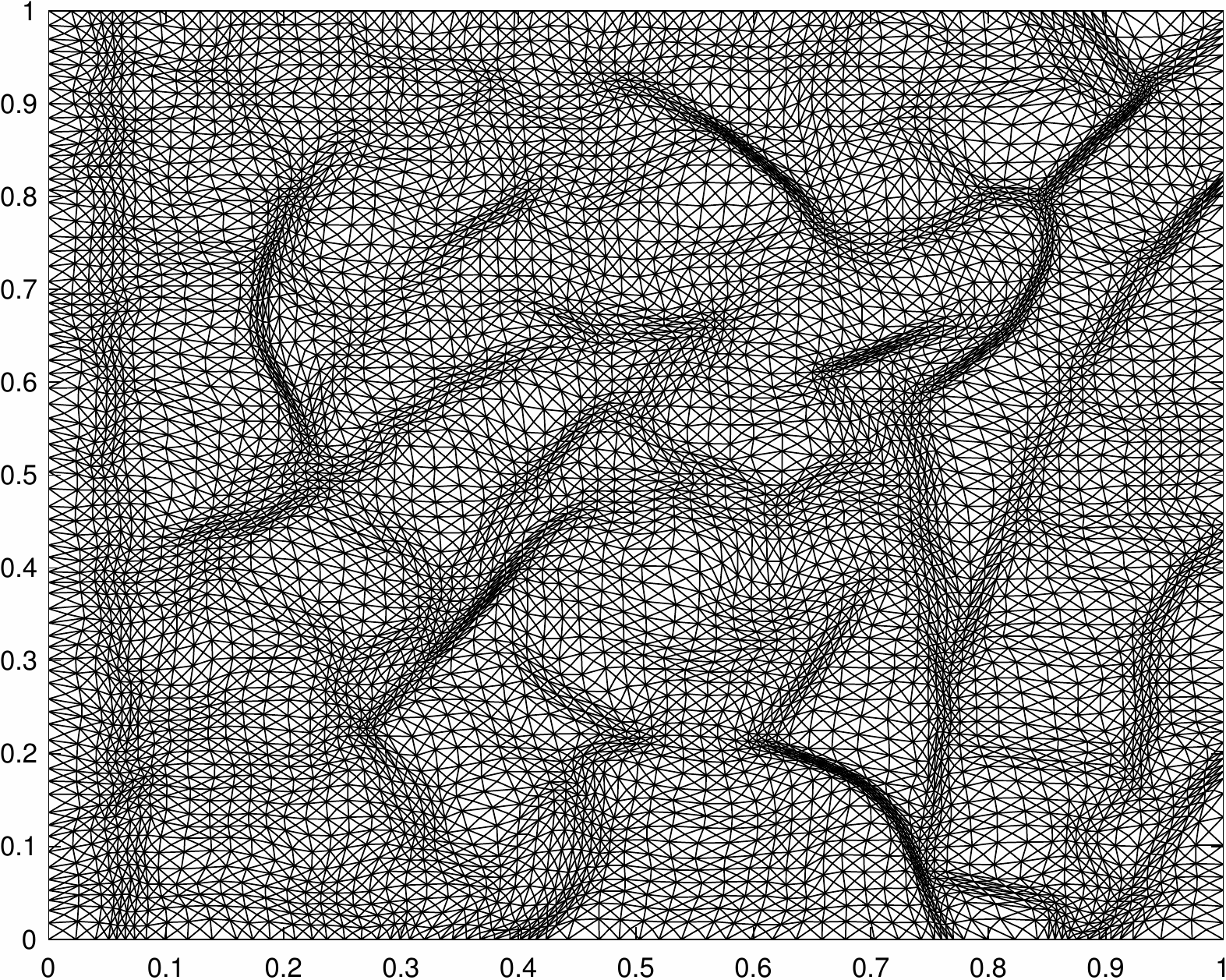}
\caption{ $t = 0.3$, $\epsilon$ is chosen by (\ref{epsilon-1})}
\end{subfigure}
\caption{The meshes corresponding to Fig.~\ref{lenaimage}.}
\label{lenamesh}
\end{figure}

\begin{figure}[htb]
\centering
\begin{subfigure}{0.32\textwidth}
\centering
\includegraphics[scale = 0.24]{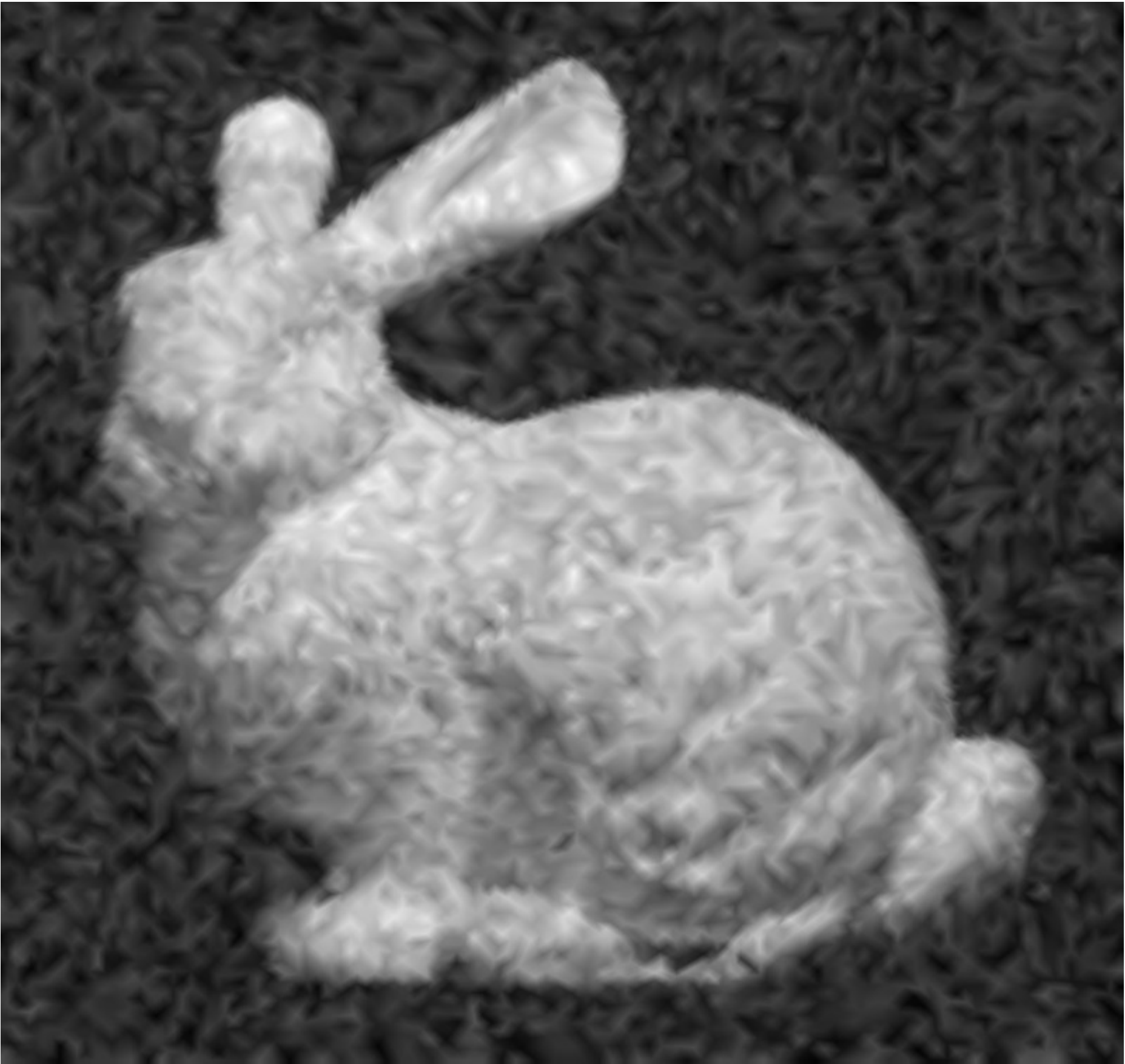}
\caption{$t = 0.002$, $\epsilon = 10^{-7}$}
\end{subfigure}
\begin{subfigure}{0.32\textwidth}
\centering
\includegraphics[scale = 0.24]{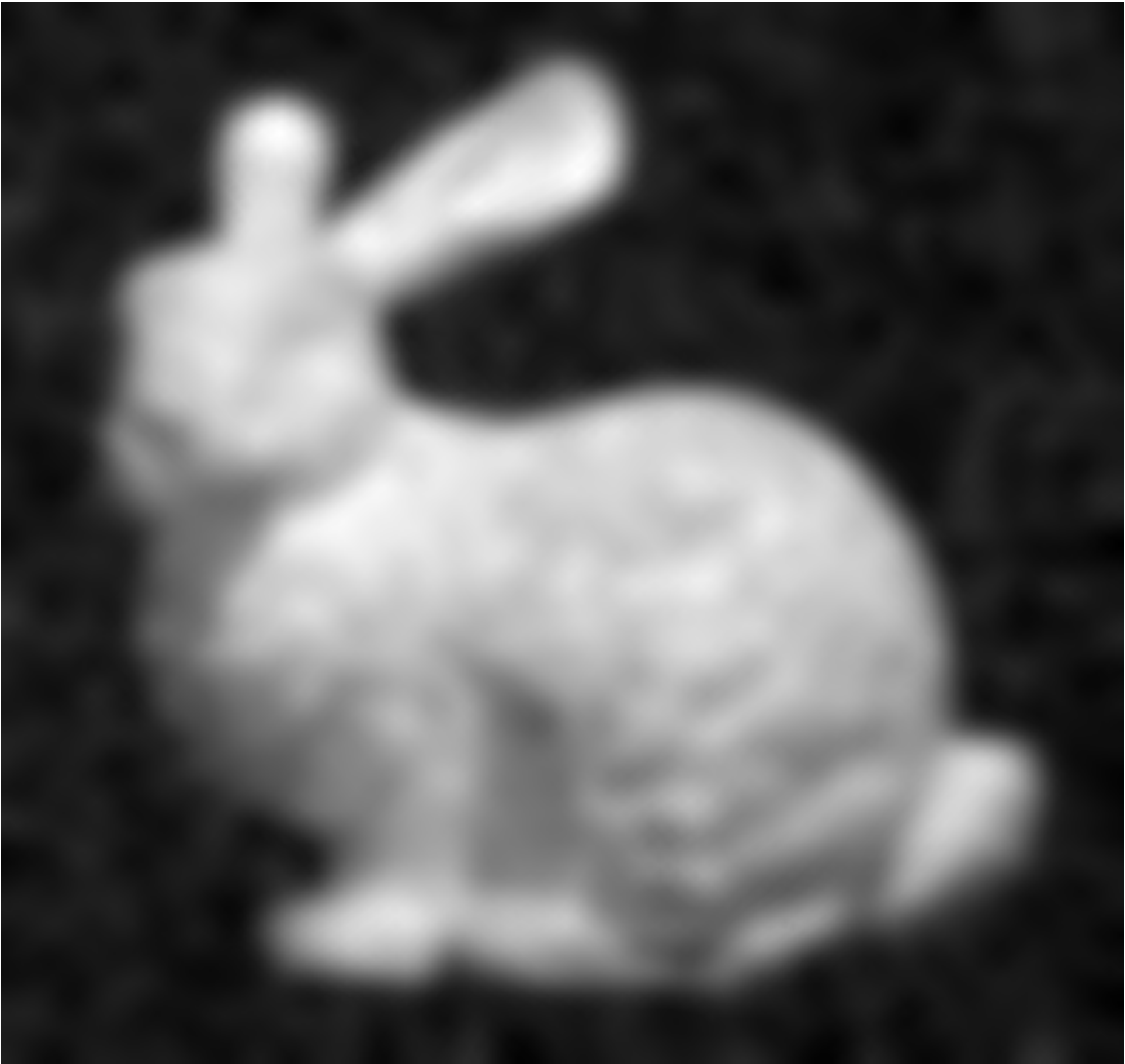}
\caption{$t = 0.008$, $\epsilon = 10^{-7}$}
\end{subfigure}
 \begin{subfigure}{0.32\textwidth}
 \centering
\includegraphics[scale = 0.24]{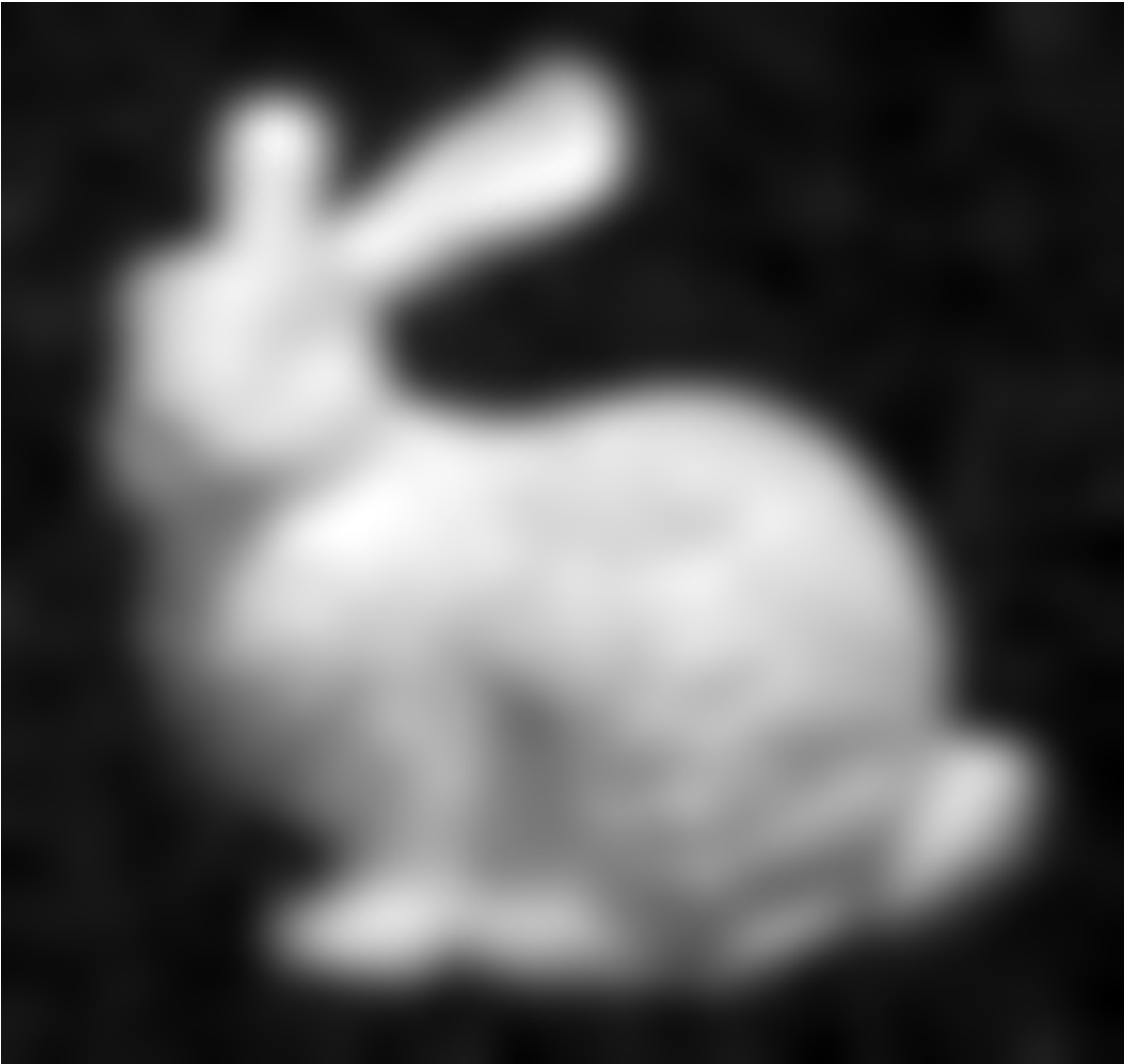}
\caption{$t = 0.2$, $\epsilon = 10^{-7}$}
\end{subfigure}
\begin{subfigure}{0.32\textwidth}
\centering
\includegraphics[scale = 0.24]{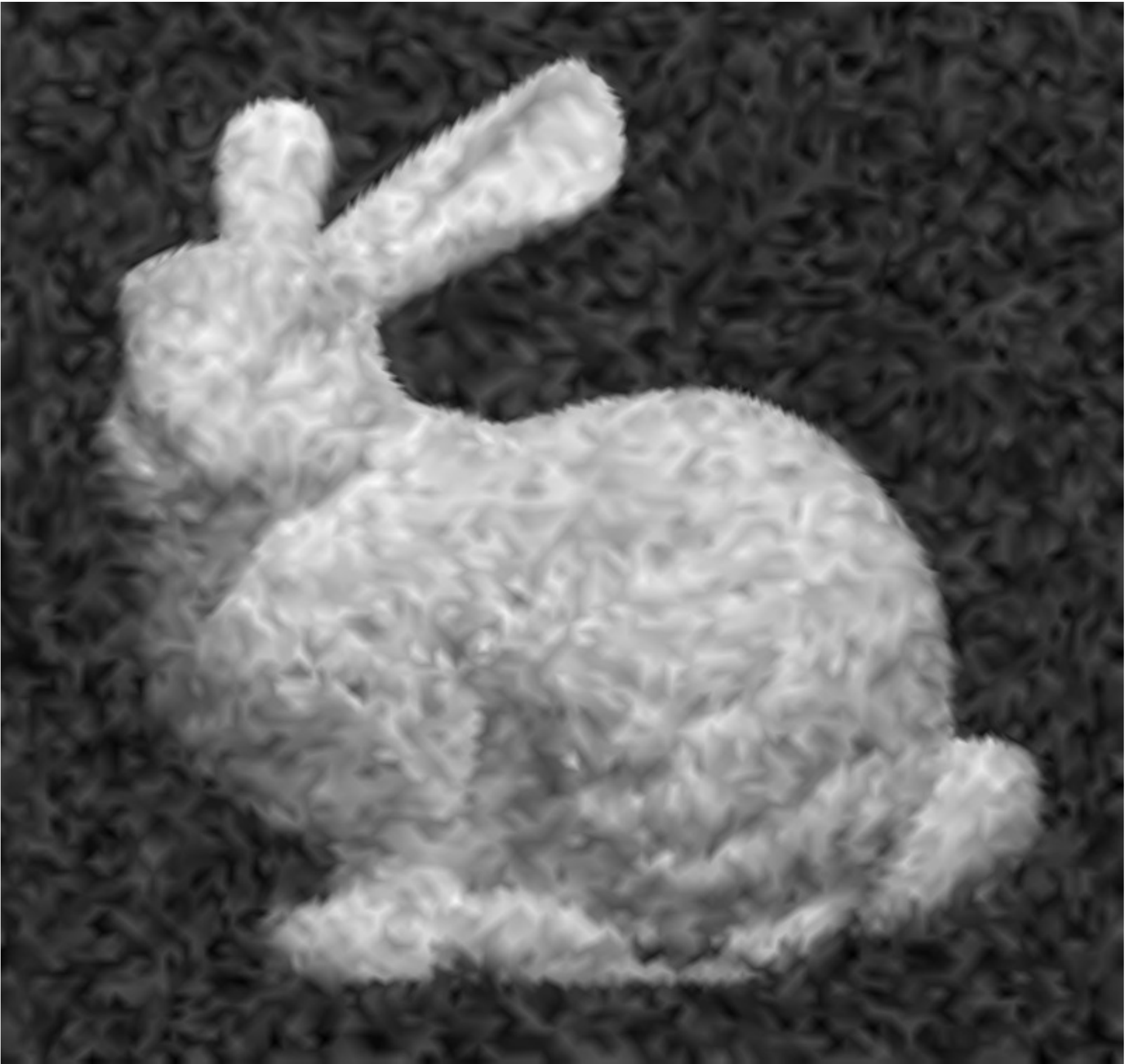}
\caption{$t = 0.002$, $\epsilon$ is chosen by (\ref{epsilon-1})}
\end{subfigure}
\begin{subfigure}{0.32\textwidth}
\centering
\includegraphics[scale = 0.24]{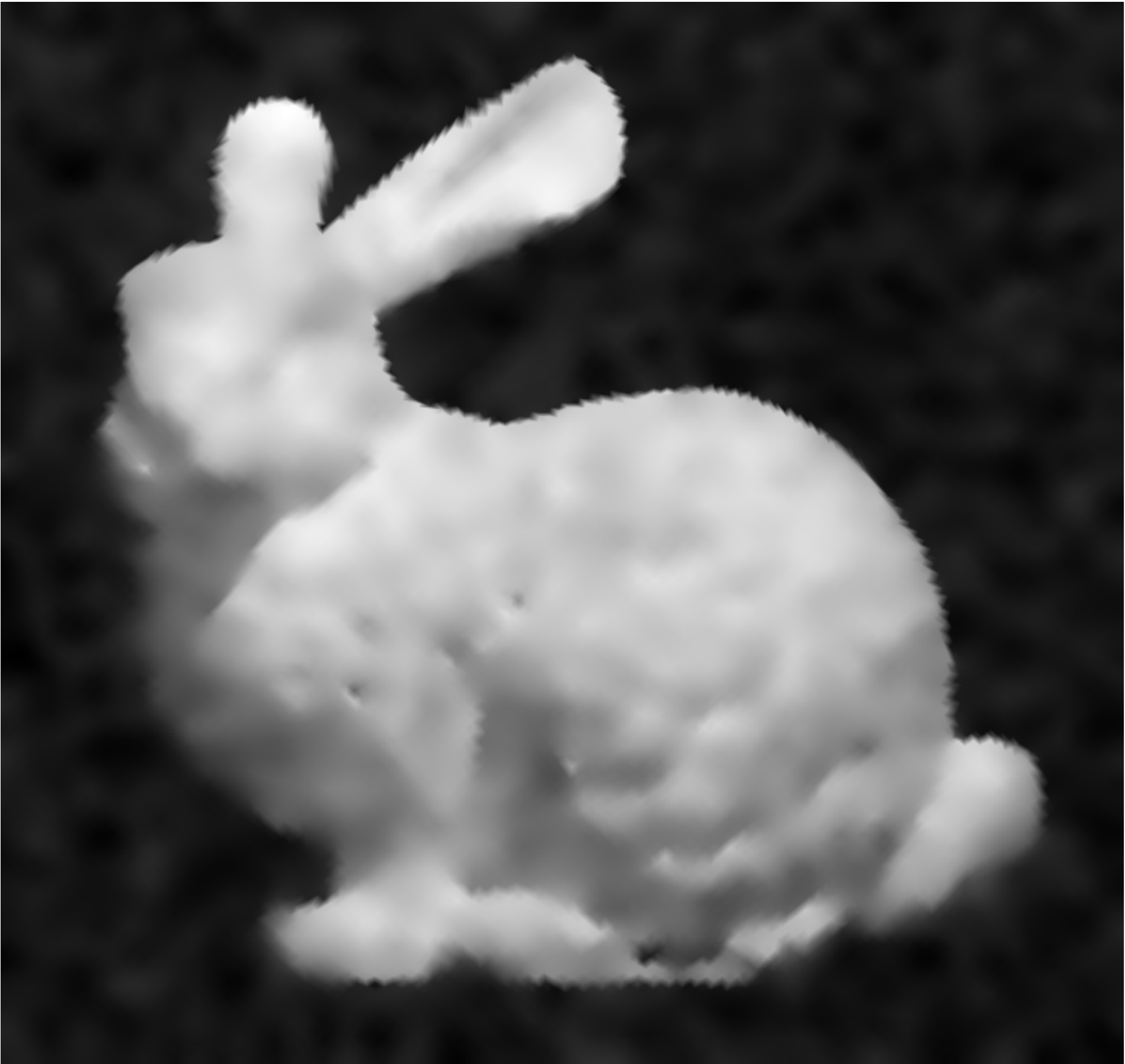}
\caption{$t = 0.008$, $\epsilon$ is chosen by (\ref{epsilon-1})}
\end{subfigure}
 \begin{subfigure}{0.32\textwidth}
 \centering
\includegraphics[scale = 0.24]{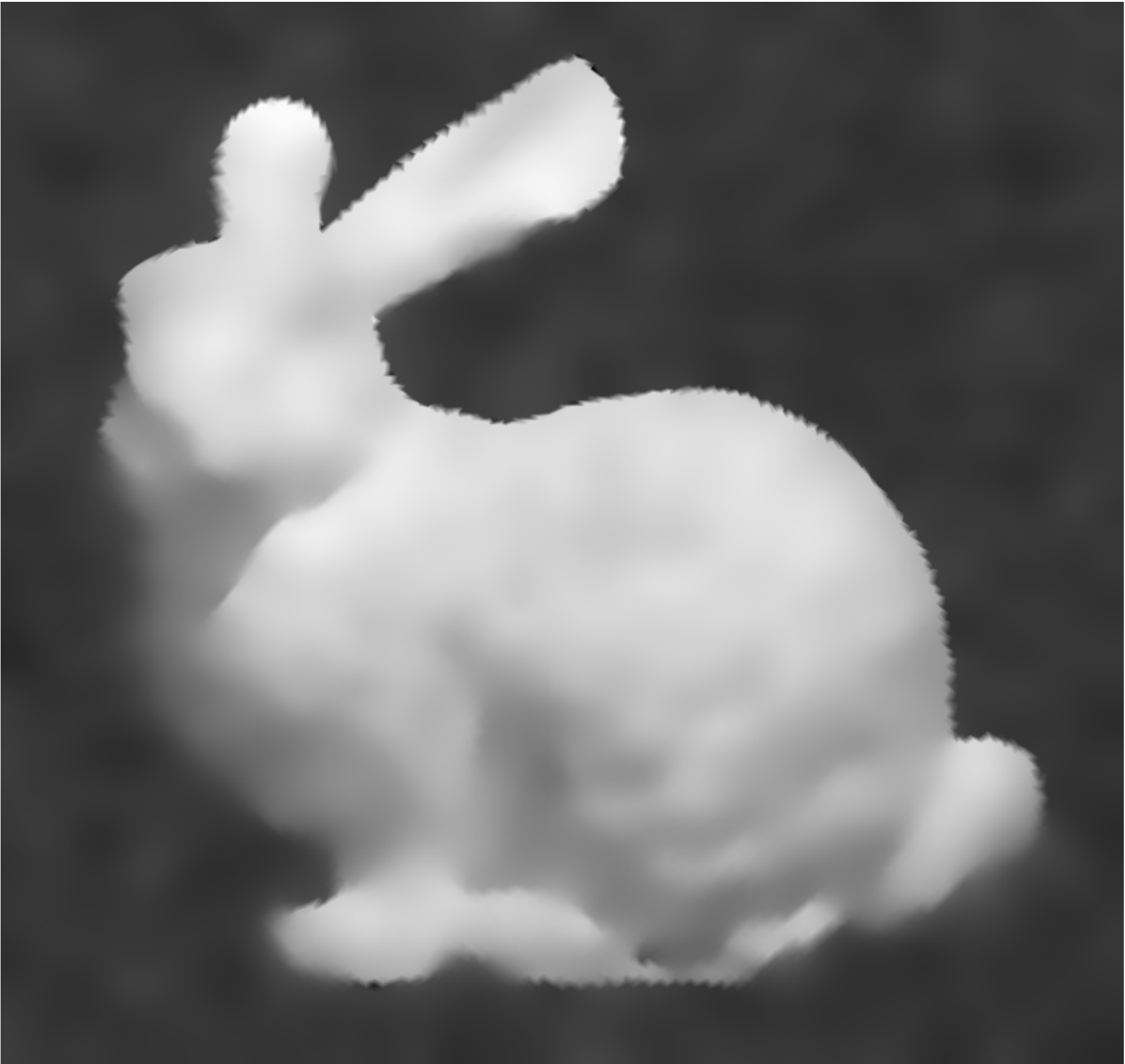}
\caption{$t = 0.2$, $\epsilon$ is chosen by (\ref{epsilon-1})}
\end{subfigure}
\caption{Evolution of the image.}
\label{bunnyimage}
\end{figure}

\begin{figure}[htb]
\centering
\begin{subfigure}{0.32\textwidth}
\centering
\includegraphics[scale = 0.3]{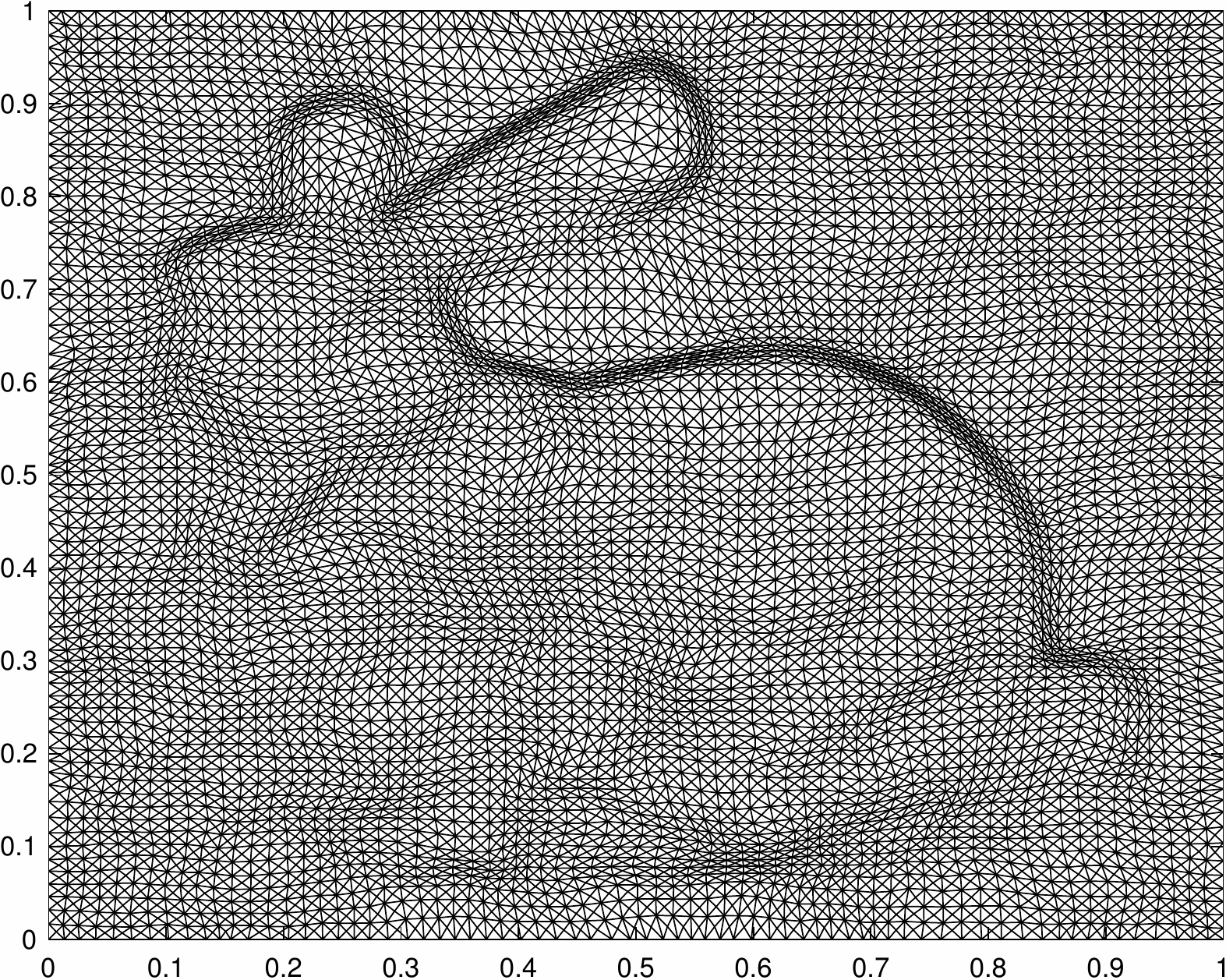}
\caption{$t = 0.002$, $\epsilon = 10^{-7}$}
\end{subfigure}
\begin{subfigure}{0.32\textwidth}
\centering
\includegraphics[scale = 0.3]{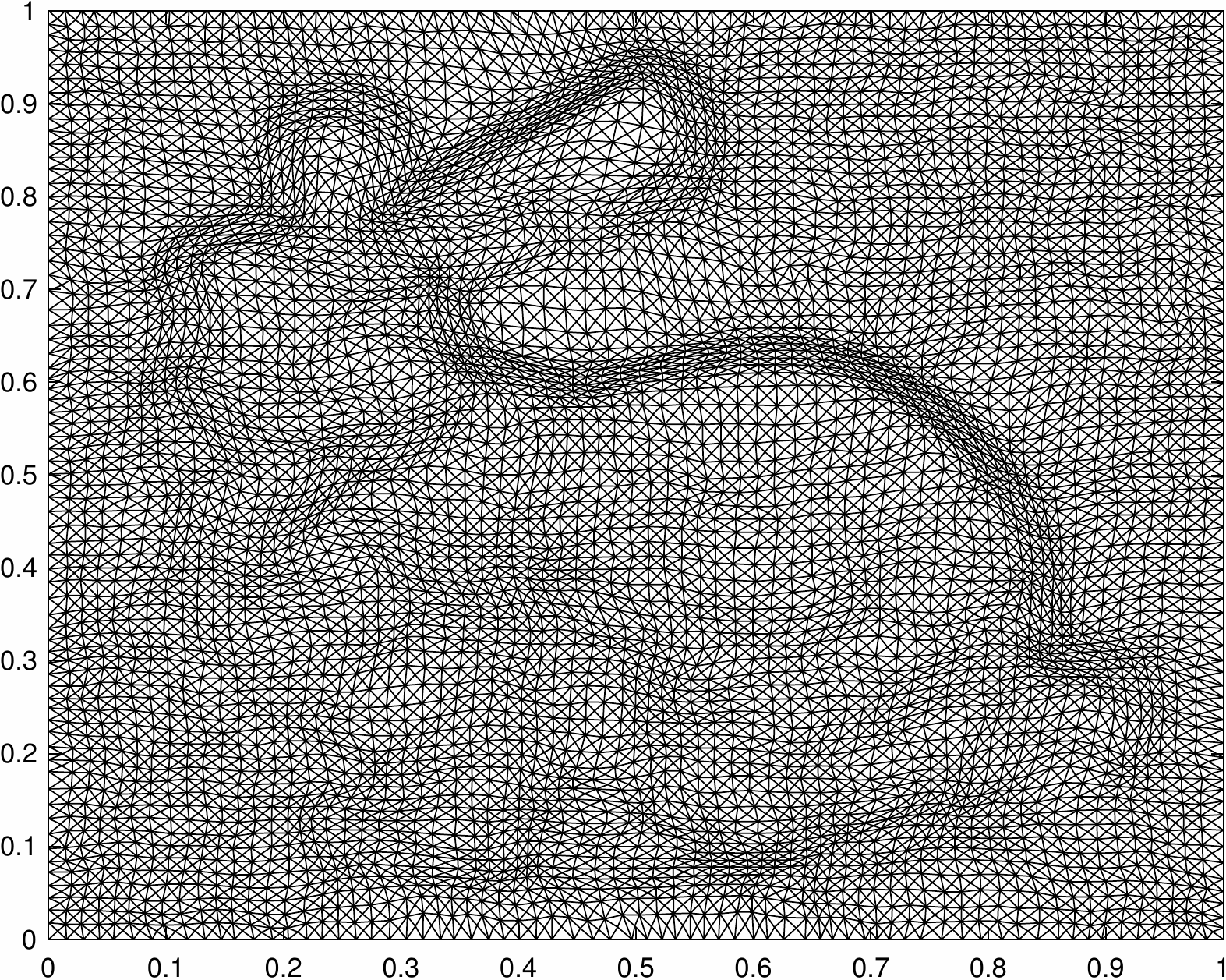}
\caption{$t = 0.008$, $\epsilon = 10^{-7}$}
\end{subfigure}
 \begin{subfigure}{0.32\textwidth}
 \centering
\includegraphics[scale = 0.3]{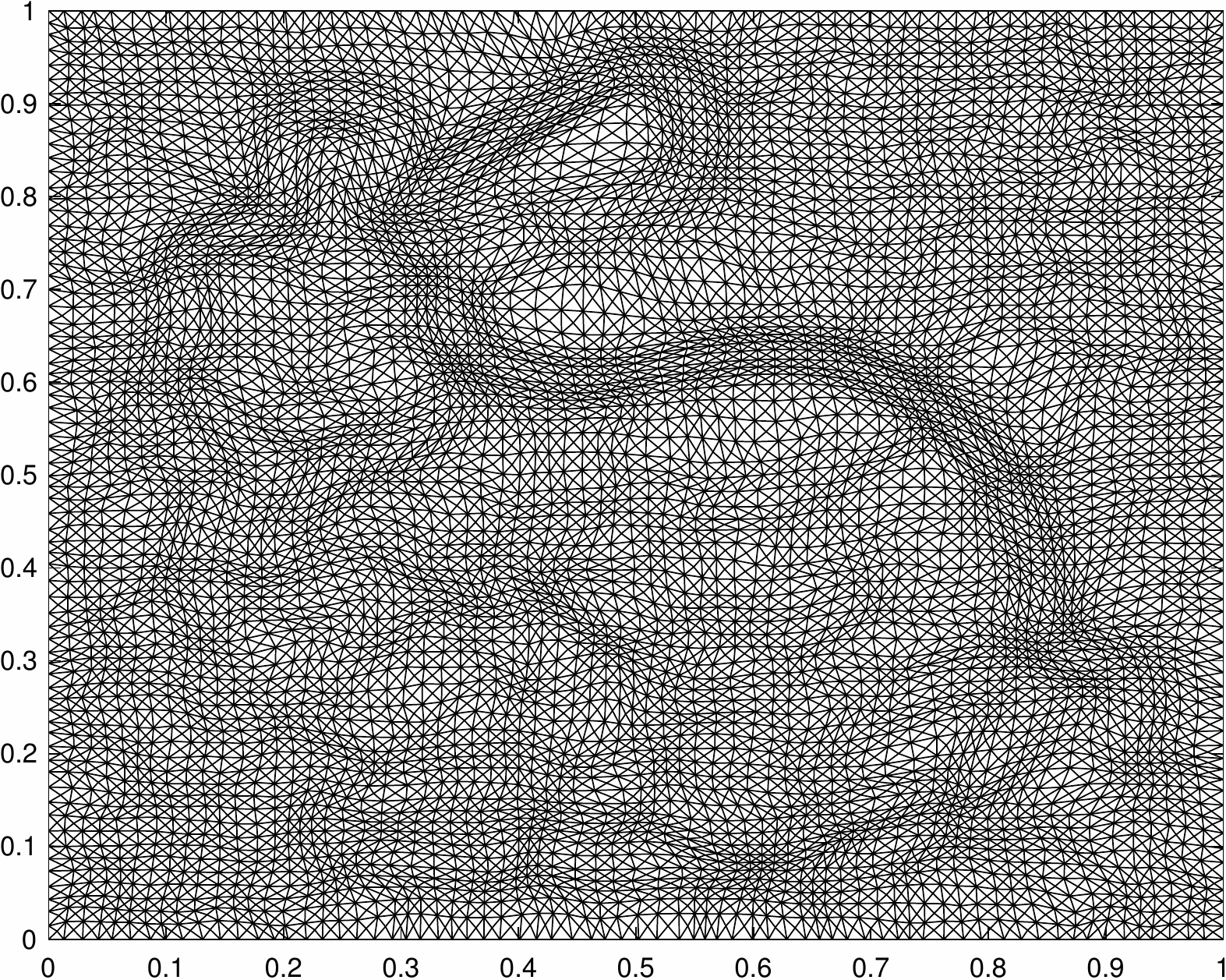}
\caption{$t = 0.2$, $\epsilon = 10^{-7}$}
\end{subfigure}
\begin{subfigure}{0.32\textwidth}
\centering
\includegraphics[scale = 0.3]{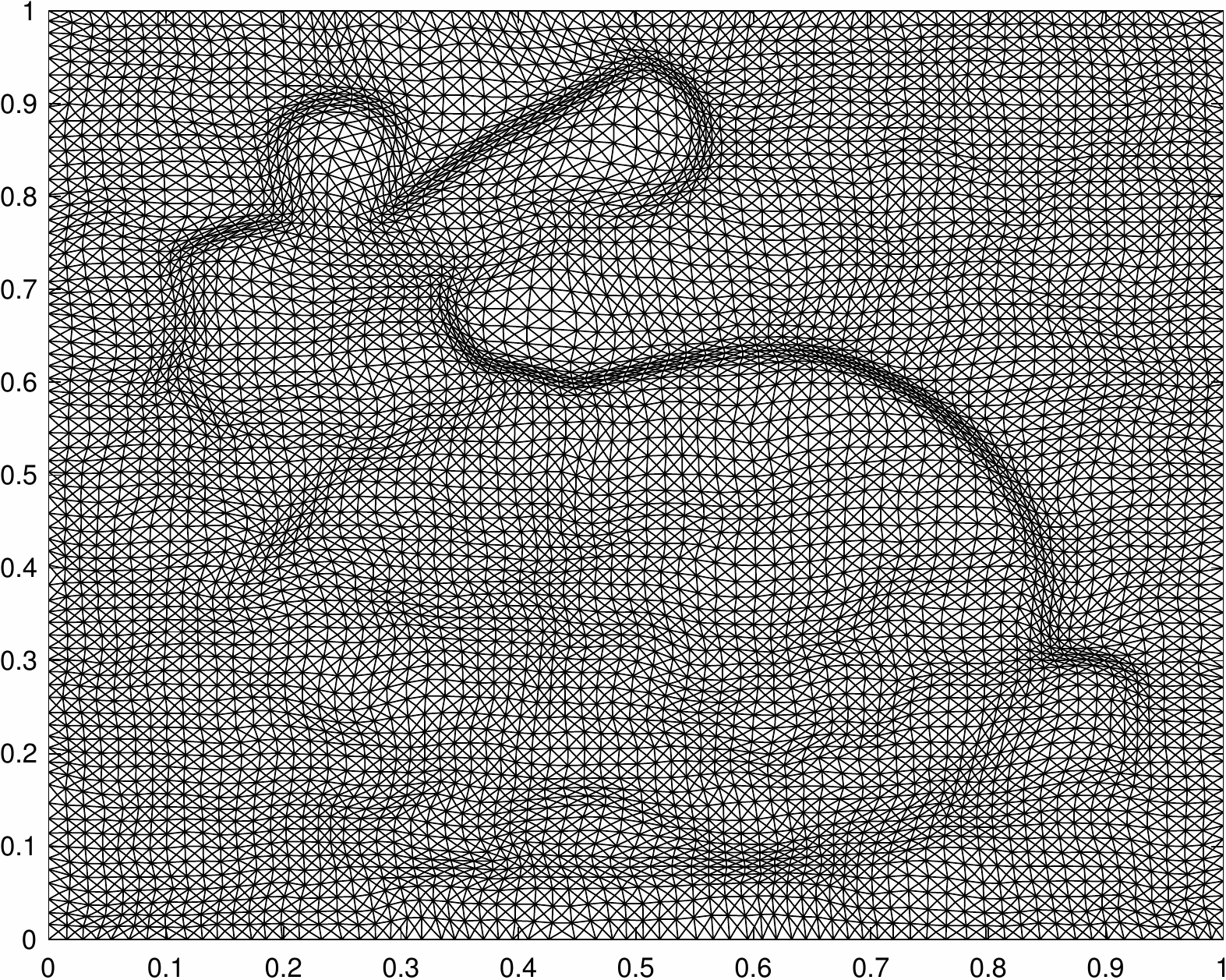}
\caption{$t = 0.002$, $\epsilon$ is chosen by (\ref{epsilon-1})}
\end{subfigure}
\begin{subfigure}{0.32\textwidth}
\centering
\includegraphics[scale = 0.3]{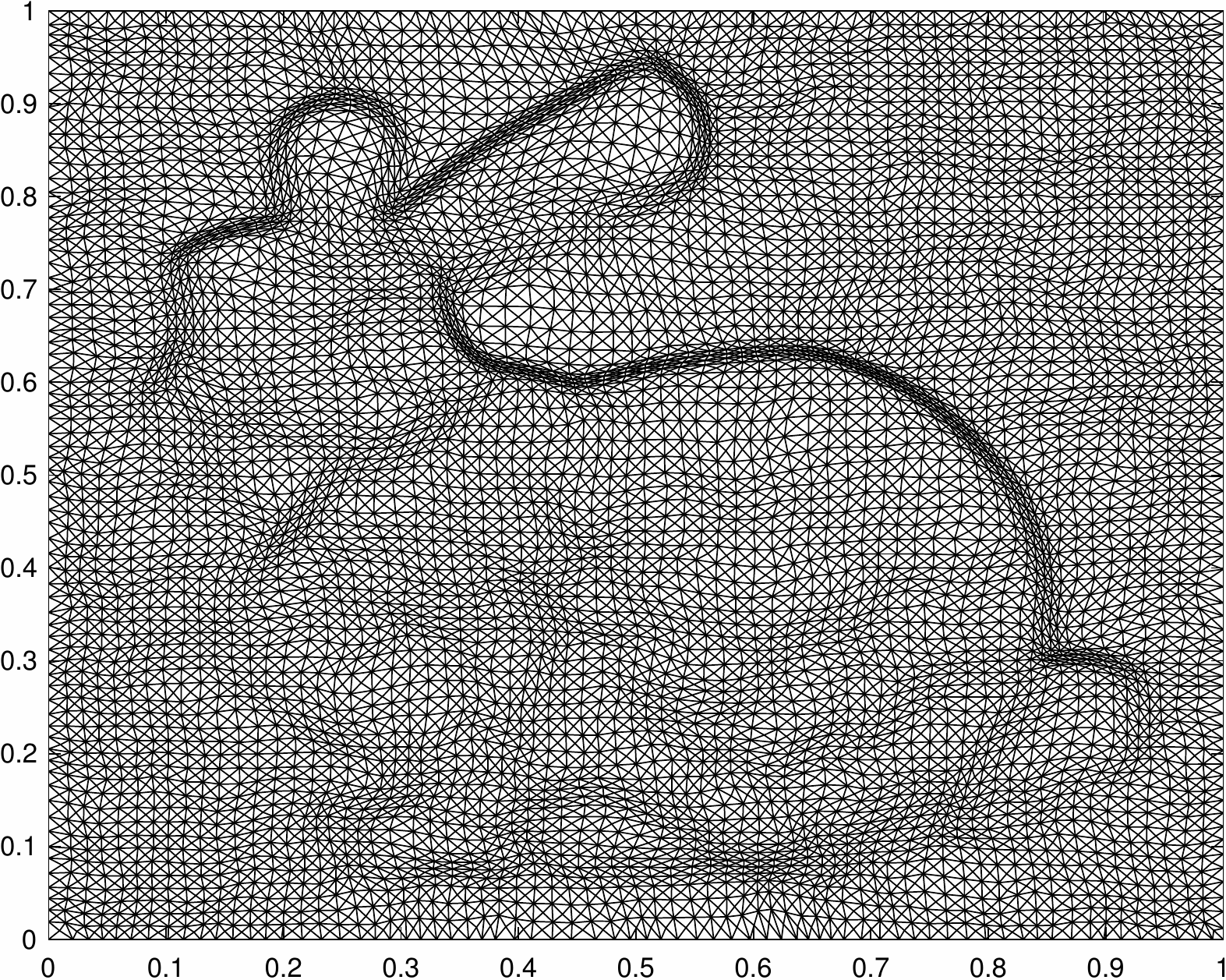}
\caption{$t = 0.008$, $\epsilon$ is chosen by (\ref{epsilon-1})}
\end{subfigure}
 \begin{subfigure}{0.32\textwidth}
 \centering
\includegraphics[scale = 0.3]{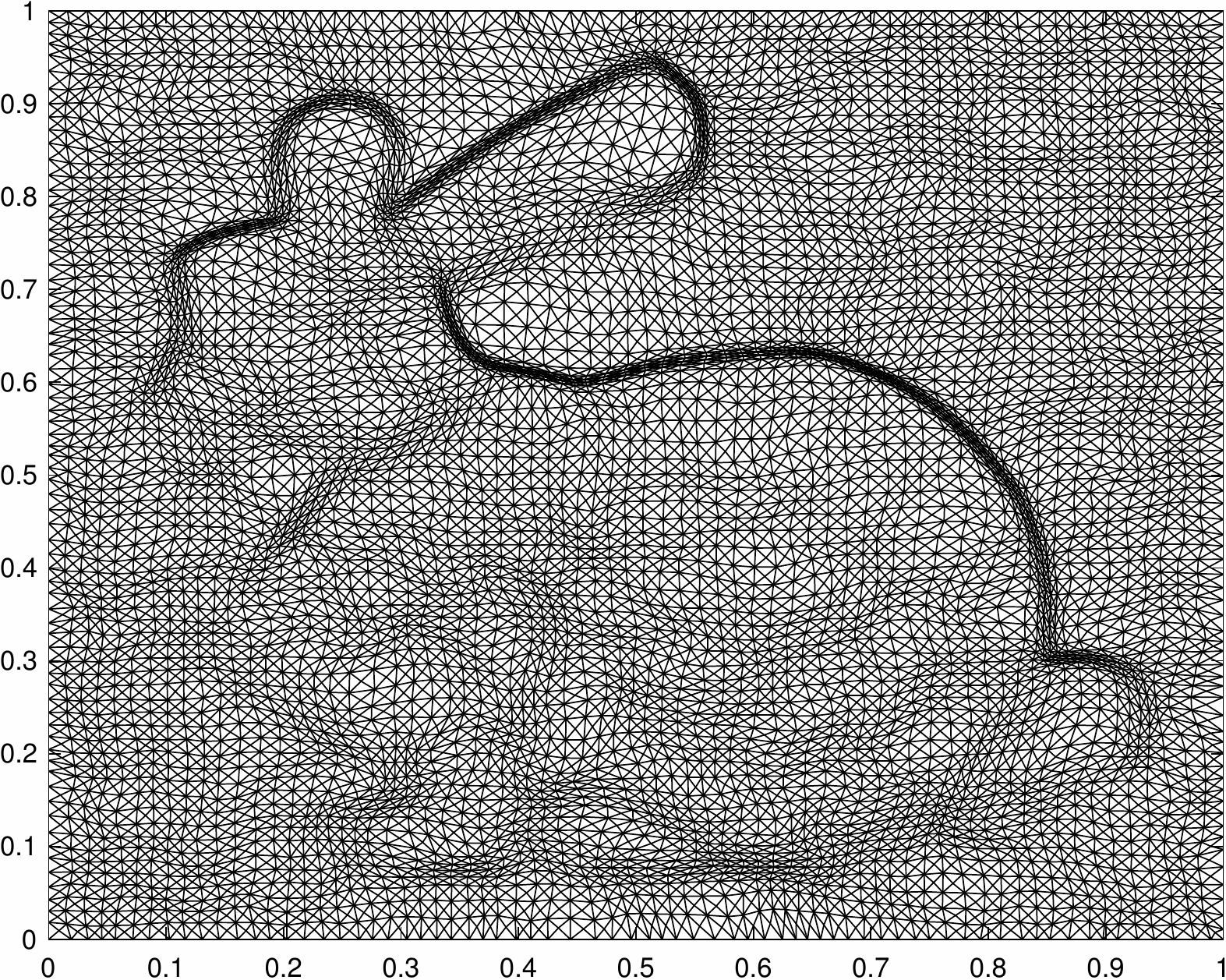}
\caption{$t = 0.2$, $\epsilon$ is chosen by (\ref{epsilon-1})}
\end{subfigure}
\caption{The meshes corresponding to Fig.~\ref{bunnyimage}.}
\label{bunnymesh}
\end{figure}

\begin{figure}[htb]
\centering
\begin{subfigure}{0.32\textwidth}
\centering
\includegraphics[scale = 0.22]{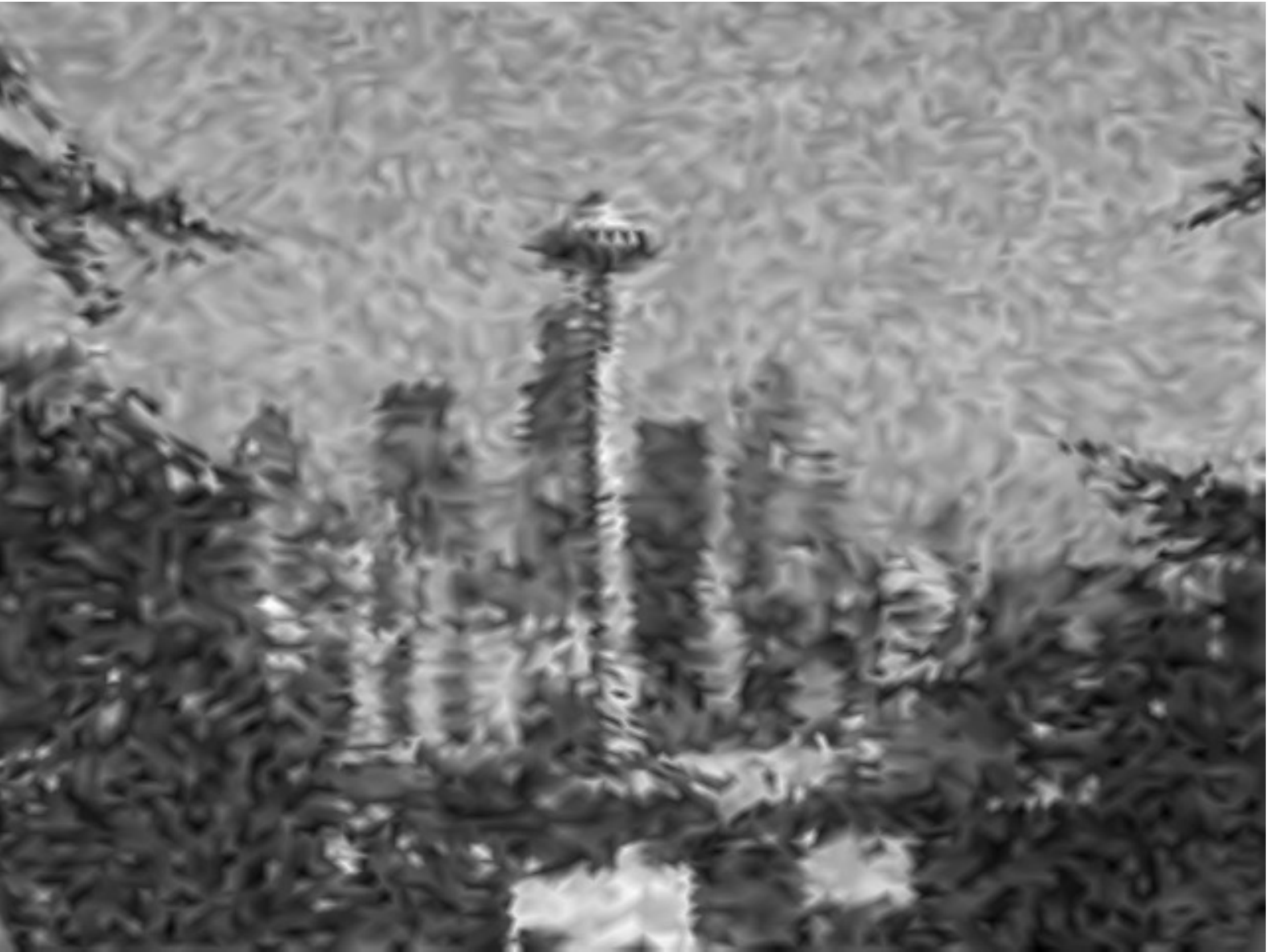}
\caption{$t = 0.00013$, $\epsilon = 10^{-7}$}
\end{subfigure}
\begin{subfigure}{0.32\textwidth}
\centering
\includegraphics[scale = 0.22]{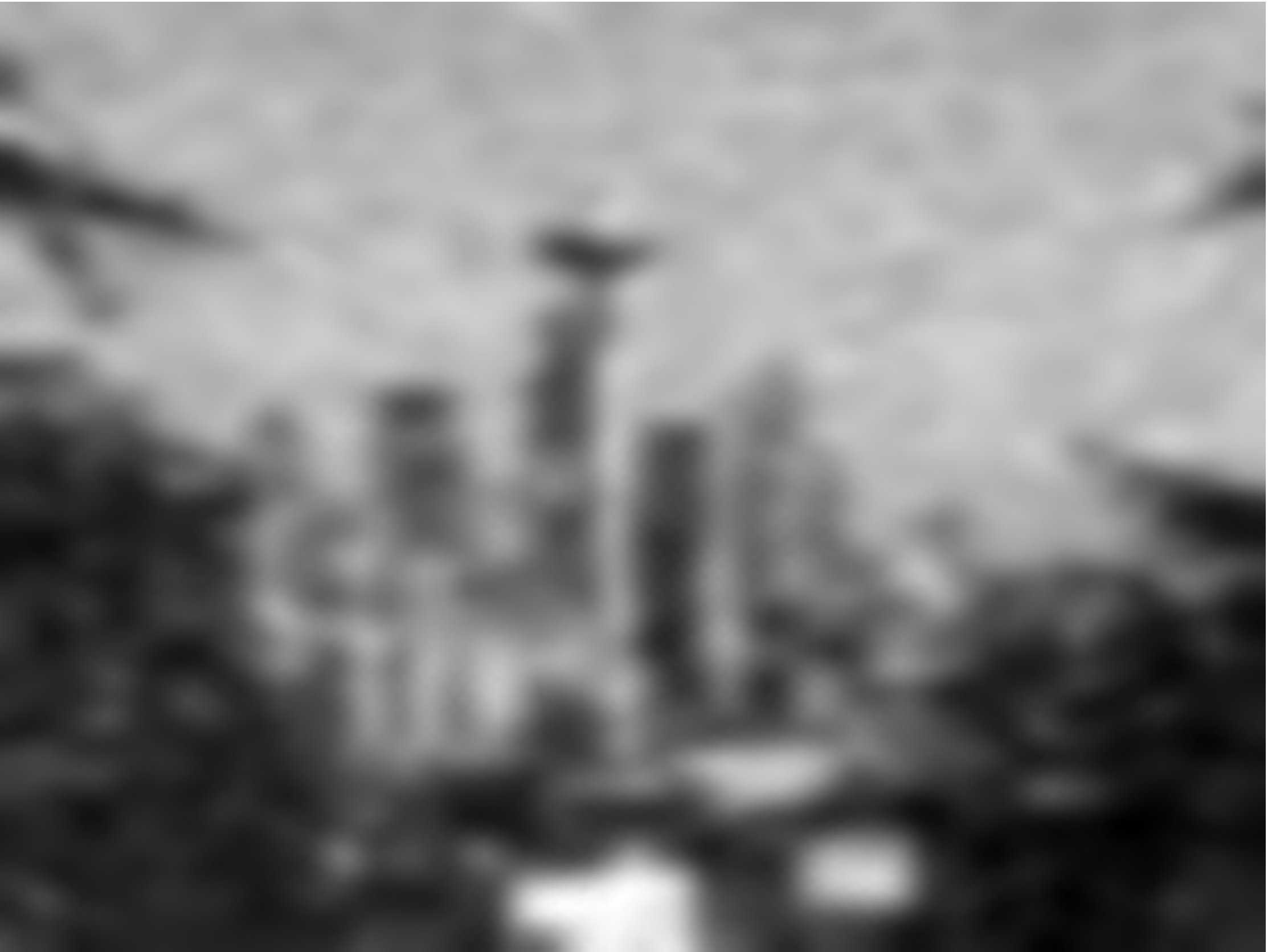}
\caption{$t = 0.08$, $\epsilon = 10^{-7}$}
\end{subfigure}
 \begin{subfigure}{0.32\textwidth}
 \centering
\includegraphics[scale = 0.22]{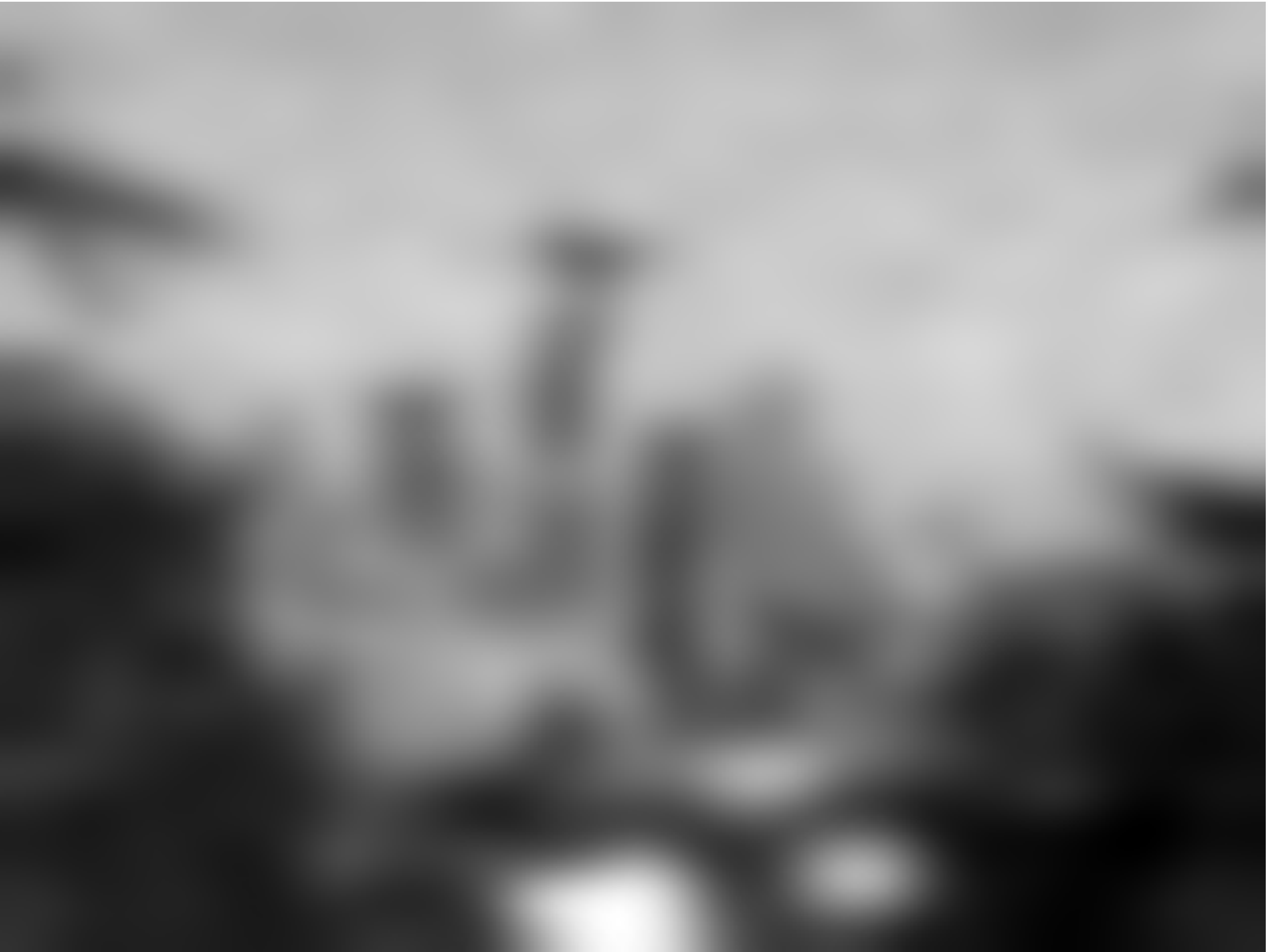}
\caption{$t = 0.3$, $\epsilon = 10^{-7}$}
\end{subfigure}
\begin{subfigure}{0.32\textwidth}
\centering
\includegraphics[scale = 0.22]{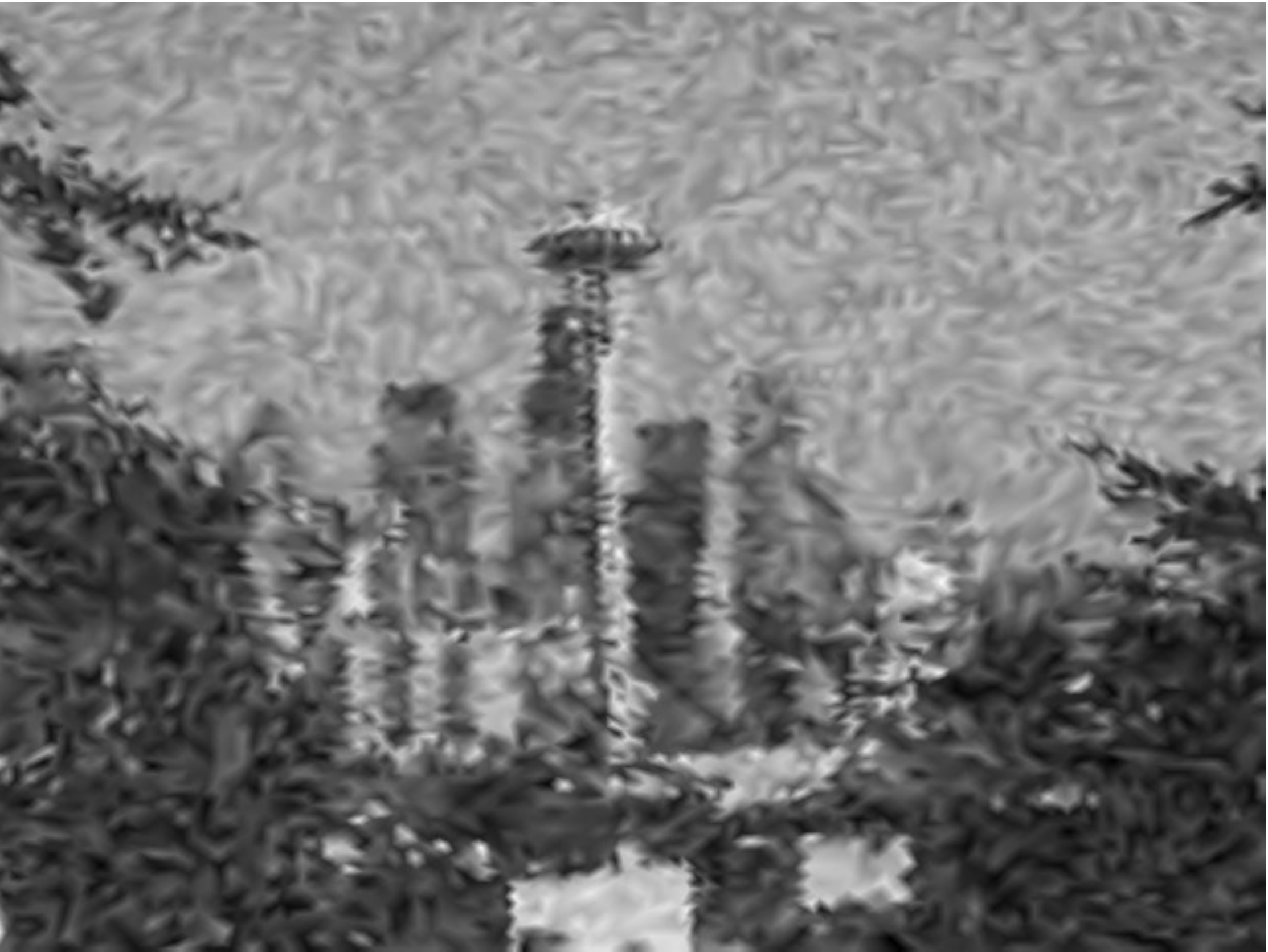}
\caption{$t = 0.00013$, $\epsilon$ is chosen by (\ref{epsilon-1})}
\end{subfigure}
\begin{subfigure}{0.32\textwidth}
\centering
\includegraphics[scale = 0.22]{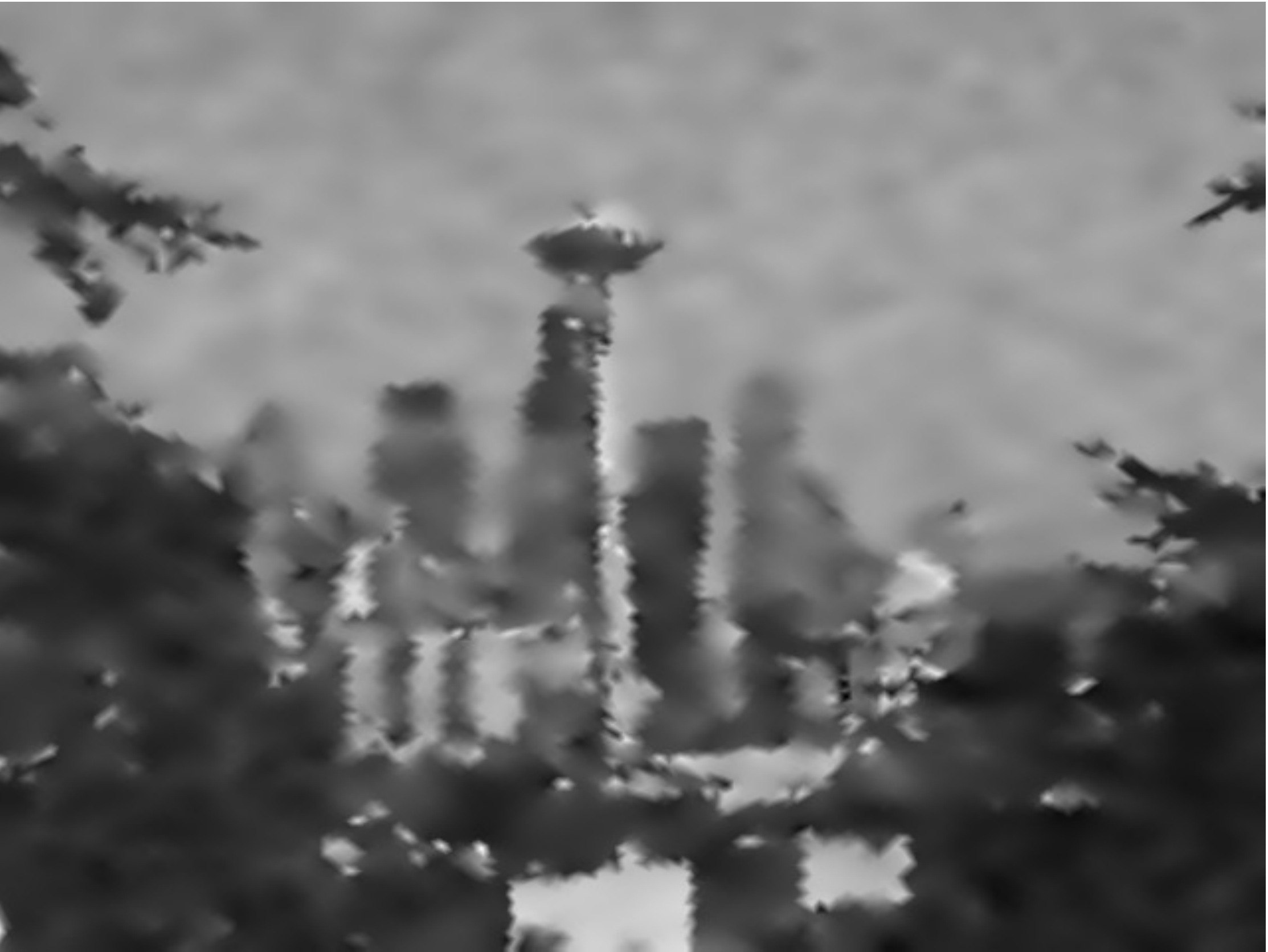}
\caption{$t = 0.08$, $\epsilon$ is chosen by (\ref{epsilon-1})}
\end{subfigure}
 \begin{subfigure}{0.32\textwidth}
 \centering
\includegraphics[scale = 0.22]{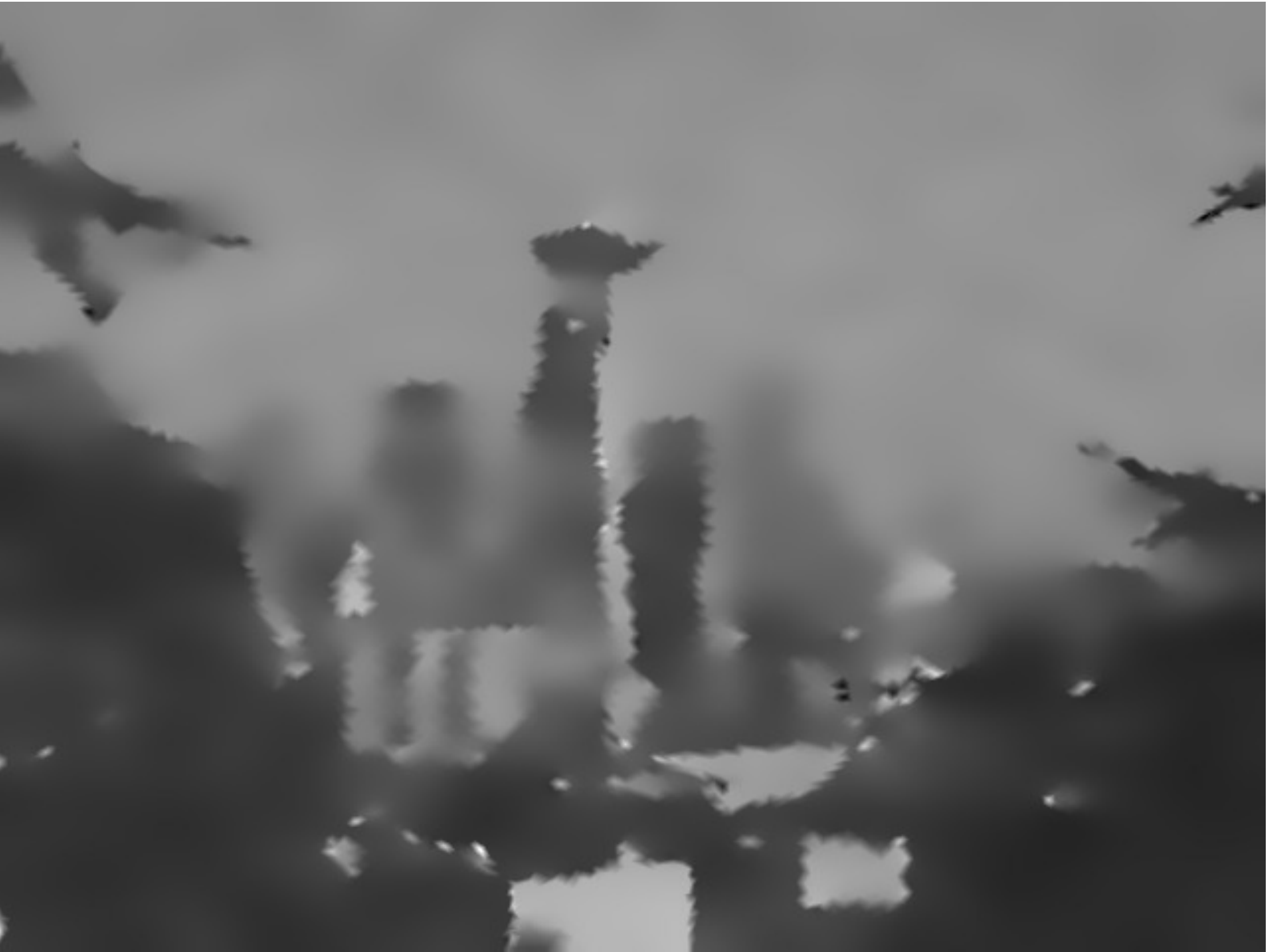}
\caption{$t = 0.3$, $\epsilon$ is chosen by (\ref{epsilon-1})}
\end{subfigure}
\caption{Evolution of the image.}
\label{usimage}
\end{figure}

\begin{figure}[htb]
\centering
\begin{subfigure}{0.32\textwidth}
\centering
\includegraphics[scale = 0.32]{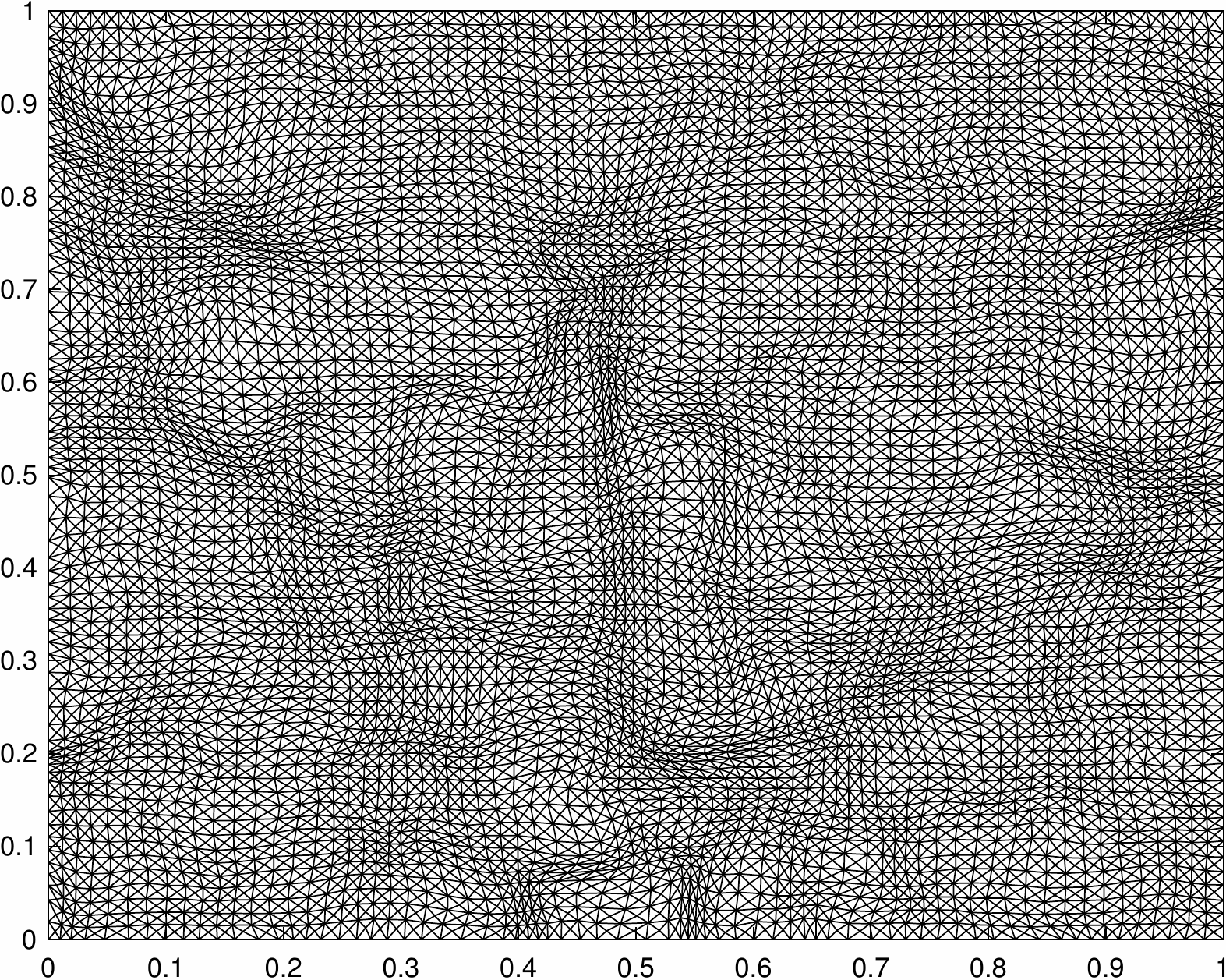}
\caption{$t = 0.00013$, $\epsilon = 10^{-7}$}
\end{subfigure}
\begin{subfigure}{0.32\textwidth}
\centering
\includegraphics[scale = 0.32]{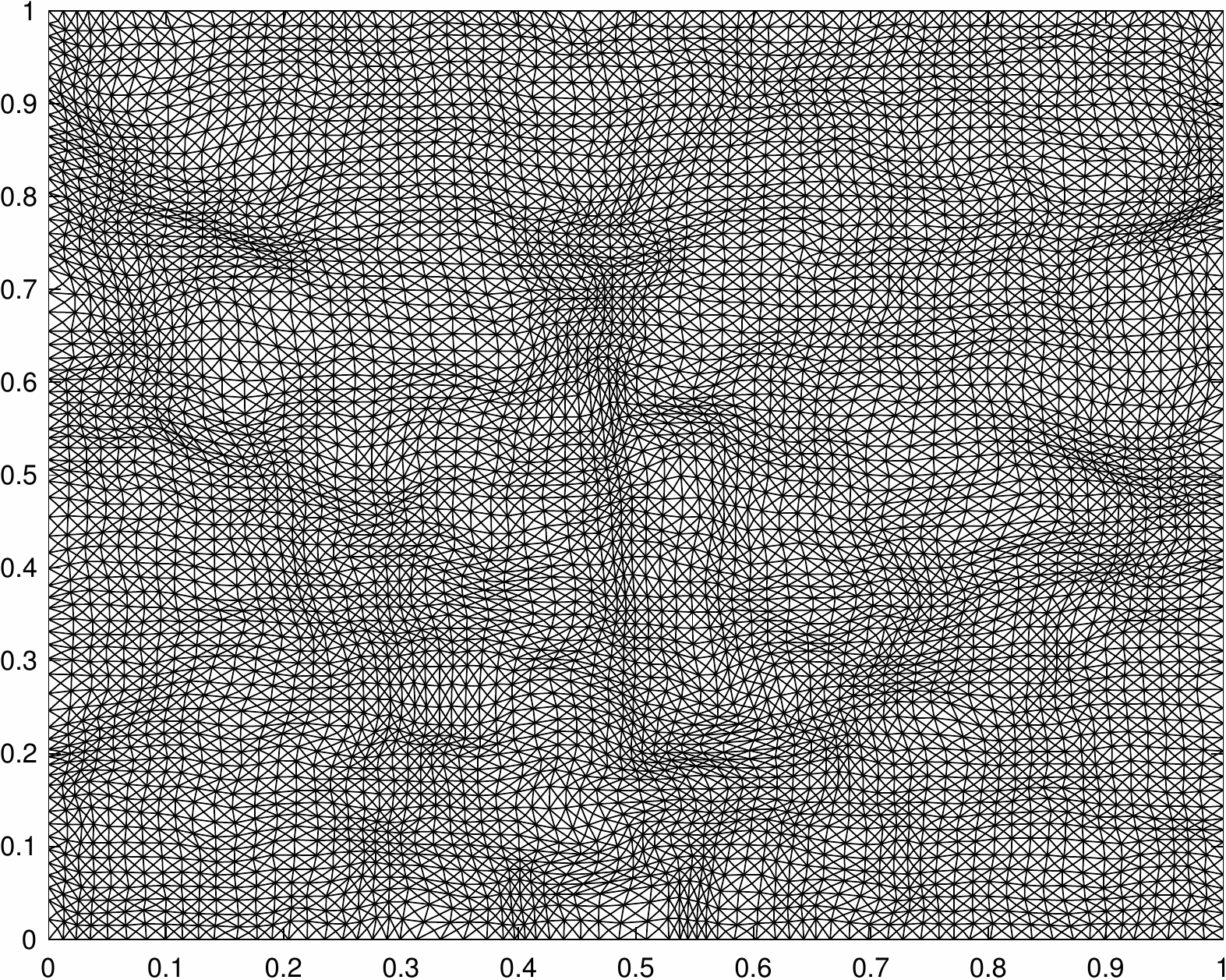}
\caption{$t = 0.08$, $\epsilon = 10^{-7}$}
\end{subfigure}
 \begin{subfigure}{0.32\textwidth}
 \centering
\includegraphics[scale = 0.32]{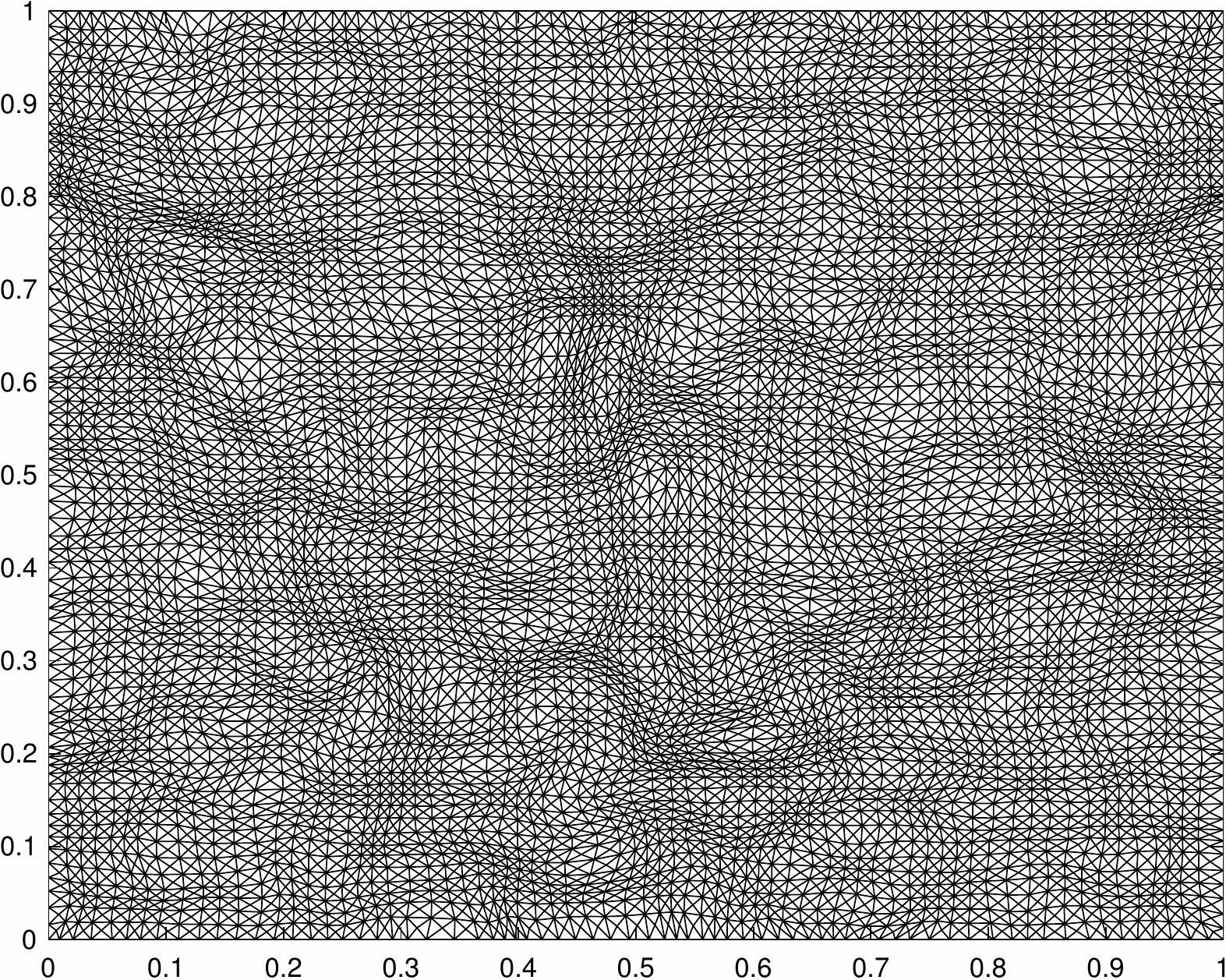}
\caption{$t = 0.3$, $\epsilon = 10^{-7}$}
\end{subfigure}
\begin{subfigure}{0.32\textwidth}
\centering
\includegraphics[scale = 0.32]{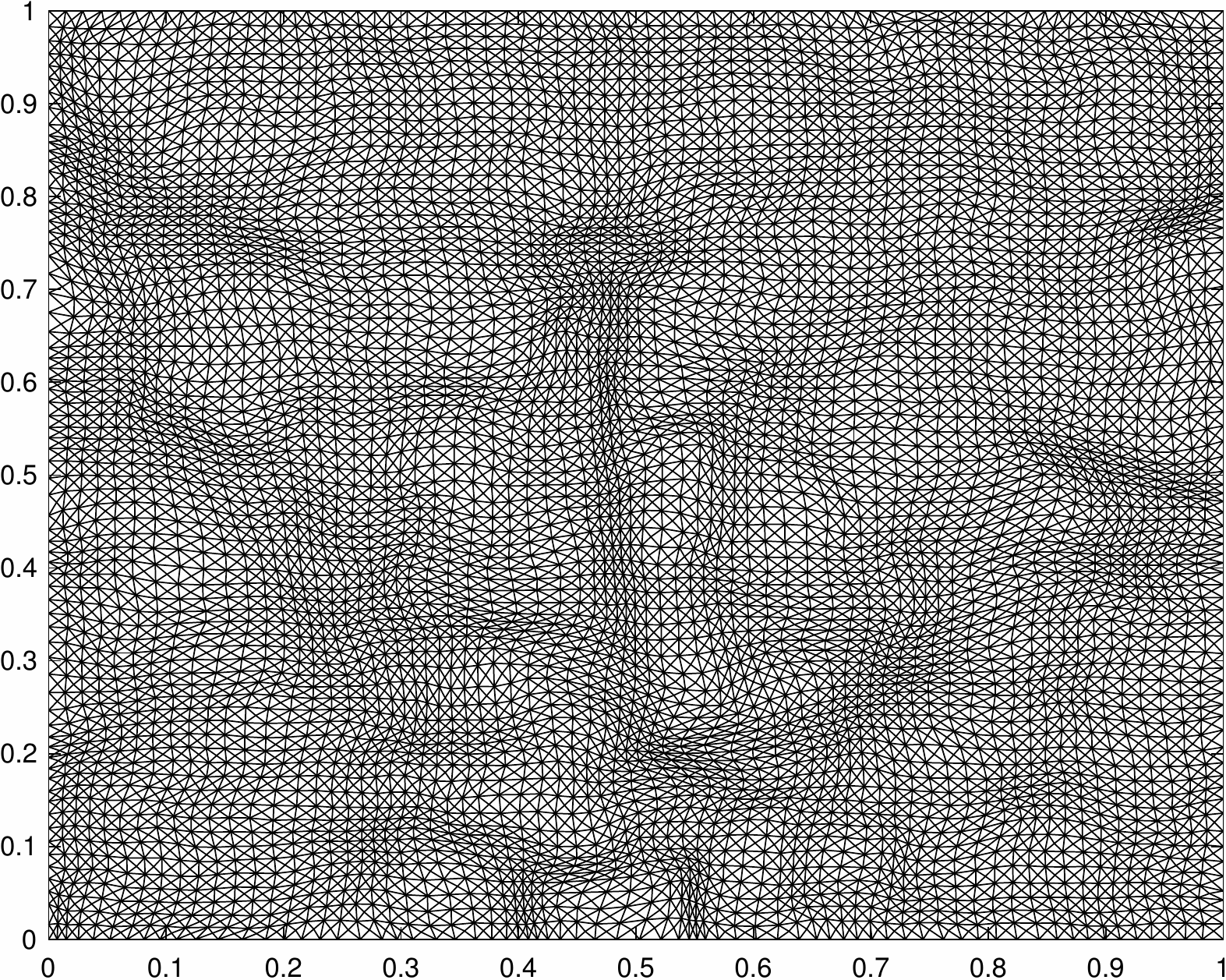}
\caption{$t = 0.00013$, $\epsilon$ is chosen by (\ref{epsilon-1})}
\end{subfigure}
\begin{subfigure}{0.32\textwidth}
\centering
\includegraphics[scale = 0.32]{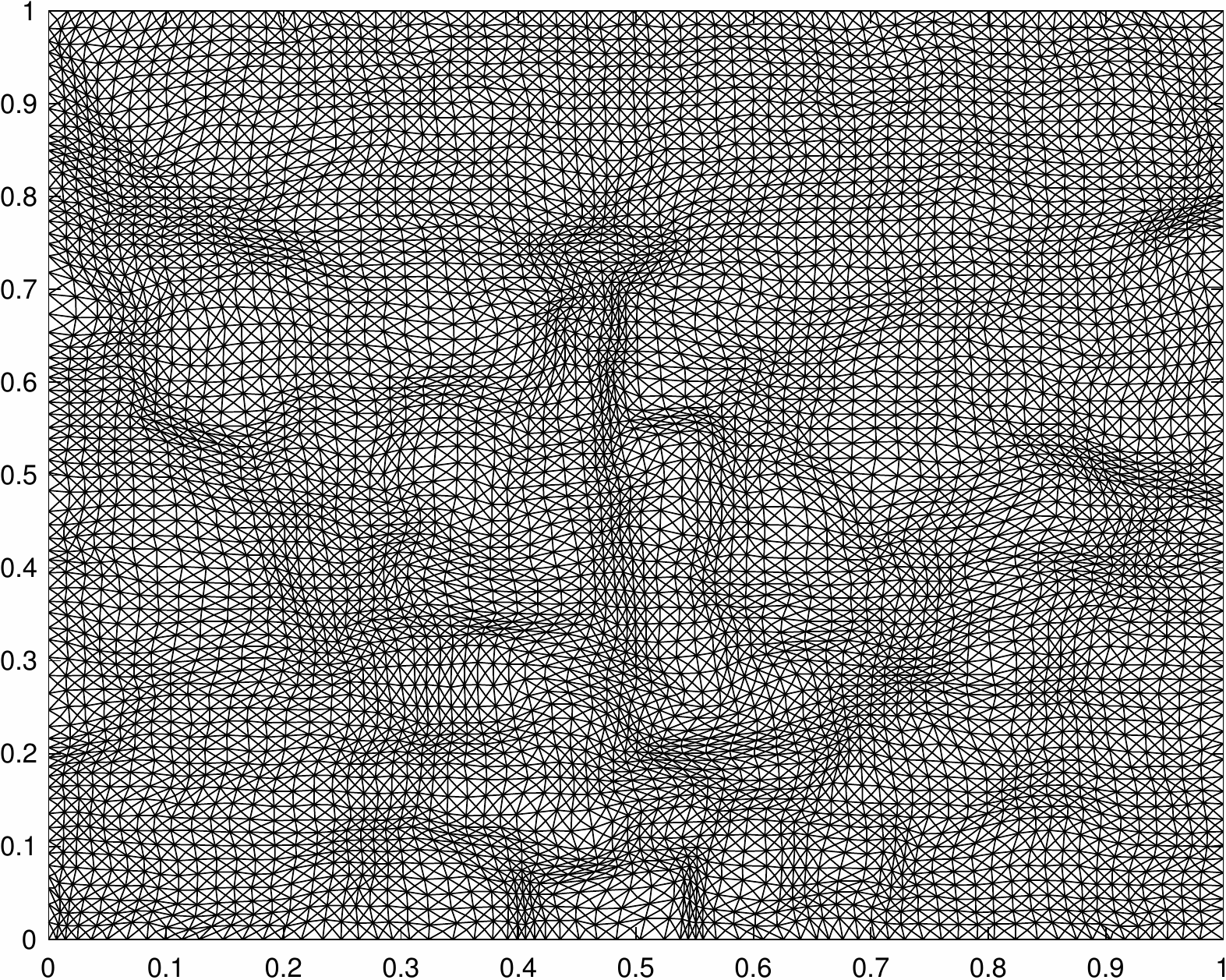}
\caption{$t = 0.08$, $\epsilon$ is chosen by (\ref{epsilon-1})}
\end{subfigure}
 \begin{subfigure}{0.32\textwidth}
 \centering
\includegraphics[scale = 0.32]{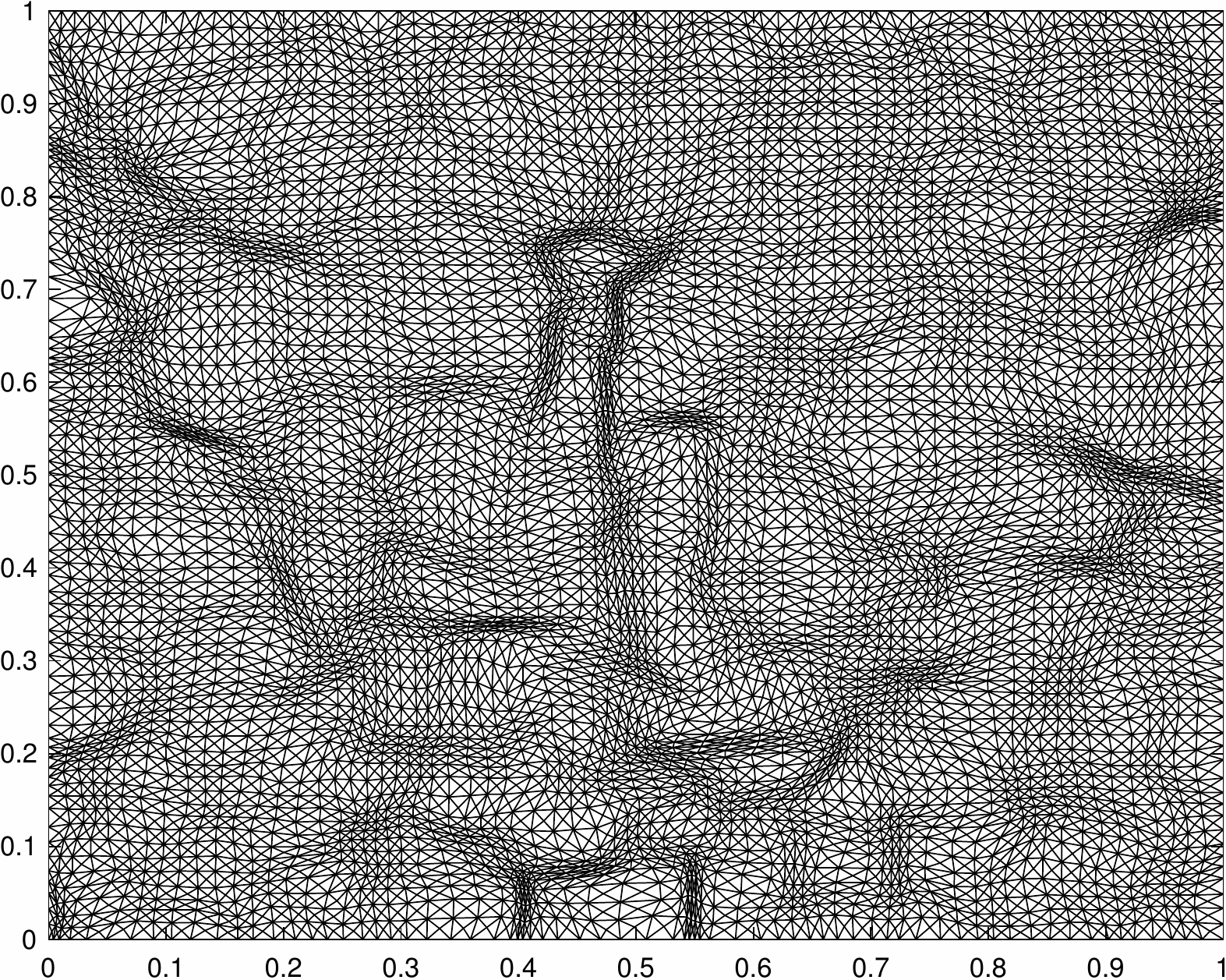}
\caption{$t = 0.3$, $\epsilon$ is chosen by (\ref{epsilon-1})}
\end{subfigure}
\caption{The meshes corresponding to Fig.~\ref{usimage}.}
\label{usmesh}
\end{figure}

To further demonstrate the effects of the selection strategy (\ref{epsilon-1}) and the scaling (\ref{L-1})
we present results obtained for four real images. The results are shown in Figs.~\ref{cameraman},
\ref{lenaimage},  \ref{bunnyimage}, and \ref{usimage} and the corresponding meshes are shown in
Figs.~\ref{cameramanmesh}, \ref{lenamesh},  \ref{bunnymesh}, and \ref{usmesh}, respectively. In these four experiments, 
$N = 2 \times 70 \times 70$, $\alpha = 10^{-3}$, $\gamma = 10^{-5}$, $\beta = 10^{-2}$, and $k_{\epsilon} = 10^{-10}$
are used. A random field in the range $(-0.25, 0.25)$ is added to $g$ as well as $u^0$.
One can observe that the selection strategy (\ref{epsilon-1}) for the regularization parameter
significantly improves segmentation for all cases.

\section{Conclusions}
\label{SEC:conclusion}

The Mumford-Shah functional has been widely used for image segmentation.
Its Ambrosio-Tortorelli Approximation has been known for its relative ease in implementation,
segmentation ability, and $\Gamma$-convergence to the Mumford-Shah functional as the regularization parameter
$\epsilon$ goes to zero.  The segmentation ability is based on the assumption that the input image
$g$ is discontinuous across the boundaries between different objects, and this discontinuity must be
maintained in the limit of $\epsilon \to 0$ during numerical computation to retain the $\Gamma$-convergence
and the segmentation ability for infinitesimal $\epsilon$ (e.g., see \cite{Bellettini1994}).
However, the maintenance of discontinuity in $g$ is often forgotten and $g$ is treated implicitly
as a continuous function in actual computation.  As a consequence, it has been observed that the segmentation
ability of the Ambrosio-Tortorelli functional varies significantly with different values of $\epsilon$
and the functional can even fail to $\Gamma$-converge to the original functional for some cases.
Moreover, there exist very few published numerical studies on the behavior of the functional as $\epsilon \to 0$.

We have presented in Section~\ref{SEC:analysis} an asymptotic analysis on the gradient flow equation
of the Ambrosio-Tortorelli functional as $\epsilon \to 0$ for continuous $g$. The analysis shows that
the functional can have different segmentation behavior for small but finite $\epsilon$ and eventually
loses its segmentation ability for infinitesimal $\epsilon$. This is consistent with the existing observations
in the literature and the numerical examples in one and two spatial dimensions presented
in Section~\ref{SEC:numerics}. Based on the analysis, we have proposed a selection strategy for
$\epsilon$ and a scaling procedure for $u$ and $g$ in Section~\ref{SEC:select}. Numerical results with real images
show that they lead to a good segmentation of the Ambrosio-Tortorelli functional.

Finally, we recall that the Ambrosio-Tortorelli functional is a special example of phase-field modeling
for image segmentation. We hope that the analysis and the selection strategy for the regularization parameter
presented in this work can also apply to other phase-field models. We are specially interested in the phase-field
modeling of brittle fracture (e.g., see \cite{Bourdin-2000,Francfort-1998,Miehe-2010}).
Investigations in this direction are currently underway.


\end{document}